\documentclass[12pt]{amsart}
\usepackage{epic,eepic}

\usepackage{mathptmx}

\newtheorem{lemma}{Lemma}[section]
\newtheorem{corollary}[lemma]{Corollary}
\newtheorem{theorem}[lemma]{Theorem}
\newtheorem{proposition}[lemma]{Proposition}

\newtheorem{definition}[lemma]{Definition}
\newtheorem{remark}[lemma]{Remark}
\newtheorem{example}[lemma]{Example}

\newtheorem{remark/definition}[lemma]{Remark/Definition}

\def\Box#1#2#3{\multiput(#1,#2)(1,0){2}{\line(0,1)1}
                            \multiput(#1,#2)(0,1){2}{\line(1,0)1}
                            \put(#1,#2){\makebox(1,1){$#3$}}}
\def\LBox#1#2#3#4{\multiput(#1,#2)(#4,0){2}{\line(0,1)1}
                            \multiput(#1,#2)(0,1){2}{\line(1,0){#4}}
                            \put(#1,#2){\makebox(#4,1){$#3$}}}
\def\Ci#1#2{\put(#1.5,#2.5){\circle{0.7}}}

\def\NZQ{\mathbb}               % the font for N,Z,Q,R,C
\def\NN{{\NZQ N }}

\def\ZZ{{\NZQ Z}}
\def\RR{{\NZQ R}}

\def\opn#1#2{\def#1{\operatorname{#2}}} % to make operators
\opn\chara{char}
\opn\gr{gr}
\opn\Rees{{\mathcal R}}
\opn\rank{rank}
\def\pot#1#2{#1[\kern-0.28ex[#2]\kern-0.28ex]}

\let\sect=\cap
\let\dirsum=\oplus

\let\iso=\cong

\let\Sect=\bigcap
\let\Dirsum=\bigoplus

\let\to=\rightarrow

\def\Implies{\ifmmode\Longrightarrow \else
      \unskip${}\Longrightarrow{}$\ignorespaces\fi}
\def\implies{\ifmmode\Rightarrow \else
      \unskip${}\Rightarrow{}$\ignorespaces\fi}

\let\:=\colon
\let\epsilon=\varepsilon
\let\phi=\varphi
\let\kappa=\varkappa

\opn\ini{in}
\opn\KRS{KRS}
\def\mm{{\mathfrak m}}
\def\pp{{\mathfrak p}}
\def\qq{{\mathfrak q}}
\opn\krs{krs}
\opn\diag{diag}
\opn\QF{QF}
\opn\DD{{\mathcal D}}
\opn\SS{{\mathcal S}}
\opn\MM{{\mathcal M}}
\opn\GL{GL}
\def\w{w}
\opn\height{height}
\opn\length{length}
\def\sep{\,|\,}
\opn\Min{Min}
\opn\Ins{Ins}

\opn\In{\textup{I}}
\opn\Ker{Ker}
\opn\CoKer{CoKer}
\opn\Tor{Tor}
\opn\Paths{Paths}
\opn\cl{cl}
\def\mod{\mathbin{{\textup{mod}}}}
\def\depth{{\operatorname{depth}}}

\def\alphah{{\widehat\alpha}}
\def\gammah{{\widehat\gamma}}

\def\XX{{\mathcal X}}

\def\addots{\mathinner{\mkern1mu\raise1pt\hbox{.}\mkern2mu\raise4pt\hbox{.}
         \mkern2mu\raise7pt\vbox{\kern7pt\hbox{.}}\mkern1mu}}

\let\oldlabel=\label
\def\prellabel{\marginparsep=1em\marginparwidth=80pt
     \def\label##1{\oldlabel{##1}\ifmmode\else\ifinner\else
          \marginpar{{\footnotesize\ \\ \tt
                     ##1}}\fi\fi}}
\begin{document}

\title{Gr\"obner bases and determinantal ideals}
\author{Winfried Bruns and Aldo Conca}
\address{Universit\"at Osnabr\"uck, FB
Mathematik/Informatik, 49069 Osna\-br\"uck, Germany}
\email{Winfried.Bruns@math.uos.de}
\address{Dipartimento di Matematica, Universit\'a di Genova,
Via Dodeca\-neso 35, 16146 Genova, Italy}
\email{conca@dima.unige.it}

\begin{abstract}
We give an introduction to the theory of determinantal ideals and
rings, their Gr\"obner bases, initial ideals and algebras,
respectively. The approach is based on the straightening law and
the Knuth-Robinson-Schensted correspondence. The article contains
a section treating the basic results about the passage to initial
ideals and algebras.
\end{abstract}


\maketitle

Let $K$ be a field and $X$ an $m\times n$ matrix of indeterminates
over $K$. For a given positive integer $t\leq \min(m,n)$, we
consider the ideal $I_t=I_t(X)$ generated by the $t$-minors
(i.~e.\ the determinants of the $t\times t$ submatrices) of $X$ in
the polynomial ring $K[X]$ generated by all the indeterminates
$X_{ij}$.

 From the viewpoint of algebraic geometry $K[X]$ should be regarded
as the coordinate ring of the variety of $K$-linear maps $f\:
K^m\to K^n$. Then $V(I_t)$ is just the variety of all $f$ such
that $\rank f<t$, and $K[X]/I_t$ is its coordinate ring.

The study of the determinantal ideals $I_t$ and the objects
related to them has numerous connections with invariant theory,
representation theory, and combinatorics. For a detailed account
we refer the reader to Bruns and Vetter \cite{BV}. A large part of
the theory of determinantal ideals can be developed over the ring
$\ZZ$ of integers (instead of a base field $K$) and then
transferred to arbitrary rings $B$ of coefficients (see
\cite{BV}). For simplicity we restrict ourselves to fields.

This article follows the line of investigation started by
Sturmfels' article \cite{Stu1} in which he applied the
Knuth-Robinson-Schensted correspondence KRS (Knuth \cite{Kn}) to
the study of the determinantal ideals $I_t$. The ``witchcraft''
(Knuth \cite[p.~60]{Kn2}) of the KRS saves one from tracing the
Buchberger algorithm through tedious inductions.

Later on the method was extended by Herzog and Trung \cite{HT} to
the so-called $1$-cogenerated ideals, ladder determinantal ideals
and ideals of pfaffians. They follow the important principle to
derive properties of $K[X]/I_t$ from the analogous properties of
$K[X]/\ini(I_t)$: the two rings appear as the generic and special
fiber of a flat $1$-parameter deformation. (By $\ini(I_t)$ we
denote the ideal of initial forms with respect to a suitable term
order.) The ring $K[X]/\ini(I_t)$ is the Stanley-Reisner ring of a
shellable simplicial complex and amenable to combinatorial methods
(see Stanley \cite{Sta2} and Bruns-Herzog \cite{BH}). In contrast
to the otherwise very elegant ASL approach, one does not replace
the indeterminates of $K[X]$ by a system of algebra generators
containing elements of degree $>1$. This is a major advantage if
one wants to investigate the Hilbert function and related
invariants.

The principle of deriving properties of ideals and algebras from
their initial counterparts was followed by the authors in
\cite{Co4}, \cite{BC1} and \cite{BC4} for the investigation of
algebras defined by minors, like the Rees algebra and the
subalgebra of $K[X]$ generated by the $t$-minors. This requires
the determination of Gr\"obner bases and initial ideals of powers
and products of determinantal ideals. On the KRS side the
necessary results are given by the theorem of Greene \cite{Gre}
and its variant found in \cite{BC1}.

In Section \ref{SectDet} we treat the straightening law of
Doubilet, Rota and Stein \cite{DRS} in the approach of De Concini,
Eisenbud, and Procesi \cite{DEP1}. It is an indispensable tool.
Moreover, we show that the residue class rings $K[X]/I_t$ are
normal domains. Section \ref{SectPow} contains the description of
the symbolic powers of the $I_t$ and the primary decomposition of
products $I_{t_1}\cdots I_{t_u}$ given in \cite{DEP1} and
\cite{BV}. While these two sections form the introduction to
determinantal ideals, Section \ref{SectGB} gives a fairly
self-contained treatment of initial ideals and algebras. Despite
the title of the article, the emphasis is on initial ideals and
not on Gr\"obner bases.

Section \ref{SectKRS} gives a short introduction to KRS (in the
``dual'' version of \cite{Kn}) and the theorems of Schensted
\cite{Sch} and Greene \cite{Gre}. These results are exploited in
Section \ref{SectGB-KRS} for determinantal ideals, their powers
and their products.

All the lines of development are brought together in Section
\ref{SectCM+HF} where we deal with the properties of $K[X]/I_t$,
especially with its Hilbert function, following Conca and Herzog
\cite{CH1}. At the end of this section we have inserted some
remarks that point out the extension to ladder determinantal
ideals and the variants for symmetric and alternating matrices.

The last two sections deal with algebras of minors, their
normality and Cohen-Macaulayness (Section \ref{SectAlgCM}) and
their canonical modules and Gorensteinness (Section
\ref{SectAlgCan}).

\section{Determinantal ideals and the straightening law}
\label{SectDet}

Almost all of the approaches one can choose for the investigation
of determinantal rings use standard bitableaux, to be defined
below, and the straightening law. In this approach one considers
all the minors of $X$ (and not just the $1$-minors $X_{ij}$) as
generators of the $K$-algebra $K[X]$ so that products of minors
appear as ``monomials''. The price to be paid, of course, is that
one has to choose a proper subset of all these ``monomials'' as a
linearly independent $K$-basis: we will see that the standard
bitableaux form a basis, and the straightening law tells us how to
express an arbitrary product of minors as a $K$-linear combination
of the basis elements. (In the literature standard bitableaux are
often called \emph{standard monomials}; however, we will have to
use the ordinary monomials in $K[X]$ so often that we reserve the
term ``monomial'' for products of the $X_{ij}$.)

Below we must often consider sequences of integers with a
monotonicity property. We say that a sequence $(r_i)$ is
\emph{increasing} if $r_i<r_{i+1}$ for all $i$. It is
\emph{non-increasing} if $r_i\ge r_{i+1}$ for all $i$.

Apart from Section \ref{SectCM+HF}, the letter $\Delta$ always
denotes a product $\delta_1\cdots\delta_w$ of minors, and we
assume that the sizes $|\delta_i|$ (i.~e.\ the number of rows of
the submatrix $X_i'$ of $X$ such that $\delta_i=\det(X_i')$) are
non-increasing, $|\delta_1|\ge \dots\ge |\delta_w|$. By
convention, the value of the empty minor $[\sep]$ is $1$. The
\emph{shape} $|\Delta|$ of $\Delta$ is the sequence
$(|\delta_1|,\dots,|\delta_w|)$. If necessary we may add factors
$[\sep]$ at the right hand side of the products, and accordingly
extend the shape by a sequence of $0$.

We denote the set of all non-empty minors of $X$ by $\MM(X)$ and
the subset of minors of size $t$ by $\MM_t(X)$. If no confusion
about the underlying matrix is possible, we will simply write
$\MM$ or $\MM_t$. Moreover,
$$
[a_1\dots a_t\sep b_1\dots b_t]
$$
stands for the determinant of the matrix $(X_{a_ib_j}\:
i=1,\dots,t,\ j=1,\dots,t)$. While we do not impose any condition
on the indices $a_i,b_i$ of $[a_1\dots a_t\sep b_1\dots b_t]$ in
general, we require that $a_1<\dots<a_t$ and $b_1<\dots<b_t$ if we
speak of a minor.

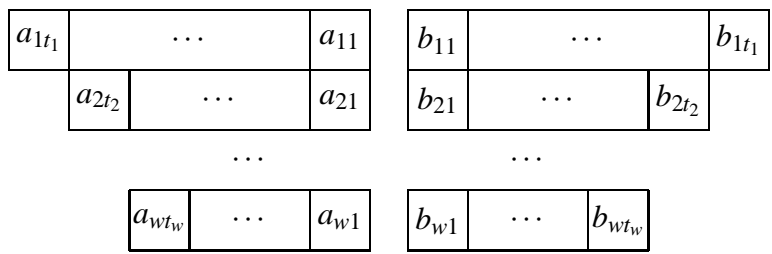
\begin{figure}[hbt]
\begin{gather*} \unitlength=0.8cm
\begin{picture}(6,4)(0,0)
\Box03{a_{1t_1}} \LBox13{\cdots}4 \Box53{a_{11}} \Box12{a_{2t_2}}
\LBox22{\cdots}3 \Box52{a_{21}} \put(2,1){\makebox(4,1){$\cdots$}}
\Box20{a_{wt_w}} \LBox30{\cdots}2 \Box50{a_{w1}}
\end{picture} \hspace{0.5cm}
\begin{picture}(6,4)(0,0)
\Box03{b_{11}} \LBox13{\cdots}4 \Box53{b_{1t_1}} \Box02{b_{21}}
\LBox12{\cdots}3 \Box42{b_{2t_2}}
\put(0,1){\makebox(4,1){$\cdots$}} \Box00{b_{w1}} \LBox10{\cdots}2
\Box30{b_{wt_w}}
\end{picture}
\end{gather*}
\caption{A bitableau}\label{Young}
\end{figure}

A product of minors is also called a \emph{bitableau}. The choice
of the term bitableau is motivated by the graphical description of
a product $\Delta$ as a pair of Young tableaux as in Figure
\ref{Young}.
Every product of minors is represented by a bitableau and,
conversely, every bitableau stands for a product of minors:
$$
\Delta=\delta_1\cdots\delta_w,\qquad \delta_i=[a_{i1}\dots
a_{it_i}\sep b_{i1}\dots b_{it_i}],\ i=1,\dots,w.
$$
According to our convention above, the indices in each row of the
bitableau are increasing from the middle to both ends. Sometimes
it is necessary to separate the two tableaux from which $\Delta$
is formed; we then write $\Delta=(R\sep C)$.

For formal correctness one should consider the bitableaux as
purely combinatorial objects (as we will do in Section
\ref{SectKRS}) and distinguish them from the ring-theoretic
objects represented by them, but since there is no real danger of
confusion, we use the same terminology for both classes of
objects.

Whether $\Delta$ is a standard bitableau is controlled by a
partial order on $\MM(X)$, namely
\begin{multline*}
[a_1\dots a_t\sep b_1\dots b_t] \preceq [c_1\dots c_u\sep d_1\dots d_u]\\
\iff\quad t\ge u\quad\text{and}\quad a_i\le c_i,\quad b_i\le d_i,\
\quad i=1,\dots,u.
\end{multline*}
A product $\Delta=\delta_1\cdots\delta_w$ is called a
\emph{standard bitableau} if
$$
\delta_1\preceq\dots\preceq\delta_w,
$$
in other words, if in each column of the bitableau the indices are
non-decreasing from top to bottom. The letter $\Sigma$ is reserved
for standard bitableaux. (The empty product is also standard.)

The fundamental straightening law of Doubilet--Rota--Stein
\cite{DRS} says that every element of $R$ has a unique
presentation as a $K$-linear combination of standard bitableaux:

\begin{theorem}\label{straight}\nopagebreak
\begin{itemize}
\item[(a)] The standard bitableaux are a $K$-vector space basis of
$K[X]$.

\item[(b)] If the product $\gamma\delta$ of minors is
not a standard bitableau, then it has a representation
$$
\gamma\delta=\sum x_i\epsilon_i\eta_i,\qquad x_i\in K,\ x_i\neq0,
$$
where $\epsilon_i\eta_i$ is a standard bitableau,
$\epsilon_i\prec\gamma,\delta\prec\eta_i$ (here we must allow that
$\eta_i=[\sep]=1$).

\item[(c)] The standard representation of an
arbitrary bitableau $\Delta$, i.e.\ its representation as a linear
combination of standard bitableaux $\Sigma$, can be found by
successive application of the straightening relations in
\emph{(b)}.

\end{itemize}
\end{theorem}

For the proof of the theorem we can assume that $m\le n$, passing
to the transpose of $X$ if necessary. We derive the theorem from
its ``restriction'' to the subalgebra $K[\MM_m]$ generated by the
$m$-minors. Each $m$-minor is determined by its column indices,
and for simplicity we set
$$
[b_1\dots b_m]=[1\dots m\sep b_1\dots b_m].
$$
The algebra $K[\MM_m]$ is the homogeneous coordinate ring of the
Grassmann variety of the $m$-dimensional vector subspaces of
$K^n$. The $m$-minors satisfy the famous \emph{Pl\"ucker
relations}. In their description we use $\sigma(i_1\dots i_s)$ to
denote the sign of the permutation of $\{1,\dots,s\}$ represented by the sequence $i_1,\dots,i_s$.

\begin{lemma}\label{Det:Pl}
For all indices $a_1,\dots,a_p$, $b_q,\dots,b_m$,
$c_1,\dots,c_s\in \{1,\dots,n\}$ such that $s=m-p+q-1>m$ and
$t=m-p>0$ one has
$$
\sum_{\scriptstyle i_1<\cdots <i_t \atop{\scriptstyle
i_{t+1}<\cdots <i_s \atop\scriptstyle \{i_1,\dots,i_s\}=
\{1,\dots,s\}}} \sigma (i_1 \dots i_s)[a_1 \dots a_p c_{i_1} \dots
c_{i_t}] [c_{i_{t+1}} \dots c_{i_s} b_q \dots b_m]=0.
$$

\end{lemma}
\begin{proof}
Let $V$ be the $K$-vector space generated by the columns $X_j$ of $X$. We define $\alpha\:
V^s \to K(X)$ by
\begin{multline*}
\alpha (y_1,\dots,y_s)= \sum \sigma (i_1 \dots i_s) \det(X_{a_1},
\dots, X_{a_p},
y_{i_1}, \dots, y_{i_t})\cdot \\
\cdot\det(y_{i_{t+1}}, \dots, y_{i_s}, X_{b_q}, \dots, X_{b_m}).
\end{multline*}
where the sum has the same range as above. It is straightforward
to check that $\alpha$ is an alternating multilinear form on
$V^s$. Since $s>\dim V=m$, one has $\alpha=0$. \end{proof}

Let $[a_1 \dots a_m]$, $[b_1 \dots b_m]$ be elements of $\MM_m(X)$
such that $a_i \le b_i$ for $i=1,\dots,p$, but $a_{p+1}>b_{p+1}$.
We put
$$
q=p+2,\quad s=m+1,\quad (c_1,\dots,c_s)=(a_{p+1},\dots,a_m,
b_1,\dots,b_{p+1}).
$$
Then, in the Pl\"ucker relation with these data, all the non-zero
terms
$$
[d_1\dots d_m][e_1\dots e_m] \neq [a_1\dots a_m][b_1\dots
b_m]
$$
have the following properties (after their column indices
have been arranged in ascending order):
$$
[d_1\dots d_m]\preceq [a_1\dots a_m],\qquad d_1 \le
e_1,\dots,d_{p+1}\le e_{p+1}.
$$
By induction on $p$ it follows that a product $\gamma\delta$ of
maximal minors $\gamma$ and $\delta$ is a linear combination of
standard bitableaux $\alpha\beta$, $\alpha\preceq\beta$, such that
$\alpha\preceq\gamma$. Note that $\alpha$ and $\beta$ arise from
$\gamma$ and $\delta$ by an exchange of indices.

Let $\gamma_1\cdots\gamma_u$ be a product of maximal minors of
length $u>2$, If it is not a standard bitableau, then we find an
index $i$ with $\gamma_i\not\preceq\gamma_{i+1}$. As just seen,
$\gamma_i\gamma_{i+1}$ can be expressed as a linear combination of
standard bitableaux. Substitution of this expression for
$\gamma_i\gamma_{i+1}$ yields a representation of
$\gamma_1\cdots\gamma_u$ as a linear combination of bitableaux. In
each of these bitableaux, indices from the $(i+1)$st row of the
bitableau of $\gamma_1\cdots\gamma_u$ have been exchanged with
\emph{larger} indices from its $i$th row. An iteration of this
procedure must eventually yield a linear combination of standard
bitableaux, since the exchange of indices can only be repeated
finitely many times.

This proves that products of maximal minors are linear
combinations of standard bitableaux. In the next step we transfer
this partial result to the full set of minors. We extend $X$ by
$m$ columns of further indeterminates, obtaining
$$
X'=\begin{pmatrix} X_{11}&\cdots
&X_{1n}&X_{1,n+1}&\cdots&X_{1,n+m}\\
\vdots&&\vdots&\vdots&&\vdots \\
X_{m1}&\cdots &X_{mn}&X_{m,n+1}&\cdots&X_{m,n+m}
\end{pmatrix}.
$$
Then $K[X']$ is mapped onto $K[X]$ by substituting for each entry
of $X'$ the corresponding entry of the matrix
$$
\begin{pmatrix}
X_{11}&\cdots&X_{1n}&0&\cdots&\cdots&0&1 \\
&&&\vdots&&\addots&\addots&0 \\
\vdots&&\vdots&\vdots &\addots&\addots &\addots&\vdots \\
&&&0&\addots&\addots&&\vdots\\
X_{m1}&\cdots&X_{mn}&1&0&\cdots&\cdots&0
\end{pmatrix}.
$$
Let $\phi:K[\MM(X')]\to K[X]$ be the induced $K$-algebra
homomorphism. Then
$$
\phi([b_1\dots b_m])=\pm [a_1\dots a_t\sep b_1\dots,b_t]
$$
where $t=\max \{i\: b_i \le n\}$ and $a_1,\dots,a_t$ have been
chosen such that
$$
\{a_1,\dots, a_t,n+m+1-b_m,\dots,n+m+1-b_{t+1}\}=\{1,\dots,m\}.
$$
Evidently $\phi$ is surjective, and furthermore it sets up a
bijective correspondence between the set $\MM_m(X')$ of $m$-minors
of $X'$ and $\MM(X)\cup\{\pm 1\}$; on
$\MM_m(X')\setminus\{\tilde\mu\}$ the correspondence is an
isomorphism of partially ordered sets. Note that the maximal
element $\tilde\mu=[n+1\dots n+m]$ of $\MM_m(X')$ is mapped to
$\pm 1$ by $\phi$, and, up to sign, standard monomials go to
standard monomials. (We leave the verification of this fact to the
reader; the details can also be found in \cite{BV}, (4.9).)

In order to represent an arbitrary element of $K[X]$ as a linear
combination of bitableaux, we lift it to $K[\MM_m(X')]$ via
$\phi$. Then the preimage is ``straightened'', and an application
of $\phi$ yields the desired expression in $K[X]$.

For part (a) of Theorem \ref{straight} it remains to show the
linear independence of the standard bitableaux. We know already
that they generate the vector space $K[X]$. Moreover, they are
homogeneous elements with respect to total degree, and as we will
see in Section \ref{SectKRS}, there are as many standard
bitableaux in every degree as there are ordinary monomials. This
implies the linear independence of the standard bitableaux and
finishes the proof of Theorem \ref{straight}(a). Therefore we may
now speak of the straightening law in $K[X]$. As we have seen,
arbitrary products of minors can be straightened by the successive
straightening of products with two factors, and so part (c) has
also been proved.

For part (b) we notice that the Pl\"ucker relations are
homogeneous of degree $2$. Therefore there are exactly two factors
$\epsilon_i$ and $\eta_i$ in each term on the right hand side of
$$
\gamma\delta=\sum x_i\epsilon_i\eta_i
$$
if $|\gamma|,|\delta|=m$. Since the straightening law in $K[X]$ is
a specialization via $\phi$ of that in $K[\MM_m(X')]$, there can be
at most $2$ factors in each of the summands on the right hand
side.

It follows easily from the straightening procedure that
$\epsilon_i\preceq\gamma$. In fact, the inequality holds for all
intermediate expressions that arise in the successive application
of the Pl\"ucker relations, as observed above. In order to see
that $\epsilon_i\preceq\delta$ as well, we simply straighten
$\delta\gamma=\gamma\delta$. By the linear independence of the
standard bitableaux, the result is the same. Also for
$\gamma,\delta\preceq\eta_i$ one has to use the linear
independence of the standard monomials (the intermediate
expressions in the straightening procedure may violate it). It is
enough to prove this relation in $K[\MM_m(X')]$. On the set
$\MM_m(X')$ we consider the reverse partial order, arising from
rearranging the columns of the matrix in the order
$m+n,m+n-1,\dots,1$. It has the same set of standard monomials as
$\preceq$, at least up to sign. By linear independence,
straightening with respect to the reverse partial order must have
the same result as that with respect to $\preceq$ (up to sign).
This concludes the proof of Theorem \ref{straight}(b).

\begin{corollary}\label{dimGrass}\nopagebreak
\begin{itemize}
\item[(a)] The kernel of $\phi:K[\MM_m(X')]\to K[X]$ is generated by
${\tilde\mu+1}$ or $\tilde\mu-1$.
\item[(b)] $\dim K[\MM_m(X)]=m(n-m)+1$.
\end{itemize}
\end{corollary}

\begin{proof}
(a) Every element $x$ of $K[\MM_m(X')]$ has a unique
representation $x=\sum p_\Delta(\tilde\mu)\Delta$ as a linear
combination of standard bitableaux $\Delta$ over
$\MM_m(X')\setminus\{\tilde\mu\}$ with coefficients
$p_\Delta(\tilde\mu)\in K[\tilde\mu]$. Clearly $\phi(x)=0$ if and
only if $\phi(p_\Delta(\tilde\mu))=p_\Delta(\pm1)=0$ for all
$\Delta$.

(b) $X$ plays the role of $X'$ for an $m\times (n-m)$ matrix of
indeterminates. Now one applies (a). \end{proof}

\begin{remark}\rm
The straightening law for $K[\MM_m(X)]$ is due to Hodge
\cite{Hodge}, \cite{HP}. The proof of the straightening law for
$K[X]$, with the exception of the linear independence of the
standard bitableaux, follows De Concini, Eisenbud, and Procesi
\cite{DEP1}. It has also been reproduced in \cite{BV}. The linear
independence can be proved without the KRS; see \cite{DEP1} (or
\cite{BV}) and \cite{HP} (or \cite{BH}) for two alternative
proofs.

A third alternative: one proves the linear independence of the
standard bitableaux in $\MM_m(X)$ by a Gr\"obner basis argument
(see Remark \ref{MaxLI}), shows \ref{dimGrass}(b), deduces
\ref{dimGrass}(a), and concludes the linear independence of all
standard bitableaux.
\end{remark}

The straightening law can be refined. Let $e_1,\dots,e_m$ and
$f_1,\dots,f_n$ denote the canonical $\ZZ$-bases of $\ZZ^m$ and
$\ZZ^n$ respectively. Clearly $K[X]$ is a
$\ZZ^m\dirsum\ZZ^n$-graded algebra if we give $X_{ij}$ the
``vector bidegree'' $e_i\dirsum f_j$. All minors and bitableaux
are homogeneous with respect to this grading. The coordinates of
the vector bidegree of a bitableau $\Delta$, usually called the
\emph{content} of $\Delta$, just count the multiplicities with
which the rows and columns of $X$ appear. The homogeneity of
bitableaux implies that \emph{straightening preserves content}.

Next we can compare the shapes of the tableaux appearing on the
left and the right hand side of a straightening relation. For a
sequence $\sigma=(s_1,\dots,s_u)$ we set
$$
\alpha_k(\sigma)=\sum_{i\le k} s_i,
$$
and define $\sigma\le \tau$ by $\alpha_k(\sigma)\le\alpha_k(\tau)$
for all $k$. It follows from Theorem \ref{straight}(b) and (c)
that \emph{straightening does not decrease shape}.

It is a natural question whether at least one standard bitableau
in a straightening relation must have the same shape as the left
hand side. This is indeed true, and for a more precise statement
we introduce the \emph{initial tableau}  $\In(\sigma)$ of shape
$\sigma=(s_1,\dots,s_u)$: in its $k$th row it contains the numbers
$1,\dots,s_k$.

\begin{theorem}\label{strongstr}
Let $\Delta$ be a bitableau of shape $\sigma$, with row tableau
$R$ and column tableau $C$.
\begin{itemize}

\item[(a)] Every bitableau in the standard representation of
$\Delta$ has the same content as $\Delta$.

\item[(b)] $(R\sep \In(\sigma))$ has a standard
representation $\sum \alpha_i(R_i\sep \In(\sigma))$, $\alpha_i\neq
0$, and $(\In(\sigma)\sep C)$ has a standard representation $\sum
\beta_j(\In(\sigma)\sep C_j)$, ${\beta_j\neq 0}$.

\item[(c)] $\Delta-\sum \alpha_i\beta_j (R_i\sep C_j)$ is a linear
combination of standard bitableaux of size $>\sigma$.
\end{itemize}
\end{theorem}

Part (a) has been justified above. Part (b) follows from the fact
that there is no tableau of shape $\ge\sigma$ that has the same
content as $\In(\sigma)$. It is much more difficult to show (c),
and we forego a proof; for example, see \cite[(11.4)]{BV}.

An ideal in a partially ordered set $(M,\le)$ is a subset $N$ such
that $N$ contains all elements $x\le y$ if $y\in N$. Let
${\mathcal N}$ be an ideal in the partially ordered set $\MM(X)$
and consider the ideal $I={\mathcal N}K[X]$ generated by
${\mathcal N}$. Every element of $I$ is a $K$-linear combination
of elements $\delta x$ with $x\in K[X]$ and $\delta\in {\mathcal
N}$. It follows from Theorem \ref{straight} that every standard
bitableaux $\Sigma=\gamma_1\cdots \gamma_v$ in the standard
representation of $\delta x$ has $\gamma_1\le \delta$, and
therefore $\gamma_1\in {\mathcal N}$. This shows

\begin{proposition}\label{po-id}
Let ${\mathcal N}$ be an ideal in the partially ordered set
$\MM(X)$. Then the standard bitableaux $\Sigma=\gamma_1\cdots
\gamma_u$ with $\gamma_1\in {\mathcal N}$ form a $K$-basis of the
ideal $I={\mathcal N}K[X]$ in the ring $K[X]$, and the (images of
the) the standard bitableaux $\Sigma'=\delta_1\cdots \delta_v$
with $\delta_j\notin {\mathcal N}$ for all $j$ form a $K$-basis of
$K[X]/I$.
\end{proposition}

\begin{corollary}\label{Itstan}
The standard bitableaux $\Sigma=\gamma_1\cdots\gamma_u$ such that
$|\gamma_1|\ge t$ form a $K$-basis of $I_t$, and (the images of)
the standard bitableaux $\Sigma'=\delta_1\cdots \delta_v$ with
$|\delta_j|\le t-1$ for all $j$ form a $K$-basis of $K[X]/I_t$.
\end{corollary}

In fact, $I_t$ is generated by all minors of size $\ge t$, and
these form the ideal
$$
\{\delta\in\MM(X): \delta\not\succeq [1\dots t-1\sep 1\dots t-1]\}.
$$
The ideals $I_t$ are special instances of the so-called
\emph{$1$-cogenerated ideals} $I_\gamma$, $\gamma\in\MM(X)$, that
are generated by all minors $\delta\not\succeq \gamma$.

\begin{remark}\rm
The straightening law and its refined version in Theorem
\ref{strongstr} can be used in various approaches to the theory of
the determinantal ideals and rings:

(a) The straightening law implies that $K[X]$ and $K[\MM_m]$ are
\emph{algebras with straightening law} on the partially ordered
sets $\MM(X)$ and $\MM_m(X)$, resp. This property is passed on the
residue class rings modulo the class of ideals considered in
Proposition \ref{po-id}. See \cite{DEP2} and \cite{BV}.

(b) The ``filtration by shapes'' as indicated in Theorem
\ref{strongstr} can be used for a deformation process. This allows
one to deduce properties of the determinantal ring $K[X)/I_t$ from
the ``semigroup of shapes'' occurring in it. See \cite{BC2}.

(c) (Related to (b).) The refined form of the straightening law is
the basis for the investigation of the determinantal rings via
representation theory. See Akin, Buchsbaum and Weyman \cite{ABW},
\cite{DEP1} or \cite{BV}.

\end{remark}

The straightening law allows us to prove basic properties about
the ideals $I_t$ without much effort. Additionally we need the
following extremely useful \emph{induction lemma}. It uses that an
ideal of type $I_t(A)$ remains unchanged if one applies elementary
row and column operations to the matrix $A$.

\begin{lemma}\label{local}
Let $X=(X_ {ij})$ and $Y=(Y_ {ij})$ be matrices of indeterminates
over $K$ of sizes $m \times n$ and $(m-1) \times (n-1)$, resp.
Then the substitution
\begin{align*}
X_ {ij} &\to Y_{ij} + X_ {mj}X_ {in}X_ {mn}^{-
1},& 1 &\le i \le m-1,&1 &\le j \le n-1, \\
X_ {mj} &\to X_ {mj},X_ {in}\to X_ {in},&1 &\le i \le m,& 1 &\le j
\le n,
\end{align*}
induces an isomorphism
\begin{equation}
K[X][X_{mn}^{-1}]\iso K[Y][X_ {m1},\dots,X_ {mn},X_ {1n},\dots,X_
{m-1,n}][X_ {mn}^{-1}]\tag{$*$}
\end{equation}
under which the extension of $I_t(X)$ is mapped to the extension
of $I_{t-1}(Y)$. Therefore one has an isomorphism
$$
(K[X]/I_t(X))[x_{mn}^{-1}]\iso (K[Y]/I_{t-1}(Y))[X_ {m1},\dots,X_
{mn},X_ {1n},\dots,X_ {m-1,n}][X_ {mn}^{-1}].
$$
(here $x_ {mn}$ denotes the residue class of $X_ {mn}$ in
$K[X]/I_t(X)$).
\end{lemma}

\begin{proof}
The substitution has an inverse, namely $Y_ {ij} \to
X_{ij}-X_{mj}X_{in}X_{mn}^{-1}$, $X_{mj}\to X_{mj}$, $X_{in}\to
X_{in}$, and so induces the isomorphism $(*)$.

Since $X_{mn}$ is invertible in $K[X][X_{mn}^{-1}]$, one can apply
elementary row and column operations to the matrix $X$ with pivot
element $X_{mn}$. Now one just ``writes'' $Y_{ij}$ for the entries
in the rows $1,\dots,m-1$ and the columns $1,\dots,n-1$ of the
transformed matrix. With this identification one has
$I_t(X)=X_{mn}I_{t-1}(Y)=I_{t-1}(Y)$ in the ring
$K[X][X_{mn}^{-1}]$. \end{proof}

\begin{theorem}\label{normal}
The ring $K[X]/I_t(X)$ is a normal domain of dimension
$(m+n-t+1)(t-1)$. Its singular locus is the variety of the ideal
$I_{t-1}/I_t$.
\end{theorem}

\begin{proof}
The case $t=1$ is trivial. Suppose that $t>1$. Set $R=K[X]/I_t$.
We claim that the residue class $x_{mn}$ of $X_{mn}$ is a
non-zero-divisor on $R$. In fact, the product of $X_{mn}=[m\sep
n]$ and a standard bitableau is a standard bitableau, and if $x$
is a linear combination of standard bitableaux without a factor of
size $\ge t$, then so is $X_{mn}x$. Corollary \ref{Itstan} now
implies our claim. The argument shows even more: $X_{mn}$ is a
non-zero-divisor modulo every ideal of $K[X]$ generated by an
ideal $\mathcal N$ in the partially ordered set $\MM(X)$ such that
$X_{mn}\notin \mathcal N$.

In order to verify that $R$ is a domain, it suffices to prove this
property for $R[x_{mn}^{-1}]$, and to the latter ring we can apply
induction via Proposition \ref{local}.

The dimension formula follows by the same induction, since an
affine domain does not change dimension upon the inversion of a
non-zero element.

Let $\mm$ be the maximal ideal of $R$ generated by the residue
classes $ x_{ij}$ of the indeterminates. Clearly $R_\mm$ is not
regular if $t\ge 2$. Every other prime ideal $\pp$ does not
contain one of the $x_{ij}$, and by symmetry we can assume
$x_{mn}\notin\pp$. Then $R_\pp$ is of the form
$S[X_{in},X_{mj},X_{mn}^{-1}]_\qq$ with $S=K[Y]/I_{t-1}(Y)$. It
follows that $R_\pp$ is regular if and only if $S_{\qq\cap S}$ is
regular. Moreover, $\pp$ contains $I_{t-1}(X)/I_t(X)$ if and only
if $\qq\cap S$ contains $I_{t-2}(Y)/I_{t-1}(Y)$. Again we can
apply induction to prove the claim about the singular locus.

For normality we use Serre's normality criterion. We know the
singular locus, and by the dimension formula it has codimension
$\ge 2$. Now it is enough to show that $\depth R_\pp\ge 2$ if
$\dim R_\pp\ge 2$. If $\pp\neq\mm$, we obtain this by induction as
above, and it remains to show that $\depth R_\mm\ge 2$. The set of
minors $\delta$ of size $<t$ has a smallest element with respect
to $\preceq$, namely $\epsilon=[1\dots t-1\sep 1\dots t-1]$. The
same argument that we have applied to $X_{mn}$ shows that
$\epsilon$ is a non-zero-divisor modulo $I_t$. Moreover, the ideal
$J=I_t+(\epsilon)$ is generated by an ideal in $\MM(X)$, and so
$X_{mn}$ is a non-zero-divisor modulo $J$. It follows that $\mm$
contains a regular $R$-sequence of length $2$. \end{proof}

We will show in Section \ref{SectCM+HF} that the rings $K[X]/I_t$
are Cohen--Macaulay.

\section{Powers and products of determinantal ideals}
\label{SectPow}

In this section we want to determine the primary decomposition of
powers and, more generally, products $J=I_{t_1}\cdots I_{t_u}$,
$t_1\ge\dots \ge t_u$, of determinantal ideals. It is easy to see
that only the ideals $I_t$ with $t\le t_1$ can be associated to
$J$. In fact, suppose $\pp$ is a prime ideal in $R=K[X]/J$
different from the irrelevant maximal ideal $\mm$. Then $\pp$ does
not contain one of the $x_{ij}$, and so we may invert $x_{ij}$
without loosing the extension of $\pp$ as an associated prime
ideal of $(R/J)[x_{ij}^{-1}]$. But now the Induction Lemma
\ref{local} applies: by symmetry we can assume $(i,j)=(m,n)$.

Thus we have to find primary components of $J$ with respect to the
ideals $I_t$. Immediate, and as we will see, optimal candidates
are the symbolic powers of the ideals $I_t$. We determine them
first.

The ideal $\pp=I_t$ is a prime ideal in the regular ring $A=K[X]$.
With each such prime ideal one associates a valuation on the
quotient field of $A$ as follows. We pass to the localization
$P=A_\pp$ and let $\qq=\pp P$. Now we set
$$
v_\pp(x)=\max\{i: x\in \qq^i\}
$$
for all $x\in P$, $x\neq 0$, and $v_\pp(0)=\infty$. The associated
graded ring $\Dirsum_{i=0}^\infty \qq^i/\qq^{i+1}$ is a polynomial
ring over the field $P/\qq$ (we only use that it is an integral
domain). This implies $v_\pp(xy)=v_\pp(x)+v_\pp(y)$ for all
$x,y\in P$, and that $v_\pp(x+y)\ge \min(v_\pp(x),v_\pp(y))$ is
clear anyway. To sum up: $v_\pp$ is a discrete valuation on $P$
and can be extended to the quotient field $\QF(P)=\QF(A)$. By
definition, the $i$th symbolic power is $\pp^{(i)}=\pp^i P\cap A$.
With the help of $v_\pp$ we can also describe it as
$\pp^{(i)}=\{x\in A: v_\pp(x)\ge i\}$. If $\pp=\mm$ is the maximal
irrelevant ideal of $K[X]$ and $f$ is a homogeneous polynomial,
then $v_\pp(f)$ is the ordinary total degree of $f$.

For the choice $A=K[X]$, $\pp=I_t$, we denote $v_\pp$ by
$\gamma_t$. We claim that
$$
\gamma_t(\delta)=\begin{cases} 0,&            |\delta|<t,\\
                                |\delta|-t+1,& |\delta|\ge t.
                  \end{cases}
$$
For $t=1$, this follows immediately, since $\gamma_1(\delta)$ is
its total degree $|\delta|$. Let $t>1$. If $|\delta|<t$, then
$\delta\notin I_t$, and so $\gamma_t(\delta)=0$. Suppose that
$|\delta|\ge t$. It is useful to note that $\gamma_t(\delta)$ does
only depend on $|\delta|$: minors $\delta$ and $\delta'$ of the
same size are conjugate under an isomorphism of $K[X]$ that leaves
$\pp$ invariant. We can therefore assume that $\delta=[m-t+1\dots
m\sep n-t+1\dots n]$. The substitution in the induction lemma
\ref{local} maps $\delta$ to a minor of size $|\delta|-1$, and
reduces $t$ by $1$, too. So induction finishes the proof.

We transfer the function $\gamma_t$ to sequences of integers:
$$
\gamma_t(s_1,\dots,s_u)=\sum_{i=1}^u \max(s_i-t+1,0).
$$
For instance, if $\lambda=(4,3,3,1)$, then $\gamma_4 (\lambda)=1$,
$\gamma_3 (\lambda)=4$, $\gamma_2 (\lambda)=7$, $\gamma_1
(\lambda)=11$.

It is clear that $\gamma_t(s_1,\dots,s_u)=
\gamma_t(\delta_1\cdots\delta_u)$ for minors $\delta_i$ with
$|\delta_i|=s_i$, $i=1\dots,u$. In Section \ref{SectDet} we have
introduced the partial order of shapes based on the functions
$\alpha_k$. We can also use the $\gamma_t$ for such a comparison,
but it yields the same partial order:

\begin{lemma}\label{al-ga}
Let $\rho=(r_1,\dots,r_u)$ and $\sigma=(s_1,\dots,s_v)$ be
non-increasing sequences of integers. Then $\rho\le\sigma$ if and
only if $\gamma_t(\rho)\le\gamma_t(\sigma)$ for all $t$.
\end{lemma}

\begin{proof}
We use induction on $u$. The case $u=1$ is trivial. Furthermore,
there is nothing to show if $r_i\le s_i$ for all $i$. It remains
the case in which $r_j>s_j$ for some $j$. Next, if $r_i=s_i>0$ for
some $i$, then we can remove $r_i$ and $s_i$ and compare the
shortened sequences by induction. Thus we may assume that $r_i\neq
s_i$ for all $i$ with $r_i>0$. Let $k$ be the smallest index with
$s_k<r_k$. It is easy to see that $k\ge 2$ if $\rho\le\sigma$ or
$\gamma_{s_1}(\rho)\le\gamma_{s_1}(\sigma)$. Extending $\sigma$ by
$0$ if necessary we may assume that $k\le v$.

We have $s_{k-1}>r_{k-1}\ge r_k>s_k$; in particular, $s_{k-1}\ge
s_k+2$. Define $s_i'=s_i$ for $i\neq k-1,k$, $s_{k-1}'=s_{k-1}-1$,
$s_k'=s_k+1$, and $\sigma'=(s_1',\dots,s_v')$.

Suppose first that $\rho\le \sigma$. Then it follows easily that
$\rho\le\sigma'$. Moreover, $\gamma_t(\sigma')\le\gamma_t(\sigma)$
for all $t$, and a second inductive argument allows us to assume
that $\gamma_t(\rho)\le\gamma_t(\sigma')$ for all $t$.

Conversely, let $\gamma_t(\rho)\le\gamma_t(\sigma)$ for all $t$.
Since $\sigma'\le\sigma$, it is enough to show that
$\gamma_t(\rho)\le \gamma_t(\sigma')$ for all $t$. One has
$\gamma_t(\sigma')=\gamma_t(\sigma)$ for $t\le s_k+1$, and
$\gamma_t(\sigma')=\gamma_t(\sigma)-1$ for
$t=s_k+2,\dots,s_{k-1}$. Since obviously
$\gamma_t(\rho)<\gamma_t(\sigma)$ for $t>r_k$, the critical range
is $s_k+2\le t\le r_k$. Suppose that
$\gamma_t(\rho)=\gamma_t(\sigma)$ for some $t$ in this range.

Exactly $s_1,\dots,s_{k-1}$ contribute to $\gamma_{t-1}(\sigma)$,
but at least $r_1,\dots,r_k$ contribute to $\gamma_{t-1}(\rho)$.
It follows that
$$
k-1=\gamma_{t-1}(\sigma)-\gamma_t(\sigma)\le
\gamma_{t-1}(\rho)-\gamma_t(\rho)\ge k,
$$
a contradiction finishing the proof. \end{proof}

Setting
$$
I^{\sigma}=I_{s_1}\cdots I_{s_u}.
$$
for $\sigma=(s_1,\dots,s_u)$ we now describe the symbolic powers
of the $I_t$.

\begin{proposition}\label{symb}
One has
$$
I_t^{(k)}=\sum_{\scriptstyle \sigma=(s_1,\dots,s_u)
\atop{\gamma_t(\sigma)\ge k}} I^{\sigma}.
$$
$I_t^{(k)}$ has a $K$-basis of the standard bitableaux $\Sigma$
with $\gamma_t(\Sigma)\ge k$.
\end{proposition}

\begin{proof}
If $t=1$, the right hand side is just the ideal of elements of
total degree $\ge k$, and there is nothing to prove. Suppose that
$t>1$. Then $X_{mn}$ is a non-zero-divisor modulo $I_t^{(k)}$ (by
definition of the symbolic power). But it is also a
non-zero-divisor modulo the right hand side, as we will see.
Therefore we can invert $X_{mn}$. This transforms all sizes in the
right way, and equality follows by induction on $t$.

We have noticed in Theorem \ref{strongstr} that straightening does
not decrease shape. This implies $\gamma_t(\Sigma)\ge
\gamma_t(\Delta)$ for all standard bitableaux in the standard
representation of a bitableau $\Delta$. Thus the right hand side
has a $K$-basis by all $\Sigma$ with $\gamma_t(\Sigma)\ge k$.
Since multiplication by $X_{mn}$ does not affect values under
$\gamma_t$, and maps standard bitableaux to standard bitableaux,
the result follows. \end{proof}

Of course, only finitely many summands are needed for $I_t^{(k)}$.
The simplest non-trivial case is $t=k=2$:
$$
I_2^{(2)}=I_2^2+I_3,
$$
since every other summand is contained in $I_2^2$ or $I_3$. If
$m\le 2$ or $n\le 2$, then $I_2^2=I_2^{(2)}$. This observation is
easily generalized to the following result of Trung \cite{Tr}.

\begin{corollary}\label{maxi}
The symbolic powers of the ideal of maximal minors coincide with
the ordinary ones.
\end{corollary}

\begin{proof}
Suppose that $m=\min(m,n)$. Then a bitableau $\Delta$ has
$\gamma_m(\Delta)\ge k$ if and only if the first $k$ factors have
size exactly $m$. \end{proof}

The primary decomposition of products of the ideals $I_t$ depends
on characteristic. This indicates that the straightening law alone
is not sufficient to prove it. Actually the straightening law
enters the proof only via Theorem \ref{symb}.

\begin{theorem}\label{decomp}
Let $\rho=(r_1,\dots,r_u)$ be a non-increasing sequence of
integers and suppose that $\chara K=0$ or  $\chara
K>\min(r_i,m-r_i,n-r_i)$ for $i=1,\dots,u$. Then
$$
I^\rho=\bigcap_{t=1}^{r_1} I_t^{(\gamma_t(\rho))}.
$$
\end{theorem}

The theorem was proved by De Concini, Eisenbud and Procesi
\cite{DEP1} in characteristic $0$ and generalized in \cite{BV}.
The inclusion $\subseteq$ is a triviality (and independent of the
hypothesis on characteristic): $\gamma_t(x)\ge \gamma_t(\rho)$ for
all $x\in I^\rho$. Before we indicate the proof of the converse
inclusion, let us have a look at the first non-trivial case,
namely
$$
I_2^2=I_1^4\cap (I_3+I_2^2).
$$
For equality we must prove that the degree $\ge 4$ elements in
$I_3$, namely those of $I_3I_1$, are contained in $I_2^2$. This
type of containment is the crucial point in the proof of the
theorem. It is the ideal-theoretic analogue of the passage from
$\sigma$ to $\sigma'$ in the proof of Lemma \ref{al-ga}.

\begin{lemma}\label{destr}
Let $u,v$ be integers, $0\le u\le v-2$. Suppose that $\chara K=0$
or $\chara K>\min(u+1,m-(u+1),n-(u+1))$. Then $I_uI_v\subseteq
I_{u+1}I_{v-1}$.
\end{lemma}

For reasons of space we refer the reader to \cite[(10.10)]{BV} for
a proof of the lemma. It is based on a symmetrization argument,
and that explains the condition on characteristic. Note that
symmetrization is, in a sense, the opposite of straightening.

In view of Theorem \ref{symb} it is enough for the proof of
Theorem \ref{decomp} to show that a product
$\Delta=\delta_1\cdots\delta_p$ is in $I^\rho$ if
$\gamma_t(\Delta)\ge\gamma_t(\rho)$ for all $t$. For the analogy
to the proof of Lemma \ref{al-ga} set $\sigma=|\Delta|$. Then the
same induction works, since Lemma \ref{destr} implies in the
critical case that $\Delta$ is a linear combination of bitableaux
$\Delta'$, $|\Delta'|=\sigma'$.

\begin{remark}\rm
(a) Since we are mainly interested in asymptotic properties, we do
not discuss when the decomposition in \ref{decomp} is irredundant;
see \cite[(10.12)]{BV} and \cite[(10.13)]{BV} for a precise
result. Roughly speaking, the $I_t$-primary component is
irredundant if $I_t$ appears in the product (evidently) or if the
number of factors $I_u$ with $t<u<\min(m,n)$ is sufficiently
large. In particular, the decomposition of $I_t^k$ is irredundant
for $k\gg 0$ if $t<\min(m,n)$.

(b) In characteristic $2$ one has $I_3I_1\not\subseteq I_2^2$ if
$m,n\ge 4$; see \cite[(10.14)]{BV}.

(c) Independently of the characteristic the intersection in
Theorem \ref{decomp} is the integral closure of $I^\sigma$; see
Bruns \cite{Br}.
\end{remark}

Because of Lemma \ref{al-ga} we can replace the $\gamma$-functions
by $\alpha$-functions in the description of the standard bases of
products and powers. For powers one obtains a very simple
statement:

\begin{proposition}\label{power}
Suppose that $\chara K=0$ or $\chara K>\min(t,m-t,n-t)$. Then
$I_t^k$ has a basis consisting of all standard bitableaux $\Sigma$
with $\alpha_k(\Sigma)\ge kt$.
\end{proposition}

It is possible to derive this proposition from Theorem
\ref{decomp} by ``diagram arithmetic'', but it easier to prove it
directly. Let $V$ be the $K$-vector space generated by all
standard bitableaux $\Sigma=\delta_1\cdots \delta_u$ with
$\alpha_k(\Sigma)\ge kt$. That $I_t^k\subseteq V$ follows again
from the fact that straightening does not decrease shape. The
converse inclusion follows immediately (and much easier than
Theorem \ref{decomp}) from Lemma \ref{destr}: if $|\delta_i|<t$
for some $i$, $1\le i\le k$, then there is also an index $j$ in
this range such that $|\delta_j|>t$. Lemma \ref{destr} allows us
to increase $|\delta_i|$ at the expense of $|\delta_j|$.

\section{Gr\"obner bases, initial ideals and initial algebras}
\label{SectGB}

The aim of this section is to recall the definitions and some
important properties of Gr\"obner bases, monomial orders, initial
ideals and initial algebras. For further information on the theory
of Gr\"obner bases we refer the reader to the books by Eisenbud
\cite{Eis}, Kreuzer and Robbiano \cite{KRo}, Sturmfels \cite{Stu2}
and Vasconcelos \cite{Va}. For the so-called Sagbi bases and
initial algebras one should consult Conca, Herzog and Valla
\cite{CHV}, Robbiano and Sweedler \cite{RS}, and \cite[Chapter
11]{Stu2}.

Throughout this section let $K$ be a field, and let $R$ be the
polynomial ring $K[X_1,\dots,X_n]$. A \emph{monomial} (or power
product) of $R$ is an element of the form $X^\alpha=\prod_{i=1}^n
X_i^{\alpha_i}$ with $\alpha\in\NN^n$. A \emph{term} is an element
of the form $\lambda m$ where $\lambda$ is a non-zero element of
$K$ and $m$ is a monomial. Let $M(R)$ be the $K$-basis of $R$
consisting of all the monomials of $R$. Every polynomial $f\in R$
can be written as a sum of terms. The only lack of uniqueness in
this representation is the order of the terms. If we impose a
total order on the set $M(R)$, then the representation is uniquely
determined, once we require that the monomials are written
according to the order, from the largest to the smallest. The set
$M(R)$ is a semigroup (naturally isomorphic to $\NN^n$) and a
total order on the set $M(R)$ is not very useful unless it
respects the semigroup structure.

\begin{definition}\label{monord} \rm
A \emph{monomial order} $\tau$ is a total order
$<_\tau$ on the set $M(R)$ which satisfies the following
conditions:
\begin{itemize}
\item[(a)] $1<_\tau m$ for all the monomials $m\in M(R)\setminus\{1\}$.

\item[(b)] If $m_1,m_2,m_3\in M(R)$ and $m_1<_\tau m_2$, then
$m_1m_3<_\tau m_2m_3$.
\end{itemize}
\end{definition}

From the theoretical as well as from the computational point of
view it is important that descending chains in $M(R)$ terminate:

\begin{remark}\label{wellorder} \rm
A \emph{monomial order} on the set $M(R)$ is a well-order, i.e.
every non-empty subset of $M(R)$ has a minimal element.
Equivalently, there are no infinite descending chains in $M(R)$.

This follows from the fact that every (monomial) ideal in $R$ is
finitely generated. Therefore a subset $N$ of $M(R)$ has only
finitely many elements that are minimal with respect to
divisibility. One of them is the minimal element of $N$.
\end{remark}

We list the most important monomial orders.

\begin{example}\label{lex-revlex} \rm
For monomials $m_1=X_1^{\alpha_1}\cdots X_n^{\alpha_n}$ and
$m_2=X_1^{\beta_1}\cdots X_n^{\beta_n}$ one defines
\begin{itemize}
\item[(a)] the \emph{lexicographic order} (Lex) by
$m_1 <_{\textup{Lex}} m_2$ iff  for some $k$ one has $\alpha_k<
\beta_k$ and $\alpha_i = \beta_i$ for $i<k$;

\item[(b)]  the \emph{degree lexicographic order} (DegLex) by $m_1
<_{\textup{DegLex}} m_2$ iff $\deg(m_1)<\deg(m_2)$ or
$\deg(m_1)=\deg(m_2)$ and $m_1 <_{Lex} m_2$;

\item[(c)] the \emph{(degree) reverse lexicographic order}
(RevLex) by $m_1 <_{\textup{RevLex}} m_2$ iff
$\deg(m_1)<\deg(m_2)$ or $\deg(m_1)=\deg(m_2)$ and for some $k$
one has $\alpha_k> \beta_k$  and $\alpha_i = \beta_i$ for $i>k$.
\end{itemize}
\end{example}

These three monomial orders satisfy $X_1>X_2>\dots>X_n$. More
generally, for every total order on the indeterminates one can
consider the Lex, DegLex and RevLex orders extending the order of
the indeterminates; just change the above definition
correspondingly.

From now on we fix a monomial order $\tau$ on (the monomials of)
$R$. Whenever there is no danger of confusion we will write $<$
instead of $<_\tau$. Every polynomial $f\neq 0$ has an unique
representation
$$
f=\lambda_1 m_1+\lambda_2 m_2 +\dots+\lambda_k m_k
$$
where $\lambda_i \in K\setminus\{0\}$ and $m_1,\dots,m_k$ are
distinct monomials such that $m_1>\dots >m_k$. The \emph{initial
monomial} of $f$ with respect to $\tau$ is denoted by
$\ini_\tau(f)$ and is, by definition, $m_1$.  Clearly one has
$$
\ini_\tau(fg)=\ini_\tau(f)\ini_\tau (g) \eqno{(1)}
$$
and $\ini_\tau(f+g)\leq \max_\tau\{\ini_\tau(f),\ini_\tau(g)\}$.
For example, the initial monomial of the polynomial
$f=X_1+X_2X_4+X_3^2$ with respect to the Lex order is $X_1$, with
respect to DegLex it is $X_2X_4$, and with respect to RevLex it is
$X_3^2$.

Given a $K$-subspace $V\neq 0$ of $R$, we define
$$
M_\tau(V)=\{ \ini_\tau(f) : f\in V\}
$$
and set
$$
\ini_\tau(V)=\mbox{ the $K$-subspace of $R$  generated by }
M_\tau(V).
$$
The space $\ini_\tau(V)$ is called the {\em space of the initial
terms of $V$.}  Whenever there is no danger of confusion we
suppress the reference to the monomial order and use the notation
$\ini(f)$, $M(V)$ and $\ini(V)$.

Any positive integral vector $a=(a_1,\dots, a_n)\in\NN^n$ induces
a graded structure on $R$, called the \emph{$a$-grading}. With
respect to the $a$-grading the indeterminate $X_i$ has degree
$a(X_i)=a_i$. Every monomial $X^\alpha$ is $a$-homogeneous of
$a$-degree $\sum \alpha_ia_i$, and the $a$-degree $a(f)$ of a
non-zero polynomial $f\in R$ is the largest $a$-degree of a
monomial in $f$. Then $R=\bigoplus_{i=0}^\infty R_i$ where $R_i$
is the \emph{$a$-graded component} of $R$ of degree $i$, i.e.\ the
span of the monomials of $a$-degree $i$. With respect to this
decomposition $R$ has the structure of a positively graded
$K$-algebra \cite[Section 1.5]{BH}. The elements of $R_i$ are
\emph{$a$-homogeneous} of $a$-degree $i$. We say that a vector
subspace $V$ of $R$ is \emph{$a$-graded} if it is generated, as a
vector space, by homogeneous elements. This amounts to the
decomposition $V=\bigoplus_{i=0}^\infty V_i$ where $V_i=V\cap
R_i$.

\begin{proposition}\label{basis}
Let $V$ be a $K$-subspace of $R$.
\begin{itemize}
\item[(a)] If $m\in M(V)$ then there exists $f_m\in V$ such that
$\ini(f_m)=m$. The polynomial $f_m$ is uniquely determined if we
further require that the support of $f_m$ intersects $M(V)$
exactly in $m$ and that $f_m$ has leading coefficient $1$.

\item[(b)] $M(V)$ is a $K$-basis of $\ini(V)$.

\item[(c)] The set $\{ f_m : m\in M(V)\}$ is a $K$-basis of $V$.

\item[(d)] If $V$ has finite dimension, then $\dim(V)=\dim(\ini(V))$.

\item[(e)] Let $a\in \NN^n$ be a positive weight vector. Suppose $V$ is
$a$-graded, say $V=\bigoplus_{i=0}^\infty V_i$. Then
$\ini(V)=\bigoplus_{i=0}^\infty \ini(V_i)$. In particular, $V$ and
$\ini(V)$ have the same Hilbert function, i.e.
$\dim(V_i)=\dim(\ini(V)_i)$ for all $i\in\NN$.

\item[(f)] Let $V_1\subseteq V_2$ be $K$-subspaces of $R$. Then
$\ini(V_1)\subseteq \ini(V_2)$ and the (residue classes of the)
elements in $M(V_2)\setminus M(V_1)$ form a $K$-basis of the
quotient space $\ini(V_2)/\ini(V_1)$. Furthermore the set of the
(residue classes of the) $f_m$ with $f_m\in V_2$ and $m \in
M(V_2)\setminus M(V_1)$ is a $K$-basis of $V_2/V_1$ (regardless of
the choice of the $f_m$).

\item[(g)] The set of the (residue classes of the) elements in
$M(R)\setminus M(V)$ is a $K$-basis of $R/V$.

\item[(h)] Let $V_1\subseteq V_2$ be $K$-subspaces of $R$. If
$\ini(V_1)=\ini(V_2)$, then $V_1=V_2$.

\item[(i)] Let $V$ be a $K$-subspace of $R$ and $\sigma,\tau$
monomial orders. If $\ini_\tau(V)\subseteq \ini_\sigma(V)$, then
$\ini_\tau(V)= \ini_\sigma(V)$.
\end{itemize}

\end{proposition}

\begin{proof} (a) and (b) follow easily from the fact that the monomials
form a $K$-basis of $R$. For (a) we have to use that descending
chains in $M(R)$ terminate.

To prove (c) one notes that the $f_m$ are linearly independent
since they have distinct initial monomials. To show that they
generate $V$, we pick any non-zero $f\in V$ and set $m=\ini(f)$.
Then $m\in M(V)$ and we may subtract from $f$ a suitable scalar
multiple of $f_m$, say $g=f-\lambda f_m$, so that
$\ini(g)<\ini(f)$, unless $g=0$. Since $g\in V$, we may repeat the
procedure with $g$ and go on in the same manner. By Remark
\ref{wellorder}, after a finite number of steps we reach $0$, and
$f$ is a linear combination of the polynomials $f_m$ collected in
the subtraction procedure.

(d) and (e) follow from (b) and (c) after the observation that the
element $f_m$ can be taken $a$-homogeneous if $V$ is $a$-graded.

The first two assertions in (f) are easy. For the last we note
that $f_m$ can be chosen in $V_1$ if $m\in\ini(V_1)$.

The residue classes of the $f_m$ with $m \in M(V_2)\setminus
M(V_1)$ are linearly independent modulo $V_1$ since otherwise
there would be a non-trivial linear combination $g=\sum \lambda_m
f_m\in V_1$. But then $\ini(g)\in \ini(V_1)$, a contradiction
since $\ini(g)$ is one of the monomials $m$ which by assumption do
not belong to $M(V_1)$.

To show that the $f_m$ with $m\in M(V_2)\setminus M(V_1)$ generate
$V_2/V_1$ take some non-zero element $f\in V_2$ and set
$m=\ini(f)$. Subtracting a suitable scalar multiple of $f_m$ from
$f$ we obtain a polynomial in $V_2$ with smaller initial monomial
than $f$ (or $0$). If $m\in M(V_1)$, then $f_m\in V_1$. Repeating
the procedure we reach $0$ after finitely many steps. So $f$ can
be written as a linear combination of elements of the form $f_m$
with $m \in M(V_2)\setminus M(V_1)$ and elements of $V_1$, which
is exactly what we want.

(g) is a special case of (f) with $V_2=R$ since in this case we
can take $f_m=m$ for all $m\in M(R)\setminus M(V)$.

(h) follows from (f) since $\ini(V_1)=\ini(V_2)$ implies that the
empty set is a basis of $V_2/V_1$.

Finally, (i) follows from (g) because an inclusion between the two
bases $\{ m\in M(R) : m\not\in M_\tau(V)\}$ and $\{ m\in M(R) :
m\not\in M_\sigma(V)\}$ of the space $R/V$ implies that they are
equal. \end{proof}

\begin{remark/definition}\label{alg-id}\nopagebreak\rm
(a) If $A$ is a $K$-subalgebra of $R$, then $\ini(A)$ is also a
$K$-subalgebra of $R$. This follows from equation (1) and from
\ref{basis}(a). The $K$-algebra $\ini(A)$ is called the
\emph{initial algebra of $A$} (with respect to $\tau$).

(b) If $A$ is a $K$-subalgebra of $R$ and $J$ is an ideal of $A$,
then $\ini(J)$ is an ideal of the initial algebra $\ini(A)$. This,
too, follows from equation (1) and from \ref{basis}(a).

(c) If $I$ is an ideal of $R$, then $\ini(I)$ is also an ideal of
$R$. This is a special case of (b) since $\ini(R)=R$.
\end{remark/definition}

\begin{definition}\label{Sag}\rm
Let $A$ be $K$-subalgebra of $R$. A subset $F$ of $A$ is said to
be a \emph{Sagbi basis} of $A$ (with respect to $\tau$) if the
initial algebra $\ini(A)$ is equal to the $K$-algebra generated by
the monomials $\ini(f)$ with $f\in F$.
\end{definition}

If the initial algebra $\ini(A)$ is generated, as a $K$-algebra,
by a set of monomials $G$, then for every $m$ in $G$ we can take a
polynomial $f_m$ in $A$ such that $\ini(f_m)=m$. Therefore $A$ has
a finite Sagbi basis iff $\ini(A)$ is finitely generated. However
it may happen that $A$ is finitely generated, but $\ini(A)$ is
not; see \cite{RS}.

\begin{definition}\label{Gro}\rm
Let $A$ be a $K$-subalgebra of $R$ and $J$ be an ideal of $A$. A
subset $F$ of $J$ is said to be a \emph{Gr\"obner basis} of $J$
with respect to $\tau$ if the initial ideal $\ini(J)$ is equal to
the ideal of $\ini(A)$ generated by the monomials $\ini(f)$ with
$f\in F$.
\end{definition}

If the initial ideal $\ini(J)$ is generated, as an ideal of
$\ini(A)$, by a set of monomials $G$, then for every $m$ in $G$ we
can take a polynomial $f_m$ in $J$ such that $\ini(f_m)=m$.
Therefore $J$ has a finite Gr\"obner basis iff $\ini(J)$ is
finitely generated. In particular, if $\ini(A)$ is a finitely
generated $K$-algebra, then it is Noetherian and so all the ideals
of $A$ have a finite Gr\"obner basis. Evidently, all the ideals of
$R$ have a finite Gr\"obner basis.

There is an algorithm to determine a Gr\"obner basis of an ideal
of $R$ starting from any (finite) system of generators, the famous
\emph{Buchberger algorithm}. Similarly there is an algorithm that
decides whether a given (finite) set of generators for a
subalgebra $A$ is a Sagbi basis. There also exists a procedure
that completes a system of generators to a Sagbi basis of $A$, but
it does not terminate if the initial algebra is not finitely
generated. If a finite Sagbi basis for an algebra $A$ is known, a
generalization of Buchberger's algorithm finds Gr\"obner bases for
ideals of $A$. We will not use these algorithms in this article
and so we refer the interested readers to the literature quoted at
the beginning of this section..
\bigskip

\noindent\emph{Initial objects with respect to weights.}\enspace
In order to present the deformation theory for initial ideals and
algebras we need to further generalize these notions and consider
initial objects with respect to weights. As pointed out above, any
positive integral weight vector $a=(a_1,\dots, a_n) \in \NN^n$
induces a structure of a positively graded algebra on $R$. Let $t$
be a new variable and set
$$
S=R[t].
$$
For $f=\sum \gamma_i m_i\in R$ with $\gamma_i\in K$ and monomials
$m_i$ one defines the \emph{$a$-homo\-ge\-niza\-tion} $\hom_a(f)$ of
$f$ to be the polynomial
$$
\hom_a(f)=\sum \gamma_i m_i t^{a(f)-a(m_i)}.
$$

Let $a'=(a_1,\dots,a_n,1)\in \NN^{n+1}$.  Clearly, for every $f\in
R$ the element $\hom_a(f)\in S$ is $a'$-homogeneous, and
$f=\hom_a(f)$ iff $f$ is $a$-homogeneous. One has
$$
\begin{aligned}
\ini_a(fg)&=\ini_a(f)\ini_a(g)\\
\hom_a(fg)&= \hom_a(f) \hom_a(g)
\end{aligned}
\qquad\text{for all } f,g\in R. \eqno{(2)}
$$
For every $K$-subspace $V$ of $R$ we set
\begin{align*}
\ini_a(V)&= \text{ the $K$-subspace of $R$ generated by
$\ini_a(f)$ with $f\in V$},\\
\hom_a(V)&=\text{ the $K[t]$-submodule of $S$ generated by
$\hom_a(f)$ with $f\in V$}.
\end{align*}

If $A$ is a $K$-subalgebra of $R$ and $J$ is an ideal of $A$, then
it follows from (2) that $\ini_a(A)$ is a $K$-subalgebra of $R$
and $\ini_a(J)$ is an ideal of $\ini_a(A)$. Furthermore
$\hom_a(A)$ is a $K[t]$-subalgebra of $S$ and $\hom_a(J)$ is an
ideal of $\hom_a(A)$. As for initial objects with respect to
monomial orders, $\ini_a(A)$ and $\hom_a(A)$ need not be finitely
generated $K$-algebras, even when $A$ is finitely generated. But
if $\ini_a(A)$ is finitely generated, we may find generators of
the form $\ini_a(f_1), \dots, \ini_a(f_k)$ with $f_1,\dots, f_k\in
A$. It is easy to see that the $f_i$ generate $A$. This follows
from the next lemma  in which we use the notation $f^\alpha=\prod
f_i^{\alpha_i}$ for a vector $\alpha \in \NN^k$ and the list
$f=f_1,\dots, f_k$.

\begin{lemma}\label{normalformweight}
Let $A$ be $K$-subalgebra of $R$. Assume that $\ini_a(A)$ is
finitely generated by $\ini_a(f_1), \dots, \ini_a(f_k)$ with
$f_1,\dots, f_k\in A$. Then every $F\in A$ has a representation
$$
F=\sum \lambda_i f^{\beta_i}
$$
where $\lambda_i\in K\setminus\{0\}$ and $a(F)\geq a(f^{\beta_i})$
for all $i$.
\end{lemma}

\begin{proof} By decreasing induction on $a(F)$. The case $a(F)=0$ being
trivial, we assume $a(F)>0$. Since $F\in A$ we have $\ini_a(F)\in
\ini_a(A)=K[\ini_a(f_1),\dots,\allowbreak\ini_a(f_k)]$. Since
$\ini_a(F)$ is an $a$-homogeneous element of the $a$-graded
algebra $\ini_a(A)$, we may write
$$
\ini_a(F)=\sum \lambda_i \ini_a(f^{\alpha_i})
$$
where $a(\ini_a(f^{\alpha_i}))=a(\ini_a(F))$ for all $i$. We set
$F_1=F-\sum \lambda_i f^{\alpha_i}$ and conclude by induction
since $a(F_1)<a(F)$ if $F_1\neq 0$. \end{proof}

The following lemma contains a simple but crucial fact:

\begin{lemma}\label{founda}
Let $A$ be a $K$-subalgebra of $R$ and $J$ be an ideal of $A$.
Assume that $\ini_a(A)$ is finitely generated by $\ini_a(f_1),
\dots, \ini_a(f_k)$ with $f_1,\dots, f_k\in A$. Let
$B=K[Y_1,\dots, Y_k]$ and take presentations
$$
\phi_1:B\to A/J \quad \mbox { and } \quad \phi:B\to
\ini_a(A)/\ini_a(J)
$$
defined by the substitutions $\phi_1(Y_i)=f_i \mod (J)$ and
$\phi(Y_i)=\ini_a(f_i) \mod(\ini_a(J))$. Set $b =(a(f_1), \dots,
a(f_k))\in \NN_+^k$. Then
$$
\ini_b(\Ker \phi_1)=\Ker \phi.
$$

\end{lemma}

\begin{proof}
As a vector space, $\ini_b(\Ker \phi_1)$ is generated by the
elements $\ini_b(p)$ with $p\in \Ker\phi_1$. Set $u=b(p)$. Then we
may write $p=\sum \lambda_i Y^{\alpha_i}+\sum \mu_j Y^{\beta_j}$
where $b(Y^{\alpha_i})=u$ and $b(Y^{\beta_j})<u$. The image
$F=\sum \lambda_i f^{\alpha_i}+\sum \mu_j f^{\beta_j}$ belongs to
$J$, and, hence, $\ini_a(F)\in \ini_a(J)$. Since
$b(Y^\gamma)=a(f^\gamma)$, it follows that $\ini_a(F)=\sum
\lambda_i \ini_a(f^{\alpha_i}$). Thus $\ini_b(p)\in \Ker\phi$, and
this proves the inclusion $\subseteq$.

For the other inclusion we lift $\phi_1$ and $\phi$ to
presentations
$$
\rho_1:B\to A \quad \mbox{  and } \quad \rho:B\to \ini_a(A),
$$
mapping $Y_i$ to $f_i$ and to $\ini_a(f_i)$, respectively. Take a
system of $b$-homogeneous generators $G_1$ of the ideal $\Ker
\rho$ of $B$ and a system of $a$-homogeneous generators $G_2$ of
the ideal $\ini_a(J)$ of $\ini_a(A)$. Every $g\in G_2$, being
$a$-homogeneous of degree $u=a(g)$, is of the form $g=\ini_a(g')$,
with $g'\in J$. Then $g'=\sum \gamma_i f^{\alpha_i}+\sum \mu_j
f^{\beta_j}$ with $a(f^{\alpha_i})=u$ and $a(f^{\beta_j})<u$.
Therefore $g=\sum \gamma_i \ini_a(f^{\alpha_i})$.

We choose the canonical preimage of the given representation of
$g$, i.e. $h_g=\sum \gamma_i Y^{\alpha_i}$. Then the set $G_1\cup
\{h_g: g \in G_2\}$ generates the ideal $\Ker \phi$. For all
$g\in G_2$ and $g'$ as above, the canonical preimage of the given
representation of $g'$, i.e. $h=\sum \gamma_i Y^{\alpha_i}+\sum
\mu_j Y^{\beta_j}$ is in $\Ker\phi_1$, and one has
$\ini_b(h)=h_g$.

It remains to show that $g\in \ini_b(\Ker\phi_1)$ for $g\in G_1$.
Every $g\in G_1$ is homogeneous, say of degree $u$, and hence
$g=\sum \lambda_i Y^{\alpha_i}$ with $b(Y^{\alpha_i})=u$. It
follows that $\sum \lambda_i \ini_a(f^{\alpha_i})=0$. Therefore
$\sum \lambda_i f^{\alpha_i}=\sum \mu_j f^{\beta_j}$ with
$a(f^{\beta_j})<u$ by Lemma \ref{normalformweight}. That is,
$g'=\sum \lambda_i Y^{\alpha_i}-\sum \mu_j Y^{\beta_j}$ is in
$\Ker\rho_1$. In particular, $g'\in \Ker\phi_1$ and
$\ini_b(g')=g$.
\end{proof}

A weight vector $a$ and a monomial order $\tau$ on $R$ define a
new monomial order $\tau a$ that ``refines'' the weight $a$ by
$\tau$:
$$
m_1>_{\tau a} m_2 \iff \left\{
\begin{aligned} a(m_1) &> a(m_2) \text{ or }\\
a(m_1)&= a(m_2) \text{ and } m_1 >_{\tau} m_2.
\end{aligned}\right.
$$
We  extend $\tau a$ to $S=R[t]$ by setting:
$$
m_1 t^i >_{\tau a' } m_2 t^j  \iff \left\{
\begin{aligned}
a'(m_1 t^i) &> a'(m_2 t^j ) \mbox{ or }\\
a'(m_1 t^i) &= a'(m_2 t^j )  \mbox{ and } i<j \mbox{ or } \\
a'(m_1 t^i) &= a'(m_2 t^j ) \mbox{ and } i=j \mbox{ and }
m_1>_{\tau} m_2.
\end{aligned}\right.
$$
By construction one has
$$
\ini_{\tau a }(f)=\ini_{\tau a' }(\hom_a(f))\quad \mbox{for all }
f\in R,\ f\neq 0.
$$
Given a $K$-subspace $V$ of $R$, we let $VK[t]$ denote the
$K[t]$-submodule of $S$ generated by the elements in $V$.

\begin{proposition}\label{deformation}
Let $a\in \NN^n$ be a positive integral vector and $\tau$ be a
monomial order on $R$. For every $K$-subspace $V$  of $R$ one has:
\begin{itemize}
\item[(a)] $\ini_{ \tau a}(V)=\ini_{ \tau a}(\ini_
a(V))=\ini_{\tau}(\ini_a(V))$,
\item[(b)] If either $\ini_\tau(V)\subseteq \ini_ a(V)$ or
$\ini_\tau(V)\supseteq \ini_ a(V)$, then $\ini_\tau(V)=\ini_
a(V)$,
\item[(c)] $\ini_{ \tau a}(V)K[t]=\ini_{ \tau a'}(\hom_a(V))$,
\item[(d)] The quotient $S/\hom_a(V)$ is a free $K[t]$-module.
\end{itemize}

\end{proposition}

\begin{proof} (a) Note that
$\ini_{\tau a}(f)=\ini_{\tau a}(\ini_a(f))=\ini_{\tau}(\ini_a(f))$
holds for every $f\in R$. It follows that the first space is
contained in the second and in the third. On the other hand, since
$\ini_a(V)$ is $a$-homogeneous, the monomials in its initial space
are initial monomials of $a$-homogeneous elements. But every
$a$-homogeneous element in $\ini_a(V)$ is of the form $\ini_a(f)$
with $f\in V$. This gives the other inclusions.

(b) If one of the two inclusions holds, then an application of
$\ini_\tau(..)$ to both sides yields that $\ini_\tau(V)$ either
contains or is contained in $\ini_\tau(\ini_a(V))$. By (a) the
latter is $\ini_{ \tau a}( V)$. Then by Proposition \ref{basis}(i)
we have that $\ini_\tau(V)=\ini_{ \tau a}( V)$. Next we may apply
\ref{basis}(h) and conclude that $\ini_\tau(V)=\ini_ a(V)$.

(c) For every $f\in R$ one has $\ini_{\tau a'}(\hom_a(f))
=\ini_{\tau a}(f)$. Thus $\ini_{ \tau a}(V)K[t] \subseteq \ini_{
\tau a'}(\hom_a(V))$. On the other hand, $\hom_a(V)$ is an
$a'$-homogeneous space. Therefore its initial space is generated
by the initial monomials of its $a'$-homogeneous elements. An
$a'$-homogeneous element of degree, say, $u$ in $\hom_a(V)$ has
the form $g=\sum_{i=1}^k \lambda_i t^{\alpha_i} \hom_a(f_i)$ where
$f_i\in V$ and $\alpha_i+a(f_i)=u$. If $\alpha_i=\alpha_j$ then
$a(f_i)=a(f_j)$ and $\hom_a(f_i+f_j)=\hom_a(f_i)+\hom_a(f_j)$. In
other words, we may assume that the $\alpha_i$ are all distinct
and, after reordering if necessary, that $\alpha_i<\alpha_{i+1}$.
Then $\ini_{\tau a'}(g)= t^{\alpha_1} \ini_{\tau a'}(\hom(f_1))=
t^{\alpha_1} \ini_{\tau a}(f_1)$. This proves the other inclusion.

(d) By (c) and Proposition \ref{basis}(b) the (classes of the)
elements $t^\alpha m$, $\alpha\in \NN$, $m\in M(R)\setminus M(V)$,
form a $K$-basis of $S/\hom_a(V)$. This implies that the set
$M(R)\setminus M(V)$ is a $K[t]$-basis of $S/\hom_a(V)$.
\end{proof}

The next proposition connects the structure of $R/I$ with that of
$R/\ini_a(R)$:

\begin{proposition}\label{flatfamily}
For every ideal $I$ of $R$ the ring $S/\hom_a(I)$ is a free
$K[t]$-module. In particular $t-\alpha$ is a non-zero divisor on
$S/\hom_a(I)$ for every $\alpha\in K$. Furthermore
$S/(\hom_a(I)+(t))\cong R/\ini_a(I)$ and
$S/(\hom_a(I)+(t-\alpha))\cong R/I$ for all $\alpha\neq 0$.
\end{proposition}

\begin{proof}
The first assertion follows from \ref{deformation}(d). It implies
that every non-zero element of $K[t]$ is a non-zero divisor on
$S/\hom_a(I)$. For $S/(\hom_a(I)\allowbreak+(t))\cong R/\ini_a(I)$
it is enough that $\hom_a(I)+(t)=\ini_a(I)+(t)$. This is easily
seen since for every $f\in R$ the polynomials $\ini_a(f)$ and
$\hom_a(f)$ differ only by a multiple of $t$. To prove that
$S/(\hom_a(I)+(t-\alpha))\cong R/I$ for every $\alpha\neq 0$, we
consider the graded isomorphism $\psi: R\to R$ induced by
$\psi(X_i)=\alpha^{-a_i}X_i$. One checks that
$\psi(m)=\alpha^{-a(m)}m$ for every monomial $m$ of $R$ and that
$\hom_a(f)-\alpha^{a(f)} \psi(f)$ is a multiple of $t-\alpha$ for
all the $f\in R$. So $\hom_a(I)+(t-\alpha)= \psi(I)+(t-\alpha)$,
which implies the desired isomorphism. \end{proof}

Now we use Proposition \ref{flatfamily} for comparing $R/I$ with
$R/\ini_a(I)$.

\begin{proposition}\label{transfer}\nopagebreak
\begin{itemize}
\item[(a)] $R/I$ and $R/\ini_a(I)$ have the same Krull dimension.
\item[(b)] The following properties are passed from $R/\ini_a(I)$
on to $R/I$: being reduced, a domain, a normal domain,
Cohen-Macaulay, Gorenstein.
\item[(c)] Suppose that $I$ is graded with respect to some positive
weight vector $b$. Then $\ini_a(I)$ is $b$-graded, too, and the
Hilbert functions of $R/I$ and $R/\ini_a(I)$ coincide.
\end{itemize}
\end{proposition}

\begin{proof}
Let us start with (b). The $K$-algebra $A=S/\hom_a(I)$ is
positively graded. Let $\mm$ denote its maximal ideal generated by
the residue classes of the indeterminates. Then $A$ has one of the
properties mentioned if and only if the localization $A'=A_\mm$
does so. In fact, all of the properties depend only on the
localizations of $A$ with respect to graded prime ideals, and such
localizations are localizations of $A'$ (see \cite[Section 1.5 and
Chapter 2]{BH}). The element $t$ is a non-zero-divisor in the
maximal ideal of the local ring $A'$. Moreover $A'/(t)$ is a
localization of $R/\ini_a(I)$, and the properties under
consideration are inherited by localizations. As just pointed out,
they ascend from $A'$ to $A$. Therefore it remains to prove that
they also ascend from $A'/(t)$ to $A'$.

It is elementary to show that $A'$ is reduced or an integral
domain if $A'/(t)$ has this property. For the Cohen-Macaulay and
Gorenstein property the same conclusion is contained in
\cite[2.1.3 and 3.1.9]{BH}.

It remains to consider normality. We show that $A'$ has the Serre
properties $(R_1)$ and $(S_2)$ if these hold for $A'/(t)$. Let
$\pp$ be a prime ideal of $A'$ with $\height\pp\le 1$. If
$t\in\pp$, then $\overline\pp=\pp/(t)$ is a minimal prime ideal of
$A'/(t)$, and the regularity of
$(A'/(t))_{\overline\pp}=A'_\pp/(t)$ implies that of $A'_\pp$. If
$t\notin\pp$, we choose a minimal prime overideal $\qq$ of
$\pp+(t)$. Since $A'$ is an integral domain and a localization of
an affine $K$-algebra, we must have $\height \qq=\height\pp+1$.
Moreover, $\height \qq/(t)=\height\qq-1=\height\pp$. It follows
that $(A'/(t))_{\overline\qq}$ is regular. So $A'_\qq$ and its
localization $A'_\pp$ are regular. Suppose now that $\height\pp\ge
2$. We must show that $\depth A'_\pp\ge 2$. If $t\in\pp$, then we
certainly have $\depth (A'/(t))_{\overline\pp}\ge 1$, since
$(A'/(t))_{\overline\pp}$ is regular or has depth at least $2$.
Otherwise we take $\qq$ as above. Then $\depth
(A'/(t))_{\overline\qq}\ge 2$, and $\depth A'_\qq\ge 3$. We choose
$u\neq 0$ in $\pp$. If $\depth A'_\pp=1$, then $\pp/(u)$ is an
associated prime ideal of $A'/(u)$. Moreover, we have $\depth
A'_\qq/(u)\ge 2$, and $\dim A'_\qq/\pp A'_\qq=1$. This is a
contradiction to \cite[1.2.13]{BH}: for a local ring $R$ one has
$\depth R\le \dim R/\pp$ for all associated prime ideals $\pp$ of
$R$.

It remains to transfer the properties listed in (b)  to
$A''=A/(t-1)\cong R/I$, the dehomogenization of $A$ with respect
to the degree $1$ element $t$. So $A''$ is the degree $0$
component of the graded ring $A[t^{-1}]$, and $A[t^{-1}]$ is just
the Laurent polynomial ring in the variable $t$ over $A''$. (This
is not hard to see; cf.\ \cite[Section 1.5]{BH}. The main point is
that the surjection $A\to A''$ factors through $A[t^{-1}]$ and
that the latter ring has a homogeneous unit of degree $1$.)
Finally, each of the properties descends from the Laurent
polynomial ring to $A''$.

For (a) one follows the same chain of descents and ascents: $\dim
R/I\allowbreak =\dim A''\allowbreak =\dim A''[t,t^{-1}]-1=\dim
A[t^{-1}]-1=\dim A-1$. For the very last equality one has to use
that $t$ is a non-zero-divisor in an affine $K$-algebra.

(c) First one should note that $\ini_a(I)$ is $b$-graded, since
the initial form of a $b$-homogeneous element is $b$-homogeneous,
too. We refine the weight $a$ by a monomial order $\tau$ and
derive the chain of equations
$$
H(R/\ini_a(I))=H(R/\ini_\tau(\ini_a(I)))=H(R/\ini_{\tau a}(I))
=H(R/I)
$$
for the Hilbert function $H(..)$ from \ref{basis}(e) and
\ref{deformation}(a). \end{proof}

Very often one wants to compare finer invariants of $R/\ini_a(I)$
and $R/I$, for example if $I$ is a graded ideal of $R$ with
respect to some other weight vector $b$. The next proposition
shows that the comparison is possible for graded components of
Tor-modules. One can prove an analogous inequality for
Ext-modules.

\begin{proposition}\label{flatfam}
Let $a,b$ positive integral vectors and let $J,J_1,J_2$ be
$b$-ho\-mo\-ge\-neous ideals of $R$ with $J \subseteq J_1$ and
$J\subseteq J_2$. Then $\ini_a(J), \ini_a(J_1), \ini_a(J_2)$ are
also $b$-homogeneous ideals, and one has
$$
\dim_K \Tor^{R/J}_i(R/J_1, R/J_2)_j\leq \dim_K
\Tor^{R/\ini_a(J)}_i ( R/\ini_a(J_1), R/\ini_a(J_2))_j
$$
where the graded structure on the $\Tor$-modules is inherited from
the $b$-graded structure of their arguments.
\end{proposition}

\begin{proof}
On $S$ we introduce a bigraded structure, setting $\deg
X_i=(b_i,a_i)$ and $\deg t=(0,1)$. The ideals $I=\hom_a(J)$,
$I_1=\hom_a(J_1)$ and $I_2=\hom_a(J_2)$ are then bigraded and so
are the algebras they define. We need a standard result in
homological algebra: if $A$ is a ring, $M,N$ are $A$-modules and
$x$ is a non-zero-divisor on $A$ as well as on $M$ then
$\Tor^A_i(M,N/xN)\cong \Tor^{A/xA}_i(M/xM, N/xN)$. (It is
difficult to find an explicit reference; for example, one can use
\cite[1.1.5]{BH}.) If, in addition, $x$ is a non-zero-divisor also
on $N$, then we have the short exact sequence $0\to N\to N\to
N/xN\to 0$. It yields the exact sequence
$$
0\to \CoKer \phi_i  \to \Tor^{A/xA}_i(M/xM, N/xN)\to \Ker
\phi_{i-1}\to 0
$$
where $\phi_i$ is multiplication by $x$ on $\Tor^A_i(M,N)$.

Set $A=S/\hom_a(J)$, $M=S/\hom_a(J_1)$, $N=S/\hom_a(J_2)$ and
$T_i=\Tor^{A}_i(M,N)$. Since the modules involved are bigraded,
so is $T_i$. Let $T_{ij}$ be the direct sum of all the components
of $T_i$ of bidegree $(j,k)$ as $k$ varies. Since $T_i$ is a
finitely generated bigraded $S$-module, $T_{ij}$ is a finitely
generated and graded $K[t]$-module (with respect to the standard
grading of $K[t]$). So we may decompose it as
$$
T_{ij}=F_{ij} \oplus G_{ij}
$$
where $F_{ij}$ is the free part and $G_{ij}$ is the torsion part,
which, being $K[t]$-graded, is a direct sum of modules of the form
$K[t]/(t^a)$ for various $a>0$. Denote the minimal number of
generators of $F_{ij}$ and $G_{ij}$ as $K[t]$-modules by $f_{ij}$
and $g_{ij}$, respectively. Now we consider the $b$-homogeneous
component of degree $j$ of the above short exact sequence with
$x=t$, which is a non-zero-divisor by Proposition
\ref{deformation}(d). It follows that
$$
\dim_K \Tor^{R/\ini_a(J)}_i ( R/\ini_a(J_1),
R/\ini_a(J_2))_j=f_{ij}+g_{ij}+g_{i-1,j}.
$$
If we take $x=t-1$ instead of $x$, then we have
$$
\dim_K \Tor^{R/J}_i ( R/J_1, R/J_2)_j=f_{ij}
$$
and this shows the desired inequality. \end{proof}

Note that one can also use Proposition \ref{flatfam} to transfer
the Cohen-Macaulay and Gorenstein properties from $R/\ini_a(I)$ to
$R/I$ if $I$ is $b$-graded.

If $I$ is graded with respect to the ordinary weight $(1,\dots,1)$
then it makes sense to ask for the Koszul property of $R/I$. By
definition, $R/I$ is Koszul if $\Tor_i^{R/I}(R/\mm,\allowbreak
R/\mm)_j$ is non-zero only for $i=j$. Backelin and Fr\"oberg
\cite{BF} give a detailed discussion of this class of rings.

\begin{corollary}\label{Koszul}
Suppose that $I$ is a graded ideal with respect to the weight
$(1,\dots,1)$. If, for some positive weight $a$, $\ini_a(I)$ is
generated by degree $2$ monomials, then $R/I$ is Koszul.
\end{corollary}

\begin{proof}
By a theorem of Fr\"oberg \cite{Fr} the algebra $R/\ini_a(I)$ is
Koszul, so that the corollary follows from \ref{flatfam}.
\end{proof}

In order to apply the previous results to initial objects defined
by monomial orders we have to approximate such orders by weight
vectors. This is indeed possible, provided only finitely many
monomials have to be considered.

\begin{proposition}\label{initbyweight}
Let $\tau$ be a monomial order on $R$.
\begin{itemize}
\item[(a)] Let $\{(m_1, n_1),\dots, (m_k,n_k)\}$ be a finite
set of pairs of monomials such that $m_i>_\tau n_i$ for all $i$.
Then there exists a positive integral weight $a\in \NN_+^n$ such
that $a(m_i)>a(n_i)$ for all $i$.

\item[(b)] Let $A$ be a $K$-subalgebra  of $R$ and $I_1,\dots,I_h$ be ideals
of $A$. Assume that $\ini_\tau(A)$ is finitely generated as a
$K$-algebra. Then there exists a positive integral weight $a\in
\NN_+^n$ such that $\ini_\tau(A)=\ini_a(A)$ and
$\ini_\tau(I_i)=\ini_a(I_i)$ for all $i=1,\dots, h$.
\end{itemize}
\end{proposition}

\begin{proof} (a) Set $m_i=X^{\alpha_i}$ and $n_i=X^{\beta_i}$ and
$\gamma_i=\alpha_i-\beta_i\in \ZZ^n$. Let $\Gamma$ be the $k\times
n$ integral matrix whose rows are the vectors $\gamma_i$. We are
looking for a positive column vector $a$ such that the
coefficients of the vector $\Gamma a$ are all $>0$. Suppose, by
contradiction, there is no such $a$. Then (one version of the
famous)  Farkas Lemma (see Schrijver \cite[Section 7.3]{Sij}) says
that there exists a linear combination $v=\sum c_i \gamma_i$ with
non-negative integral coefficients $c_i\in \NN$  such that $v\leq
0$, that is $v=(v_1,\dots, v_n)$ with $v_i\leq 0$. Then it follows
that $\prod_i m_i^{c_i}X^{-v}=\prod_i n_i^{c_i}$, which
contradicts our assumptions because the monomial order is
compatible with the semigroup structure.

(b) Let $F_0$ be a finite Sagbi basis of $A$, let $F_i$ be a
finite Gr\"obner basis of $I_i$ and set $F= \bigcup_i F_i$.
Consider the set $U$ of pairs of monomials $(\ini(f),m)$ where
$f\in F$ and $m$ is any non-initial monomial of $f$. Since $U$ is
finite, by (a) there exists $a\in \NN_+^n$ such that
$\ini_a(f)=\ini_\tau(f)$ for every $f\in H$. We show $a$ has the
desired property. Set $V_0=A$ and $V_i=I_i$. By construction the
(algebra for $i=0$ and ideal for $i>0$) generators of the
$\ini_\tau(V_i)$ belong to $\ini_a(V_i)$ so that
$\ini_\tau(V_i)\subseteq \ini_a(V_i)$. But then, by Proposition
\ref{deformation}(b), we may conclude that
$\ini_\tau(V_i)=\ini_a(V_i)$. \end{proof}

The main theorem of this section summarizes what we can say about
the transfer of ring-theoretic properties from initial objects.
For the Koszul property of subalgebras we must allow a
``normalization'' of degree. Suppose that $b$ is a positive weight
vector $b$, and suppose that a subalgebra $A$ is generated by
elements $f_1,\dots,f_s$ of the same $b$-degree $e\in\NN$. Then
every element $g$ of $A$ has $b$-degree divisible by $e$, and
dividing the $b$-degree by $e$ we obtain the \emph{$e$-normalized
$b$-degree} of $g$.

\begin{theorem}\label{reprbywei}
Let $\ini(..)$ denote the initial objects with respect to a
positive integral vector $a\in \NN^n$ or to a monomial order
$\tau$ on $R$. Let $A$ be a $K$-subalgebra of $R$ and $J$ be an
ideal of $A$. Suppose that $\ini(A)$ is finitely generated.
\begin{itemize}
\item[(a)] One has $\dim A/J=\dim \ini(A)/\ini(J)$.
\item[(b)]  If $\ini(A)/\ini(J)$ is reduced, a domain, a
normal domain, Cohen-Macaulay, or Gorenstein, then so is $A/J$.
\item[(c)] Let $b$ be a positive weight vector, and suppose that
$A$ and $J$ are $b$-graded. Then $A/J$ and $\ini(A)/\ini(J)$ have
the same Hilbert function.
\item[(d)] If, in addition to the hypothesis of (c), $\ini(A)/\ini(J)$
is Koszul with respect to $e$-normalized $b$-degree for some $e$,
then so is $A/J$.
\end{itemize}
\end{theorem}

\begin{proof} If the initial objects are formed with respect to a
monomial order then, by \ref{initbyweight}, we may represent them
as initial objects with respect to a suitable  positive integral
weight vector. Therefore in both cases the initial objects are
taken with respect to a positive integral weight $a$. By Lemma
\ref{founda} there exist a polynomial ring, say $B$, an ideal $H$,
and a positive weight $c$ such that $B/H\cong A/J$ and
$B/\ini_c(H)\cong \ini(A)/\ini(J)$. Furthermore, under the
hypothesis of (c), the weight $b$ can be lifted from the
generators of $\ini(A)$ to the indeterminates of $B$. Now the
theorem follows from Proposition \ref{transfer} and Lemma
\ref{Koszul}. \end{proof}

The theorem is usually applied in two extreme cases. In the first
case $A=R$, so that $\ini(A)=R$, and in the second case $H=0$, so
that  $\ini(J)=0$. There is a special instance that deserves a
separate statement.

\begin{corollary}\label{mainc2}
Let $A$ be $K$-subalgebra of $R$, and suppose that $\ini(A)$ is
finitely generated. If it is generated by monomials (e.g.\ if the
initial algebra is taken with respect to a monomial order) and
 normal, then $A$ is normal and Cohen-Macaulay.
\end{corollary}

\begin{proof}
By a theorem of Hochster \cite[6.3.5]{BH} the normal semigroup
algebra $\ini(A)$ is Cohen-Macaulay. \end{proof}

Sometimes one of the implications in Theorem \ref{reprbywei} can
be reversed:

\begin{corollary}\label{GorIni}
Let $b$ be a positive weight vector, and suppose that the
$K$-subalgebra $A$ is $b$-graded and has a Cohen-Macaulay initial
algebra $\ini(A)$. Then $A$ is Gorenstein iff $\ini(A)$ is
Gorenstein.
\end{corollary}

\begin{proof}
Since $\ini(A)$ is Cohen-Macaulay, $A$ is Cohen-Macaulay as well.
So both algebras are positively graded Cohen-Macaulay domains. By
a theorem of Stanley \cite[4.4.6]{BH}, the Gorenstein property of
such rings depends only on their Hilbert function, and both
algebras have the same Hilbert function. \end{proof}

We want to extend Theorem \ref{reprbywei} in such a way that it
allows us to determine the canonical module of $A/I$.

\begin{theorem}\label{CanIni}
Let $A$ be a subalgebra of $R$ as in Theorem \ref{reprbywei}, and
$I\subseteq J$ ideals of $A$. Suppose that $\ini(A)/\ini(I)$ and,
hence, $A/I$ are Cohen-Macaulay.
\begin{itemize}
\item[(a)] If $\ini(J)/\ini(I)$ is the canonical module of
$\ini(A)/\ini(I)$, then $J/I$ is the canonical module of $A/I$.
\item[(b)] Suppose in addition that $A,I,J$ are $b$-graded with
respect to a positive weight and $\ini(J)/\ini(I)$ is the
canonical module of $\ini(A)/\ini(I)$ (up to a shift). Then $J/I$
is the graded canonical module (up to the same shift).
\end{itemize}
\end{theorem}

\begin{proof}
For the sake of simplicity and since it is sufficient for our
applications, we restrict ourselves to the graded case in (b) and
$I=0$. Since $A$ is a Cohen-Macaulay positively graded $K$-algebra
which is a domain, to prove that $J$ is the canonical module of
$A$ it suffices to show that $J$ is a maximal Cohen-Macaulay
module whose Hilbert series satisfies the relation
$H_J(t)=(-1)^dt^kH_A(t^{-1})$ for some integer $k$ and $d=\dim A$
\cite[Thm.\ 4.4.5, Cor.\ 4.4.6]{BH}.

The relation $H_J(t)=(-1)^dt^kH_A(t^{-1})$ holds since, by
hypothesis, the corresponding relation holds for the initial
objects, and Hilbert series do not change by taking initial terms;
see Theorem \ref{reprbywei}(c). So it is enough to show that $J$
is a maximal Cohen-Macaulay module. But $\ini(J)$ is a height $1$
ideal since it is the canonical module \cite[Prop.\ 3.3.18]{BH},
and hence $J$, too, has height $1$. Therefore it suffices that
$A/J$ is a Cohen-Macaulay ring. But this follows again from
\ref{reprbywei}(b) since $\ini(A)/\ini(J)$ is Cohen-Macaulay (even
Gorenstein) by \cite[Prop.\ 3.3.18]{BH}. \end{proof}

In order to prove the general version, one chooses representations
$A/I\cong B/I_1$, $A/J\cong B/I_2$, $\ini(A)/\ini(I)\cong
B/\ini(I_1)$, $\ini(A)/\ini(J)\cong B/\ini(I_2)$ as in Lemma
\ref{founda}. For the application of \ref{flatfamily} one notes
that $t$ is a non-zero-divisor on all the residue class rings to
be considered and that $C$ is the canonical module of a positively
graded ring $R$ if $C/tC$ is the canonical module of $R/(t)$ for a
homogeneous non-zero-divisor of $R$ and $C$. (\cite{BH} contains
all the tools one needs to prove this claim.)

\section{The Knuth--Robinson--Schensted correspondence }
\label{SectKRS}

The Knuth--Robinson--Schensted correspondence (in our context)
sets up a bijection between standard bitableaux and monomials in
the ring $K[X]$. The passage from bitableaux to monomials is based
on the \emph{deletion} algorithm.

\begin{definition}\label{Deletion}\rm
\emph{Deletion} takes a standard tableau $A=(a_{ij})$, say of
shape ${(s_1,s_2,\dots)}$, and an index $p$ such that
$s_p>s_{p+1}$, and constructs from them a standard tableau $B$ and
a number $x$, determined as follows:
\begin{itemize}

\item[(1)] Define the sequence $k_p,k_{p-1},\dots,k_1$ by setting
$k_p=s_p$ and choosing $k_i$ for $i<p$ to be the largest integer
$\leq s_i$ such that $a_{i k_i}\leq a_{i+1, k_{i+1}}$.

\item[(2)] Define $B$ to be the standard tableau obtained from $A$
by
 \begin{itemize}
 \item removing $a_{ps_p}$ from the $p$th row, and
 \item replacing the entry $a_{i k_i}$ of the $i$th row
      by $a_{i+1, k_{i+1}}$, $i=1,\dots, p-1$
 \end{itemize}

\item[(3)] Set $x=a_{1k_1}$.

\end{itemize}
\end{definition}

The reader should check that $B$ is a again a standard tableau. It
has the same shape as $A$, except that its row $p$ is shorter by
one entry. Deletion has an inverse:

\begin{definition}\rm
\emph{Insertion} takes a standard tableau $A=(a_{ij})$, say of
shape ${(s_1,s_2,\dots)}$, and an integer $x$, and constructs from
them a standard tableau $B$ and an index $p$ determined as
follows:
\begin{itemize}
\item[(1)] Set $i=1$ and $B=A$.

\item[(2)] If $s_i=0$ or $x>a_{is_i}$, then add $x$ at the end of the $i$th row
of $B$, set $p=i$ and terminate.

\item[(3)] Otherwise let $k_i$ be the smallest $j$ such that $x\leq a_{js_i}$,
replace $b_{k_is_i}$ with $x$, set $x=a_{k_is_i}$ and $i=i+1$.
Then go to (2).
\end{itemize}
\end{definition}

Again it is easily checked that $B$ is a standard tableau whose
shape coincides with that of $A$, except that the row $p$ of $B$
is longer by one entry.

Deletion and Insertion are clearly inverse to each other: if
Deletion applied to input $(A,p)$ gives output $(B,x)$, then
Insertion applied to $(B,x)$ gives output $(A,p)$ and viceversa.

The \emph{Knuth--Robinson--Schensted correspondence}, KRS for
short, is at first defined as a bijective correspondence between
the set of the standard bitableaux (as combinatorial objects) and
the set of the two-line arrays of a certain type. The two-line
array $\KRS(\Sigma)$ is constructed from the standard bitableau
$\Sigma$ by an iteration of the following KRS-step:

\begin{definition}\label{1-step}\rm
Let $\Sigma=(A\mid B)=(a_{ij} \sep b_{ij})$ be a non-empty
standard bitableau. Then \emph{KRS-step} constructs a pair of
integers $(\ell,r)$ and a standard bitableau $\Sigma'$ as follows.

\begin{itemize}
\item[(1)] Choose the largest entry $\ell$ in the \emph{left}
tableau of $\Sigma$; suppose that $\{ (i_1,j_1),\allowbreak
\dots,\allowbreak (i_u,j_u)\}$, $i_1<\dots<i_u$, is the set of
indices $(i,j)$ such that $\ell=a_{ij}$. Set $p=i_u$ and $q=j_u$.
(We call $(p,q)$ the \emph{pivot position}.)

\item[(2)] Let $A'$ be the standard tableau obtained by removing
$a_{pq}$ from $A$.

\item[(3)] Apply Deletion to the pair $(B,p)$. The output is a
standard tableau $B'$ and an element $r$.

\item[(4)] set $\Sigma'=(A',B')$.
\end{itemize}
\end{definition}

Now $\KRS(\Sigma)$ is constructed from the outputs of a sequence
of KRS-steps:

\begin{definition}\label{2line}\rm
Let $\Sigma$ be a non-empty standard bitableau of shape
$s_1,\allowbreak s_2,\allowbreak \dots s_p$. Set $k=s_1+\dots+s_p$
and define the two-line array
$$
\KRS(\Sigma)=\begin{pmatrix}
\ell_1 & \ell_2 & \dots \ell_{k-1} & \ell_k \\
 r_1 & r_2 & \dots r_{k-1} & r_k
\end{pmatrix}
$$
as follows. Starting from $\Sigma_k=\Sigma$, the KRS-step
\ref{1-step}, applied to $\Sigma_i$ for $i=k,k-1,\dots, 1$,
produces the bitableau $\Sigma_{i-1}$ and the pair $(\ell_i,
r_i)$.
\end{definition}

We give an example in Figure \ref{Del}. The circles in the left
tableau mark the pivot position, those in the right mark the
chains of ``bumps'' given in \ref{Deletion}(2):
\begin{figure}[hbt]
\begin{gather*}
\begin{gathered}
\begin{picture}(4,2)(0,0)
\Box311 \Box213 \Box114 \Box015 \Box302 \Box206 \Ci20
\end{picture}
\hspace{0.5cm}
\begin{picture}(4,2)(0,0)
\Box011 \Box112 \Box213 \Box316 \Box004 \Box105 \Ci10 \Ci21
\end{picture}
\\
\begin{picture}(4,2)(0,0)
\Box311 \Box213 \Box114 \Box015 \Box302 \Ci01
\end{picture}
\hspace{0.5cm}
\begin{picture}(4,2)(0,0)
\Box011 \Box112 \Box215 \Box316 \Box004 \Ci31
\end{picture}
\\
\begin{picture}(4,2)(0,0)
\Box311 \Box213 \Box114 \Box302 \Ci11
\end{picture}
\hspace{0.5cm}
\begin{picture}(4,2)(0,0)
\Box011 \Box112 \Box215 \Box004 \Ci21
\end{picture}
\end{gathered}
\qquad\qquad\qquad
\begin{gathered}
\begin{picture}(2,2)(0,0)
\Box111 \Box013 \Box102 \Ci01
\end{picture}
\hspace{0.5cm}
\begin{picture}(2,2)(0,0)
\Box011 \Box112 \Box004 \Ci11
\end{picture}
\\
\begin{picture}(2,2)(0,0)
\Box111 \Box102 \Ci10
\end{picture}
\hspace{0.5cm}
\begin{picture}(2,2)(0,0)
\Box011 \Box004 \Ci00 \Ci01
\end{picture}
\\
\begin{picture}(2,2)(0,0)
\Box111 \Ci11
\end{picture}
\hspace{0.5cm}
\begin{picture}(2,2)(0,0)
\Box014 \Ci01
\end{picture}
\end{gathered}\\[12pt]
\KRS(\Sigma)=\begin{pmatrix}1&2&3&4&5&6\\4&1&2&5&6&3\end{pmatrix}
\longleftrightarrow X_{14}X_{21}X_{32}X_{45}X_{56}X_{63}.
\end{gather*}
\caption{The KRS algorithm}\label{Del}
\end{figure}

The two-line array $\KRS(\Sigma)$ has the following properties:
\begin{itemize}
\item[(a)] $\ell_{i}\leq \ell_{i+1}$ for all $i$,\nopagebreak
\item[(b)] if $\ell_{i}=\ell_{i+1}$ then $r_{i}\geq r_{i+1}$.
\end{itemize}

Property (a) is clear since the algorithm chooses $\ell_{i+1}\ge
\ell_i$. If $\ell_{i}=\ell_{i+1}$ then the pivot position of the
$(i+1)$th deletion step lies left of (or above) the pivot position
of the $i$th deletion step. Now it is easy to see that the element
pushed out by the $(i+1)$th step is not larger than that pushed
out by the $i$th step.

KRS gives a correspondence between standard bitableaux and
two-line arrays with properties (a) and (b). It is bijective since
it has an inverse. For the inversion of KRS one just applies the
Insertion algorithm to the bottom line of the array to build the
right tableau: at step $i$ it inserts $r_i$ in the tableau
obtained after the $(i-1)$th step. Simultaneously one accumulates
the left tableau by placing the element $\ell_i$ in the position
which is added to the right tableau by the $i$th insertion.

It remains to explain how we can interpret any two-line array
satisfying (a) and (b) as a monomial: we associate the monomial
$$
X_{\ell_1r_1} X_{\ell_2r_2\cdots } X_{\ell_kr_k},
$$
to it, clearly establishing the desired bijection. To sum up, we
have constructed a bijective correspondence between standard
bitableaux and monomials. If we restrict our attention to standard
bitableaux and monomials where the entries and the indeterminates
come from an $m\times n$ matrix, we get

\begin{theorem}\label{MainKRS}
The map $\KRS$ is a bijection between the set of standard
bitableaux on $\{1,\dots,m\}\times\{1,\dots,n\}$ and the monomials
of $K[X]$.
\end{theorem}

This theorem proves the second half of the straightening law: the
KRS correspondence says that in every degree $d$ there are as many
standard bitableaux as monomials. Since we know already that the
standard bitableaux generate the space of homogeneous polynomials
of degree $d$, we may conclude that they must be linearly
independent.

\begin{remark}\rm
In the fundamental paper \cite{Kn} Knuth extensively treats the
KRS correspondence for \emph{column standard} bitableaux with
increasing columns and non-decreasing rows. (Deletion still bumps
the entries row-wise, but the condition $a_{i k_i}\leq a_{i+1,
k_{i+1}}$ in \ref{Deletion}(1) must be replaced by $a_{i k_i}<
a_{i+1, k_{i+1}}$). The same point of view is taken in Fulton
\cite{Fu}, Knuth \cite{Kn2}, Sagan \cite{Sa} and Stanley
\cite{Sta3}. The version we are using is the ``dual'' one (see
\cite[Section 5 and p.\ 724]{Kn}). The notes to Chapter 7 of
\cite{Sta3} contain a detailed historical discussion of the
correspondence.

Below we will consider decompositions of sequences into increasing
and non-increasing sequences. For column standard bitableaux these
attributes must be exchanged.
\end{remark}

We have just seen that KRS is a bijection between two bases of the
same vector space. Therefore we can \emph{extend KRS to a
$K$-linear automorphism} of $K[X]$. The automorphism KRS does not
only preserve the total degree, but even the $\ZZ^m\dirsum\ZZ^n$
degree introduced above: in fact, no column or row index gets
lost. However, note that KRS it is not a $K$-algebra isomorphism:
it acts as the identity on polynomials of degree $1$ but it is not
the identity map. It would be interesting to have some insight in
the properties of KRS as a linear map, like, for instance, its
eigenvalues and eigenspaces.

\begin{remark}\label{KRSprop}\rm
 We note some important properties of KRS:

(a) KRS commutes with transposition of the matrix $X$: Let $X'$ be
a $n\times m$ matrix of indeterminates, and let $\tau:K[X]\to
K[X']$ denote the $K$-algebra isomorphism induced by the
substitution $X_{ij}\mapsto X'_{ji}$; then
$\KRS(\tau(f))=\tau(\KRS(f))$ for all $f \in K[X]$. Note that it
suffices to prove the equality when $f$ is a standard
bi\-ta\-bleau. Herzog and Trung \cite[Lemma 1.1]{HT} point out how
to translate Knuth's argument from column to row standard
tableaux.

(b) All the powers $\Sigma^k$ of a standard bitableau are again
standard, and one has
$$
\KRS(\Sigma^k)=\KRS(\Sigma)^k.
$$
This is not hard to check: $k$ successive deletion steps on
$\Sigma^k$ act like a single deletion step on $k$ copies of
$\Sigma$.

(c) If $\Sigma$ is a minor $[a_1 a_2 \dots a_t \sep b_1 b_2 \dots
b_t]$ then $\KRS(\Sigma)$ is just (the product of the elements on)
the main diagonal of $\Sigma$. More generally, if one the two
tableaux of $\Sigma$ is ``nested'', i.e.\ the set of entries in
each row contains the entries in the next row, then $\KRS(\Sigma)$
is the product of the main diagonals of $\Sigma$. (This is easy to
see if the right tableau is nested; one uses (a) for
transposition.)

Note, however, that in general $\KRS(\Sigma)$ need not to be one
of the monomials which appear in the expansion of $\Sigma$. In
other words, $\KRS$ does not simply select a monomial of the
polynomial $\Sigma$; there is no algebraic relation between
$\KRS(\Sigma)$ and the monomials appearing in $\Sigma$.
\end{remark}

In the application of KRS to Gr\"obner bases of determinantal
ideals it will be important to relate the shape of a standard
bitableau $\Sigma$ to ``shape'' invariants of $\KRS(\Sigma)$. In
the KRS correspondence the right tableau, and hence its shape, is
determined solely by the bottom line of the corresponding two-line
array. Therefore we are lead to the following problem: Let
$r=r_1,r_2,\dots, r_k$ be a sequence of integers and let $\Ins(r)$
be the standard tableau determined by the iterated insertions of
the $r_i$. What is the relationship between the shape of $\Ins(r)$
and the sequence $r$?

A subsequence $v$ of $r$ is determined by a subset $U$ of
$\{1,\dots, k\}$: if $U=\{i_1,\dots, i_t\}$ with $i_1<\dots, <i_t$
then $v=v(U,r)=r_{i_1}, \dots, r_{i_t}$.  The \emph{length} of $v$
is simply the cardinality of $U$. A first answer to the question
above was given by Schensted and had indeed been a motivation of
his studies.

\begin{theorem}[Schensted]\label{schen}
The length of the first row of $\Ins(r)$ is the length of the
longest increasing subsequence of $r$.
\end{theorem}

Does the length of the $i$th row for $i>1$ have a similar meaning?
Actually, these lengths cannot be interpreted individually, but
some of their combinations reflect properties of the
decompositions of the sequence $r$ into increasing subsequences.
This is the content of Greene's extension of Schensted theorem,
which we will now explain. A decomposition of $r$ into
subsequences corresponds to a partition $U=(U_1,U_2,\dots, U_s)$
of the set $\{1,\dots, k\}$. The shape of the decomposition is
$(|U_1|,|U_2|,\allowbreak \dots,|U_s|)$: we always assume that
$|U_1|\geq |U_2|\geq \dots\geq |U_s|$. An \emph{inc-decomposition}
of $r$ is given by a partition $(U_1,U_2,\dots, U_s)$ of the set
$\{1,\dots, k\}$ such that the associated subsequences are
increasing.

In Section \ref{SectDet} we have defined the functions $\alpha_k$,
namely $\alpha_k(\lambda)=\sum_{i=1}^k \lambda_i$. Now We
introduce a variant $\alphah_k$ for sequences of integers $r$,
setting
$$
\alphah_k(r)=\max\{ \alpha_k(\lambda) : \mbox{ $r$ has a
inc-decomposition of shape $\lambda$} \}.
$$
Similarly we set
$$
\alpha_k(P)=\alpha_k(\lambda) \qquad\text{and}\qquad
\alpha_k(U)=\alpha_k(\lambda)
$$
for every tableau $P$ and every inc-decomposition $U$ of shape
$\lambda$.

Inc-decompositions are crucial for us since they describe
realizations of a monomial as an initial monomial of a product of
minors. This will be made more precise in the next section.
However, in Section \ref{SectAlgCM} it will turn out useful to
have also a measure for decompositions into \emph{non-increasing
subsequences}. For a shape $\lambda=(s_1,\dots,s_t)$ we define
$\lambda^*$ to be the \emph{dual} shape: the $i$th component of
$\lambda^*=(s_1^*,\dots,s_{s_1}^*)$ counts the number of boxes in
the $i$th column of a tableau of shape $\sigma$; formally
$$
s_i^*=|\{k: s_k\ge i\}|.
$$
The functions $\alpha_k$ are dualized to
$$
\alpha_k^*(\lambda)=\alpha_k(\lambda^*)  \qquad\text{and}\qquad
\alpha_k^*(P)=\alpha_k^*(\lambda)
$$
if $P$ is a tableau of shape $\lambda$. Analogously one defines
$\alpha_k^*(U)$ for a non-inc-decom\-po\-sition. However, the
passage from inc-decompositions to non-inc-decomposi\-tions
contains already a dualization, and we set
$$
\alphah_k^*(r)=\max\{\alpha_k(\lambda) : \mbox{ $r$ has a
non-inc-decomposition of shape $\lambda$} \}.
$$

We can rephrase the definition of $\alphah_k$ and $\alphah_k^*$
for sequences $r$ as follows: $\alphah_k(r)$ (respectively,
$\alphah_k^*(r)$) is the length of the longest subsequence of $r$
that can be decomposed into $k$ increasing (non-increasing)
subsequences.

\begin{theorem}[Greene] \label{Greene} For every sequence of
integers $r$ and every $k\ge 0$ we have
$$
\text{(a)}\quad\alphah_k(r)=\alpha_k(\Ins(r)),\qquad
\text{(b)}\quad\alphah_k^*(r)=\alpha_k^*(\Ins(r)).
$$
\end{theorem}

For $k=1$ one obtains Schensted's theorem from (a). (Schensted
also proved (b) for $k=1$.) The proof of Greene's theorem is based
on Knuth's basic relations: two sequences of integers $r$ and $s$
differ by a \emph{Knuth relation} if one of the following
conditions holds:
$$
\begin{array}{ll} (1) &
r=x_1,\dots,x_{i},y,z,w,x_{i+4},\dots,x_k
\mbox{\ \ and\ \ }\\
  & s=x_1,\dots,x_{i},y,w,z,x_{i+4},\dots,x_k,\\
  & \mbox{with\ \ } z\leq y<w;
\end{array}
$$
$$
\begin{array}{ll} (2) &
r=x_1,\dots,x_{i},z,w,y,x_{i+4},\dots,x_k
\mbox{\ \ and\ \ }\\
  & s=x_1,\dots,x_{i},w,z,y,x_{i+4},\dots,x_k, \\
  & \mbox{with\ \ } z < y \leq w.
\end{array}
$$

For every standard tableau $P$ there is a \emph{canonical
sequence} $r_P$ associated with $P$ such that $\Ins(r_P)=P$,
defined as follows: $r_P$ is obtained by listing the rows of $P$
from bottom to top, i.e.\
$$
r_P=p_{t1}p_{t2}\dots p_{ts_t} p_{t-1 1}p_{t-1 2}\dots
p_{t-1s_{t-1}}\dots p_{11}p_{12}\dots p_{1s_1};
$$
here $(s_1,\dots, s_t)$ is the shape of $P$.

\begin{theorem}[Knuth]
Let $r$ and $s$ be sequences. Then $\Ins(r)=\Ins(s)$ iff $r$ is
obtained from $s$ by a sequence of Knuth relations. In particular,
the canonical sequence associated with $\Ins(r)$ is obtained from
$r$ by a sequence of Knuth relations.
\end{theorem}

A detailed proof of the theorem for column standard tableaux is
contained in \cite{Kn}. On \cite[p.~724]{Kn} one finds an
explanation how to modify the Knuth relations and statements for
row standard tableaux.

Now Greene's theorem is proved as follows: for (a) one shows that
\begin{itemize}
\item[(i)]  $\alpha_k(P)=\alphah_k(r_P)$ for the canonical
sequence $r_P$ associated to a standard ta\-bleau $P$, and
\item[(ii)] $\alphah_k$ has the same value on sequences that
differ by a Knuth relation.
\end{itemize}
The same scheme works for (b). We confine ourselves to (a),
leaving (b) to the reader.

For the proof of (i) let $P=(p_{ij})$ and $\sigma=(s_1,\dots,s_u)$
be the shape of $P$. First of all, $r_P$ has the inc-decomposition
$p_{i1},\dots,p_{is_i}$, $i=1,\dots,r$. The decomposition has the
same shape as $P$. Hence $\alpha_k(P)\leq \alphah_k(r_P)$. On the
other hand, note that the columns of $P$ partition $r_P$ into
non-increasing subsequences. Therefore an increasing subsequence
can contain at most one element from each column. and the total
number of elements in a disjoint union of $k$ increasing
subsequences is at most
$\sum_{i=1}^{s_1}\min(k,s_i^*)=\alpha_k(P)$.

For the proof of (ii) we consider sequences $a$ and $b$ of
integers that differ by a Knuth relation. In order to prove that
$\alphah_k(a)=\alphah_k(b)$ it suffices to show that for every
inc-decomposition $G$ of $a$ there exists an inc-decomposition $H$
of $b$ such that $\alpha_k(G)\leq \alpha_k(H)$. So let $G$ be an
inc-decomposition of $a$ and let $z,w$ and $y$ as in the
definition of the Knuth relations. If $z$ and $w$ belong to
distinct subsequences in the decomposition $G$, then the
increasing subsequences of $G$ are not affected by the Knuth
relation, and we may take $H$ equal to $G$. (Strictly speaking, we
must change the partition of the index set underlying the
sequences by exchanging the positions of $z$ and $w$.)

It remains the case in which $z$ and $w$ belong to the same
increasing subsequence. Then $a$ must play the role of $r$, and
$y$ must belong to another subsequence in its decomposition $G$.
Let $u=p_1,z,w,p_2$ and $v=p_3,y,p_4$ denote those subsequences in
$g$ that contain $z$, $w$ and $y$. Here the $p_i$ are increasing
subsequences of $a$. Assume first the Knuth relation is of type
(1). We can rearrange the elements of the sequences $u$ and $v$
into increasing subsequences of $b$ in three ways:
$$ (\mbox{I}) \left\{
\begin{array}{l} u'=p_1,z,p_4\\ v'=p_3,y,w,p_2
\end{array}
\right.\quad (\mbox{II}) \left\{
\begin{array}{l} u'=p_1,w,p_2\\ v'=p_3,y,p_4\\ w'=z
\end{array}
\right.\quad (\mbox{III}) \left\{
\begin{array}{l} u'=p_1,y,w,p_2\\ v'=p_3,p_4\\ w'=z
\end{array}
\right.
$$

Suppose that both $u$ and $v$ contribute to $\alpha_k(G)$. Then we
replace $u$ and $v$ by the sequences $u'$ and $v'$ defined in (I),
obtaining an inc-decomposition $H$ of $b$ with $\alpha_k(H)\ge
\alpha_k(G)$: in fact, $H$ contains $k$ subsequences whose lengths
sum up to $\alpha_k(G)$.

If $u$ does not contribute to $\alpha_k(G)$, then we replace $u$
and $v$ by the sequences $u',v',w'$ in (II). The inc-decomposition
$H$ of $b$ consisting of the remaining subsequences of $a$ and
$u',v',w'$ has $\alpha_k(H)=\alpha_k(G)$.

Now suppose that $u$ contributes to $\alpha_k(G)$, but $v$ does
not. Then we replace $u$ and $v$ by the three sequences defined in
(III), with the same result as in the previous case.

The dual argument works in the case of a Knuth relation of type
(2). It completes the proof of Greene's theorem.
\medskip

Another series of very useful functions is given by the $\gamma_t$
introduced in Section~\ref{SectPow}, $ \gamma_t(\lambda) =\sum_i
\max(\lambda_i-t+1,0)$, where $t$ is a non-negative integer and
$\lambda$ is a shape. We extend the $\gamma_t$ to sequences in the
same way as the $\alpha_k$:
$$
\gammah_t(r)=\max\{\gamma_t(\lambda): r \text{ has an
inc-decomposition of shape} \lambda\}.
$$
Like the $\alpha_k$, the $\gamma_t$ are invariant under KRS:

\begin{theorem} \label{gamma}
For every sequence of integers $r$ we have $\gammah_t(r)=
\gamma_t(\Ins(r))$.
\end{theorem}

One can prove Theorem \ref{gamma} by arguments completely
analogous with those leading to Greene's theorem. This approach
has been chosen in \cite{BC1}. Alternatively one can use the
following lemma:

\begin{lemma}\label{a->g}
For each shape $\lambda$ the following holds:
$$
\gamma_t(\lambda)\ge u\quad\iff\quad \alpha_k(\lambda)\ge (t-1)k+u
\quad\text{for some }k,\ 1\le k\le t.
$$
\end{lemma}

We leave the easy proof to the reader, as well as the dual version
of Theorem \ref{gamma}.

\begin{remark}\label{not=shape}\rm
In spite of the Theorems \ref{Greene} and \ref{gamma}: in general
there does not exist an inc-decomposition of a sequence $r$ with
the same shape as $\Ins(r)$. The sequence $4,1,2,5,6,3$ as in
Figure \ref{Del} has no inc-decomposition of shape $(4,2)$ but the
shapes $(4,1,1)$ and $(3,3)$ occur, and this is enough for the
invariance of the functions $\alpha_k$ and $\gamma_t$.
\end{remark}

\section{KRS and Gr\"obner bases of ideals}
\label{SectGB-KRS}

Once and for all we now introduce a \emph{diagonal term order} on
the polynomial ring $K[X]$. With respect to such a term order, the
initial monomial $\ini(\delta)$ is the product of the elements on
the main diagonal of $\delta$ (for brevity we call this monomial
the main diagonal). There are various choices for a diagonal term
order. For example, one can take the lexicographic order induced
by the total order of the $X_{ij}$ that coincides with the
lexicographic order of the positions $(i,j)$.

\begin{remark}\label{MaxLI}\rm
It is not hard to show that distinct standard bitableaux $\Sigma$
of maximal minors have distinct initial monomials with respect to
a diagonal term order on $K[X]$. This proves the linear
independence of these standard bitableaux by \ref{basis}(d). Also
KRS is ``trivial'' for such $\Sigma$, since
$\KRS(\Sigma)=\ini(\Sigma)$. See \ref{KRSprop}(c) for a more
general statement.
\end{remark}

The power of KRS in the study of Gr\"obner bases for determinantal
ideals was detected by Sturmfels \cite{Stu1}. His simple, but
fundamental observation is the following. Assume that $I$ is an
ideal of $K[X]$ which has a $K$-basis of standard bitableaux $B$.
Then $\KRS(I)$ is a vector space of $K[X]$ that has two of the
properties of an initial ideal: it has a basis of monomials and it
has the same Hilbert function as $I$. Can we conclude that
$\KRS(I)$ is the initial ideal of $I$? Not in general, but if
$\KRS(I)\subseteq \ini(I)$ or the other way round, then equality
is forced by the Hilbert function. This argument yields

\begin{lemma}\label{KRS1}\nopagebreak
\begin{itemize}
\item[(a)] Let $I$ be an ideal of $K[X]$ which has a $K$-basis, say
$B$, of standard bitableaux, and let $S$ be a subset of $I$.
Assume that for all $\Sigma\in B$ there exists $s\in S$ such that
$\ini(s)\sep\KRS(\Sigma)$. Then $S$ is a Gr\"obner basis of $I$
and $\ini(I)=\KRS(I)$.

\item[(b)] Let $I$ and $J$ be homogeneous ideals such that
$\ini(I)=\KRS(I)$ and $\ini(J)=\KRS(J)$. Then
$\ini(I)+\ini(J)=\ini(I+J)=\KRS(I+J)$ and $\ini(I)\sect
\ini(J)=\ini(I\sect J)=\KRS(I\sect J)$.
\end{itemize}
\end{lemma}

\begin{proof}
(a) Let $J$ be the ideal generated by the monomials $\ini(s)$ with
$s\in S$. The hypothesis implies that $\KRS(I)\subseteq J\subseteq
\ini(I)$. Since the first and the third term have the same Hilbert
function it follows that $\KRS(I)=J=\ini(I)$. For (b) one uses
\begin{align*}
\KRS(I+J)&=\KRS(I)+\KRS(J)=\ini(I)+\ini(J)\subseteq \ini(I+J),\\
\KRS(I\sect J)&=\KRS(I)\sect \KRS(J)=\ini(I)\sect \ini(J)\supseteq
\ini(I\sect J),
\end{align*}
and concludes equality from the Hilbert function argument.
\end{proof}

Sturmfels applied Schensted's theorem to prove

\begin{theorem}\label{Sturm}
The $t$-minors of $X$ form a Gr\"obner basis of $I_t$, and
$\KRS(I_t)=\ini(I_t)$.
\end{theorem}

\begin{proof}
The set $B$ of standard bitableaux whose first row has length
$\geq t$ is a $K$-basis of $I_t$. Let $S$ be the set of the minors
of size $t$. By \ref{KRS1}, it is enough to show that for every
$\Sigma\in B$ there exists $\delta\in S$ such that $\ini(\delta)
\sep \KRS(\Sigma)$. Let $\ell$ and $r$ be the top and bottom
vector of the KRS image of $\Sigma$ so that $\KRS(\Sigma)=\prod
X_{\ell_i r_i}$. By Schensted's theorem \ref{schen} we can find an
increasing subsequence $r$ of length $t$, say $r_{i_1}<r_{i_2}<
\dots< r_{i_t}$ with $i_1<i_2< \dots< i_t$. But then
$\ell_{i_1}<\ell_{i_2}< \dots< \ell_{i_t}$ follows from property
$(b)$ of the KRS-image. In other words, the factor $\prod
X_{\ell_{i_j} r_{i_j}}$ of $\KRS(\Sigma)$ is the main diagonal and
hence the initial monomial of a $t$-minor. \end{proof}

\begin{remark}\rm
That the $t$-minors of $X$ form a Gr\"obner basis of $I_t$ has
been proved by several authors. To the best of our knowledge, the
result was first published by Narasimhan \cite{Na}. Independently,
a proof was given by Caniglia, Guccione, and Guccione \cite{CGG}.
The result was re-proved by Ma \cite{MaY}.

The ideal of maximal minors has better properties than the $I_t$
in general, not only in regard to its primary decomposition (see
Corollary \ref{maxi}), but also in regard to the Gr\"obner basis:
the maximal minors form a universal Gr\"obner basis, i.e.\ a
Gr\"obner basis for every term order on $K[X]$. This difficult
result was proved by Bernstein and Zelevinsky \cite{BeZ}.

Also in the case $t=2$ a universal Gr\"obner basis of $I_t$ is
known. It consists of binomials; see Sturmfels \cite{Stu2}.
\end{remark}

We fix an important observation, already used in the proof above:

\begin{remark}\label{variedec}\rm
Let $\Sigma$ is a standard bitableau, $\ell$ and $r$ be the top
and bottom vector of the KRS image of $\Sigma$. Then any
decomposition of the monomial $\KRS(\Sigma)$ into product of main
diagonals corresponds to a decomposition of $r$ into increasing
subsequences.

Consequently we extend the definition of $\alphah_k$ and
$\gammah_t$ to monomials by setting
$$
\alphah_k(M)=\alphah_k(r)\quad\text{and}\quad\gammah_t(M)=\gammah_t(r)
$$
where $r$ denotes the bottom row in the two line array
representing the monomial $M$, as discussed between Definition
\ref{2line} and Theorem \ref{MainKRS}.
\end{remark}

Lemma \ref{KRS1} leads us to introduce the following definition:

\begin{definition}\rm
Let $I$ be an ideal with a basis of standard bitableaux. Then we
say that $I$ is \emph{in-KRS} if $\ini(I)=\KRS(I)$; if, in
addition, the bitableaux (standard or not) $\Delta\in I$ form a
Gr\"obner basis, then $I$ is \emph{G-KRS}. In slightly different
words, an ideal $I$ with a basis of standard bitableaux is in-KRS
if for each $\Sigma\in I$ there exists $x\in I$ with
$\KRS(\Sigma)=\ini(x)$; it is G-KRS, if $x$ can always be chosen
as a bitableau.
\end{definition}

As a consequence of Lemma \ref{KRS1} one obtains

\begin{lemma}\label{GandIn} Let $I$ and $J$ be ideals with a basis of
standard bitableaux.
\begin{itemize}
\item[(a)] If $I$ and $J$ are G-KRS, then $I+J$ is also G-KRS.
\item[(b)] If $I$ and $J$ are in-KRS, then $I+J$ and $I\cap J$ are also in-KRS.
\end{itemize}
\end{lemma}

In general the property of being G-KRS is not inherited by
intersections as we will see below.

Now we are in the position to use the information we have
accumulated on determinantal ideals and on the KRS map to describe
Gr\"obner bases and/or initial ideals of powers, products and
symbolic powers of determinantal ideals.

\begin{theorem}\label{GB-symbPow}
For every $k\in \NN$ the symbolic power $I_t^{(k)}$ of $I_t$ is a
G-KRS ideal. Its initial ideal is generated, as a vector space, by
the monomials $M$ with $\gammah_t(M)\geq k$. In particular, a
Gr\"obner basis of $I_t^{(k)}$ is given by the set of bitableaux
$\Sigma$ with $\gamma_t(\Sigma)=k$ and no factor of size $<t$.
\end{theorem}

\begin{proof}
Let $S_1$ be the set of the products of minors $\Delta$ with
$\gamma_t(\Delta)\geq k$. One has $S_1\subseteq I_t^{(k)}$. By
virtue of Theorem \ref{gamma} and Remark \ref{variedec} we know
that for all standard bitableau $\Sigma$ with
$\gamma_t(\Sigma)\geq k$ there exists $\Delta $ in $S_1$ with
$\ini(\Delta)\mid\KRS(\Sigma)$. Thus it follows from \ref{KRS1}
that $S_1$ is a Gr\"obner basis of $I_t^{(k)}$ and $I_t^{(k)}$ is
G-KRS.

It remains to show that the initial term of any product of minors
$\Delta$ with $\gamma_t(\Delta)\geq k$ is divisible by the initial
term of a product of minors $\Delta_1$ without factors of size
$<t$ and with $\gamma_t(\Delta_1)=k$. If $\Delta$ has factors of
size $<t$, we simply get rid of them. If $\gamma_t(\Delta)>k$,
then we cancel $\gamma_t(\Delta)-k$ boxes in the bitableau with
the corresponding entries. In this way we get $\Delta_1$.
\end{proof}

Another important consequence of Lemma \ref{GandIn} and Theorem
\ref{decomp} is:

\begin{theorem}\label{in-prod}
Let $ t_1,\dots,t_r$ be positive integers and set $I=I_{t_1}\cdots
I_{t_r}$ and $g_i=\gamma_i(t_1, ,\dots,t_r)$. If $\chara K=0$ or
$\chara K>\min(t_i,m-t_i,n-t_i)$ for all $i$, then $I$ is in-KRS
and $\ini(I)$ is generated, as a K-vector space, by the monomials
$M$ with $\gammah_i(M)\geq g_i$ for all $i$.
\end{theorem}

Theorem \ref{in-prod} is satisfactory if one only wants to
determine the initial ideal of the product $I_{t_1}\cdots
I_{t_r}$, but it does not tell us how to find a Gr\"obner basis. A
natural guess is that any such ideal is G-KRS, i.e.\ a Gr\"obner
basis of $I_{t_1}\cdots I_{t_r}$ is given by the products of
minors (standard or not) which are in $I_{t_1}\cdots I_{t_r}$.
Unfortunately this is wrong in general.

\begin{example}
\label{Ex1}\rm Suppose that $m\ge 4$ and $\chara K=0$ or $>3$, and
consider the ideal $I_4I_2$. The monomial
$M=X_{11}X_{13}X_{22}X_{34}X_{43}X_{45}$ has $\gammah_4(M)=1$,
$\gammah_3(M)=2$, $\gammah_2(M)=4$, $\gammah_1(M)=6$. (We have
seen a similar example already in Remark \ref{not=shape}.) Hence,
by virtue of \ref{in-prod}, we know that $M\in \ini(I_4I_2)$. The
products of minors of degree $6$ in $I_4I_2$ have the shapes $6$
or $(5,1)$, or $(4,2)$. Clearly $M$ is not the initial monomial of
a product of minors of shape $6$ or of shape $5,1$. The only
initial monomial of a $4$-minor that divides $M$ is
$X_{11}X_{22}X_{34}X_{45}$ but the remaining factor $X_{13}X_{43}$
is not the initial monomial of a $2$-minor. Hence $M$ is not the
initial monomial of a product of minors that belongs to $I_4I_2$.
\end{example}

Nevertheless, if we confine our attention to powers of
determinantal ideals, the result is optimal.

\begin{theorem}
\label{g-for-powers} Suppose that $\chara K=0$ or $\chara
K>\min(t,m-t,n-t)$. Then $I_t^k$ is G-KRS and $\ini(I_t^k)$ is
generated, as a K-vector space, by the monomials $M$ with
$\alphah_k(M)\geq kt$. In particular, a Gr\"obner basis of $I_t^k$
is given by the products of minors $\Delta$ such that $\Delta$ has
at most $k$ factors, $\alpha_k(\Delta)=kt$, and $\deg \Delta=kt$.
Therefore $I_t^k$ has a minimal system of generators which is a
Gr\"obner basis.
\end{theorem}

\begin{proof}
Let $S_1$ be the set of the products of minors $\Delta$ with
$\alpha_k(\Delta)\geq kt$. By Proposition \ref{power} we know that
$S_1\subseteq I_t^{k}$. Greene's theorem \ref{Greene} and Remark
\ref{variedec} imply that for all standard bitableau $\Delta$ with
$\alpha_k(\Delta)\geq kt$ there exists $\Delta$ in $S_1$ with
$\ini(\Delta)\mid\KRS(\Sigma)$. Thus it follows from Theorem
\ref{KRS1} that $S_1$ is a Gr\"obner basis of $I_t^{k}$ and
$I_t^{k}$ is G-KRS.

It remains to show that for every product of minors $\Delta$ with
$\alpha_k(\Delta)\geq kt $ there exists a product of minors
$\Delta_1$ with at most $k$ factors, degree $kt$ and
$\alpha_k(\Delta_1)=kt$ such that $\ini(\Delta_1) | \ini(\Delta)$.
This step is as easy as the corresponding one in the proof of
Theorem \ref{GB-symbPow}: $\Delta_1$ is obtained from $\Delta$ by
skipping the rows of index $>k$ (if any) and deleting
$\alpha_k(\Delta)-kt$ boxes from the first $k$ rows (in any way).
\end{proof}

\begin{remark}\label{Messina}\rm
(a) We can obviously generalize \ref{in-prod} and
\ref{g-for-powers} as follows: Let $c_1,\dots,c_t\in\NN$,
$t=\min(m,n)$, and $V$ be the vector space spanned defined by all
standard monomials $\Sigma$ with $\gamma_i(\Sigma)\ge c_i$ (or
$\alpha_i(\Sigma)\ge c_i$) for all $i$; then $V$ is an ideal and
in-KRS. In fact, each of the inequalities defines a G-KRS ideal in
$K[X]$. We can even intersect $V$ with a homogeneous component of
$K[X]$ (with respect to the total degree or the
$\ZZ^m\dirsum\ZZ^n$-grading) to obtain an in-KRS vector space.

(b) In \cite{BC3} we have further analyzed the properties of being
G-KRS or in-KRS. If $\chara K=0$ or $>\min(m,n)$, then all ideals
$I$ defined by shape have a standard monomial basis and are
in-KRS: that $I$ is defined by shape means that it is generated by
products of minors and, for such a product $\Delta$, it depends
only on $|\Delta|$ whether $\Delta$ belongs to $I$.

Furthermore, an ideal defined by shape is G-KRS exactly if is the
sum of ideals $J(k,d)\cap I_1^u$ and, if $m=n$, $(J(k,d)\cap
I_1^u)I_m^v$ where the $J(k.d)$ play the same role for the
$\alpha$-functions as the symbolic powers do for the
$\gamma$-functions: $J(k,d)$ is generated by all bitableaux
$\Delta$ with $\alpha_k(\Delta)\ge d$.
\end{remark}

\begin{remark}\label{Kwie}\rm
Let $T$ be a new indeterminate, and consider the polynomial ring
$R'=K[X.T]$ where $X$, as usually, is an $m\times n$ matrix of
indeterminates with $m\le n$. The KRS-invariance of the functions
$\alpha_k$ has found another application to the ideal
$$
J=I_m+I_{m-1}T+\dots+I_1T^{m-1}+(T^m).
$$
and its powers; see Bruns and Kwieci\'nski \cite{BK}. With the
results accumulated so far, the reader can easily show that
$$
J^k=R'\biggl(\sum_{d=0}^{km} J(k,d)T^{km-d}\biggr).
$$
Let us extend the diagonal term order from $K[X]$ to $R'$ by first
comparing total degrees and, in the case of equal total degree,
the $X$-factors of the monomials. It follows that
$\ini(J^k)=\ini(J)^k$ and $J^k$ has a Gr\"obner basis of products
$\Delta T^{km-d}$ where $\Delta$ is a bitableaux of total degree
$d$ such that $\alpha_k(\Delta)\ge d$. The technique by which we
explore the Rees algebra of the ideal $I_t$ in Section
\ref{SectAlgCM} can also been applied to the Rees algebra of $J$;
see \cite{BK}.
\end{remark}

\section{Cohen-Macaulayness and Hilbert series of determinantal rings}
\label{SectCM+HF}

Hochster and Eagon \cite{HE} proved that the determinantal ring
$K[X]/I_t$ is Cohen-Macaulay. Their proof is based on the notion
of principal radical system; for this and several other approaches
see \cite{BV}. Abhyankar \cite{Abh} presented a formula for the
Hilbert function of $K[X]/I_t$ obtained by enumerating the
standard bitableaux in the standard basis of $K[X]/I_t$.

The goal of this section is to show how these results can be
proved by Gr\"obner deformation, i.e.\ by the study of the ring
$K[X]/\ini(I_t)$. By \ref{Sturm} we know that $\ini(I_t)$ is a
square-free monomial ideal. There are special techniques available
for the study of such ideals. We  briefly recall the main
properties and notions to be used; for more details we refer the
reader to \cite[Chapter 5]{BH}.

A simplicial complex on a set of vertices $V=\{1,\dots,n\}$ is a
set $\Delta$ of subsets $F$ of $V$ such that $G\in \Delta$
whenever $G\subseteq F$ and $F\in \Delta$. To any square-free
monomial ideal $I$ in a polynomial ring $R=K[X_1, \dots, X_n]$ one
can associate the (abstract) simplicial complex
$$
\Delta=\{ F\subseteq \{1,\dots, n\} : X_F\not\in I\}
$$
where $X_F=\prod_{i\in F} X_i$. Conversely, to any simplicial
complex $\Delta$ on the vertices $\{1,\dots, n\}$ one associates a
square-free monomial ideal $I$ by setting
$$
I=( X_F : F \not\in \Delta).
$$
The ring $K[\Delta]=K[X]/I$ is called the \emph{Stanley-Reisner
ring} associated to $\Delta$. One can study the homological
properties and the numerical invariants of $K[\Delta]$ by
analyzing the combinatorial properties and invariants of $\Delta$.
An element $F$ of $\Delta$ is called a \emph{face}; its
cardinality is denoted by $|F|$. The \emph{dimension} of $F$ is
$|F|-1$ and the dimension of $\Delta$ is $\max \{ \dim F | F \in
\Delta\}$. By ${\bf F}(\Delta)$ we denote the set of the facets of
$\Delta$, i.e.\ the maximal elements of $\Delta$ under inclusion.
Then $\Delta$ is said to be \emph{pure} if every facet has maximal
dimension, in other words, if $\dim F=\dim \Delta$ for all $F\in
{\bf F}(\Delta)$.

\begin{lemma}
\label{mult} The Krull dimension of $K[\Delta]$ is $\dim\Delta+1$,
and the multiplicity of $K[\Delta]$ equals the number of facets of
maximal dimension of $\Delta$.
\end{lemma}

\begin{proof}
This follows from the fact that the defining ideal of $K[\Delta]$
is radical and its minimal primes are of the form $(X_i : i\not\in
F)$ where $F$ is a facet of $\Delta$. \end{proof}

\begin{definition}\rm
A simplicial complex $\Delta$ is said to be \emph{shellable} if it
is pure and if its facets can be given a total order, say
$F_1,F_2,\dots, F_e$, so that the following condition holds: for
all $i$ and $j$ with $1\leq j<i\leq e$ there exist $v\in
F_i\setminus F_j$ and an index $k$, $1\leq k<i$, such that
$F_i\setminus F_k=\{v\}$. A total order of the facets satisfying
this condition is called a \emph{shelling} of $\Delta$.
\end{definition}

Shellability is a strong property, very well suited for inductive
arguments. Suppose that $F_1,\dots, F_e$ is a shelling of a
simplicial complex $\Delta$, let $\Delta_i$ denote the smallest
simplicial complex containing $F_1,\dots, F_i$, and $\Delta_i^*$
the smallest simplicial complex containing $F_j\sect F_i$ for all
$j<i$. Then $\Delta_i$ is obviously shellable and one has a short
exact sequence
$$
0\to K[\Delta_i]\to K[\Delta_{i-1}]\dirsum K[F_i]\to
K[\Delta_i^*]\to 0
$$
where $K[F_i]$ is the Stanley-Reisner ring of the simplex defined
by $F_i$, i.e.\ $K[F_i]$ is the polynomial ring on the set of
vertices of $F_i$. The fact that the given order of the facets is
a shelling translates immediately into an algebraic property:
$K[\Delta_i^*]$ is defined by a single monomial. Its degree is the
cardinality of the set
$$
c(F_i)= \{ v\in F_i : \mbox{ there exists } k<i \mbox{ such that }
F_i\setminus F_k=\{v\} \}.
$$
This implies

\begin{theorem}
\label{simp1} Let $\Delta$ be a shellable simplicial complex of
dimension $d-1$ with shelling $F_1,\dots, F_e$. Then:
\begin{itemize}
\item[(a)] The Stanley-Reisner ring $K[\Delta]$ is Cohen-Macaulay.
\item[(b)] The Hilbert series of $K[\Delta]$ has the form
$h(z)/(1-z)^d$ with $h(z)=\sum h_j z^j\in \ZZ[x]$, $h_0=1$ and
$h_j=\bigl|\{ i \in \{1,\dots, e\} : |c(F_i)|=j\}\bigr|$.
\end{itemize}
\end{theorem}

\begin{proof}
In view of the exact sequence above, one easily proves both
statements by induction on $e$, using the behavior of
Cohen-Macaulayness and Hilbert series along short exact sequences.
\end{proof}

Let us return to determinantal rings. As pointed out, the ideal
$\ini(I_t)$ is generated by square free monomials, namely the main
diagonal monomials of the $t$-minors of $X$. The corresponding
simplicial complex $\Delta_t$ consists of all the subsets of
$$
V=\{1,\dots, m\} \times \{1,\dots, n\}
$$
that do not contain a $t$-diagonal. (Note that the usual notation
of matrix positions differs from Cartesian coordinates by a
$90^\circ$ rotation!) The facets of $\Delta_t$ can be described in
terms of family of non-intersecting paths. To do this, we give $V$
a poset structure (certainly not the most natural one). We set
$$
(i,j)\leq (h,k) \iff i\leq h \mbox{ and } j\geq k.
$$
A subset $A$ of $V$ is said to be a \emph{chain} if each two
elements of $A$ are comparable in the poset $V$, and it is an
\emph{antichain} if it does not contain a pair of comparable
elements. It is easy to see that an antichain with $t$ elements
corresponds to the main diagonal of a $t$-minor. For $t=2$ the
simplicial complex coincides with the \emph{order }(or
\emph{chain}) complex of $V$: its faces are the chains and its
facets are the maximal chains of the poset $V$. For general $t$
the simplicial complex $\Delta_t$ is the set of those subsets of
$V$ which do not contain antichains of $t$ elements; it is called
the \emph{step $t$ order complex} of $V$.

A maximal chain of $V$ can described as a path in $V$. A
\emph{path} $P$ in $V$ from point $A$ to point $B$, with $A\leq
B$, is, by definition, an unrefinable chain with minimum $A$ and
maximum $B$. It can be written as a sequence
$$
P:\ \ A=(a_1,b_1),(a_2,b_2), \dots, (a_d,b_d)=B
$$
where
$$
(a_{i+1}, b_{i+1})-(a_i,b_i)=(1,0) \mbox{ or } (0,-1)
\quad\mbox{for all } i.
$$
A point $(a_k,b_k)$ is said to be a \emph{right-turn} of the path
$P$ if $1<k<d$ and
$$
(a_{k+1},b_{k+1})-(a_k,b_k)=(0,-1), \quad
(a_{k},b_{k})-(a_{k-1},b_{k-1})=(1,0).
$$
If one describes the lattice $V$ using either the Cartesian or the
matrix notation, then a right-turn of $P$ is exactly a point where
the path turns to the right. Given two sets of $s$ points
${\mathcal A}=A_1,\dots,A_s$ and ${\mathcal B}=B_1,\dots, B_s$ of
$V$, a set $F$ is said to be a \emph{family of non-intersecting
paths} from ${\mathcal A}$ to ${\mathcal B}$ if it can be
decomposed as $F=P_1\cup P_2\dots \cup P_{s}$ where $P_i$ is a
path from $A_i$ to $B_i$ and $P_i\sect P_j=\emptyset$ if $i\neq
j$. (We identify a family of paths with the set of points on its
paths. In the present setting this is allowed because the
decomposition above is unique.) We will then say that a point
$C\in F$ is a right turn of $F$ is it is a right turn of the path
to which it belongs.
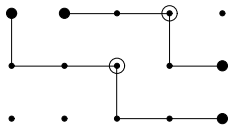
\begin{figure}[hbt]
\unitlength=0.7cm
\def\vertex{\circle*{0.10}}
\begin{center}
\begin{picture}(4,2)(0,0)
 \multiput(0,0)(0,1){3}{\multiput(0,0)(1,0){5}{\vertex}}
 \path(0,2)(0,1)(1,1)(2,1)(2,0)(4,0)
 \path(1,2)(3,2)(3,1)(4,1)
  \put(0,2){\circle*{0.20}}
 \put(1,2){\circle*{0.20}}
 \put(4,0){\circle*{0.20}}
 \put(4,1){\circle*{0.20}}
 \put(2,1){\circle{0.3}}
 \put(3,2){\circle{0.3}}
\end{picture}
\end{center}
\caption{A pair of non-intersecting paths with $2$ right-turns
($m=5$, $n=3$, $t=3$)}
\end{figure}
The \emph{length} of a path form $A=(x_1,x_2)$ to a point
$B=(y_1,y_2)$ depends only on $A$ and $B$ and is equal to
$y_1-x_1+x_2-y_2+1$.

\begin{proposition}
\label{facets} The facets of $\Delta_t$ are exactly the families
of non-intersecting paths from $(1,n),(2,n),\dots, (t-1,n)$ to
$(m,1),(m,2), \dots, (m,t-1)$.
\end{proposition}

\begin{proof}
A family of non-intersecting paths is in $\Delta_t$ since an
antichain intersects a chain in at most one point. That it is a
facet can be easily proved directly, but follows also from the
fact that any such family has dimension $(m+n-t+1)(t-1)-1$, which
is, by \ref{normal}, the dimension of $\Delta_t$ since
$\dim\Delta_t=\dim K[\Delta]-1=\dim K[X]/I_t-1$.

It remains to show that every facet $G\in \Delta_t$ is a family of
non-intersecting paths from $(1,n),(2,n),\dots, (t-1,n)$ to
$(m,1),(m,2), \dots, (m,t-1)$. The points of
$$
W=\{(a,b): b-a\geq n-t+1\text{ or }a-b\geq m-t+1\}
$$
do not belong to any $t$-antichain: there is not enough room. So
$G$ must contain $W$. We  put
$$
(a,b)\prec (c,d)  \iff a<c \mbox{  and }  b<d.
$$
By construction, distinct points $P,Q$ are either comparable with
respect to $<$ or comparable with respect to $\prec$, but not
both. So a chain with respect to $\prec$ is an antichain with
respect to $<$, and viceversa. For a set of points $A$ we denote
by $\Min_\prec(A)$ the set of the elements of $A$ which are
minimal with respect to $\prec$, i.e.\ the elements $P\in A$  such
that there is no $Q\in A$ with $Q\prec P$.  We then define:
$$
G_1=\Min_\prec(G) \quad\mbox {and}\quad
G_i=\Min_\prec\biggl(G\setminus \bigcup_{j=1}^{i-1} G_j\biggr)
\mbox{ for }\  i>1.
$$
This is called the \emph{light and shadow decomposition} (the
light here comes from the point $(1,1)$). The family of the $G_i$
satisfies the following conditions:
\begin{itemize}
\item[(a)] Every $G_i$ is a chain since otherwise $G_i$ would contain
two $<$-incomparable elements $P,Q$ and then either $P\prec Q$ or
$Q\prec P$ which is impossible.
\item[(b)] For every $P$ in $G_i$ there exist $P_1\in G_1,\dots,
P_{i-1}\in G_{i-1}$ such that $P_1\prec P_2\dots \prec P_{i-1}
\prec  P$. This is clear by construction.
\item[(c)] For $i=1,\dots t-1$ the set $G_i$ contains the points
$(i,n)$ and $(m,i)$. This follows from the fact that $W\subseteq
G$.
\item[(d)]  $G_i$ is empty for $i>t-1$ since otherwise we
would get a $t$-antichain in $G$ by (b).
\end{itemize}

It follows that $G$ is the disjoint union of the chains
$G_1,\dots, G_{t-1}$ and that each chain $G_i$ contains $(i,n)$
and $(m,i)$. We prove that each $G_i$ is indeed a path from
$(i,n)$ and $(m,i)$. Clearly $G_i$ cannot contain points which are
smaller than $(i,n)$ or larger than $(m,i)$ since those points
belong already to the $G_j$ with $j<i$. So it remains to show that
$G_i$ is saturated. Recall that $Q$ is said to be an upper
neighbor of $P$ if $P<Q $ an there is not an $H$ with $P<H<Q$. We
have to show the following
\medskip

\noindent {\em Claim}.\enspace  If $P=(a,b)$ and $Q=(c,d)$ belong
to $G_i$ and $Q>P$, but $Q$ is not an upper neighbor, then there
exists $H_0\in G_i$ such that $P<H_0<Q$.
\medskip

We set
$$
H=\left\{
\begin{array}{ll}
(c,b) & \mbox{ if } c\neq a \mbox{ and } b\neq d,\\
(c,d+1) &  \mbox{ if } c=a \mbox{ and } b\neq d,\\
(a+1,b) & \mbox{ if } c\neq a \mbox{ and } b= d.
\end{array}
\right.
$$

Since $P<H<Q$, if $H\in G_i$, then we are done; just set $H_0=H$.
If $H \not\in G_i$, then there are three possible cases:
\smallskip

\noindent (1)  If $H\in G_j$ for some $j<i$, then, by (b), we may
find $T_1, T_2 \in G_j$ such that $T_1\prec P$ and $T_2\prec Q$.
But $T_1$ and $T_2$ must be $<$-comparable with $H$.  This is a
contradiction because, by the choice of $H$, either $T_1\prec H$
or $T_2\prec H$. \smallskip

\noindent (2) If $H\in G_j$ for some $j>i$, then by $(b)$ there is
a $T\in G_i$ such that $T\prec H$. If $P<T<Q$, then we set
$H_0=T$. Otherwise either $T\leq P$ or $T\geq Q$. This is
impossible since $P<H<Q$. \smallskip

\noindent (3) If $H\not\in G$ then $G\cup \{H\}$ does not contain
a $t$-antichain since it has a decomposition into $t-1$ chains;
just add $H$ to $G_i$. This contradicts the maximality of $G$ and
concludes the proof. \end{proof}

It follows from the above description that every facet of
$\Delta_t$ has exactly $(m+n+1-t)(t-1)$ elements, so that
$\Delta_t$ is pure. Therefore the multiplicity of $K[\Delta_t]$
and, hence, that of $K[X]/I_t$ is given by the number of families
of non-intersecting paths from $(1,n),(2,n),\dots, (t-1,n)$ to
$(m,1),(m,2), \dots, (m,t-1)$. This number can be computed by the
Gessel-Viennot determinantal formula \cite{GV}: given two sets of
points ${\mathcal A}=A_1,\dots, A_s$ and ${\mathcal B}=B_1,\dots,
B_s$, the number $\Paths({\mathcal A} ,{\mathcal B})$ of families
of non-intersecting paths from ${\mathcal A}$ to ${\mathcal B}$ is
$$
\Paths({\mathcal A} ,{\mathcal B}) = \det ( \Paths(A_i, B_j)
)_{i,j=1,\dots, s}
$$
provided there is no family of non-intersecting paths from
${\mathcal A}$ to any non-trivial permutation of ${\mathcal B}$.
Here $\Paths(A_i, B_j)$ denotes the number of paths from $A_i$ to
$B_j$.

For $A_i=(i,n)$ and $B_j=(m,j)$ a simple inductive argument gives
$\Paths(A_i,\allowbreak B_j)= {m-i+n-j \choose m-i}$ and hence it
yields the formula
$$
e(K[X]/I_t)= \det \left( {m-i+n-j \choose m-i}
\right)_{i,j=1,\dots, t-1.}
$$
After some row and column operations one can evaluate the
determinant using Vandermonde's formula to obtain

\begin{theorem}\label{MultDet}
$$
e(K[X]/I_t)=\prod_{i=0}^{n-t}\frac{\binom{m+i}{t-1}}{\binom{t+i-1}{t-1}}
$$
\end{theorem}

The formula for $e(K[X]/I_t)$ is due to Giambelli (1903). The
proof above has been given by Herzog and Trung \cite{HT}. They
have generalized this approach (including Schensted's theorem
\ref{schen}) to the $1$-cogenerated ideals introduced in Section
\ref{SectDet}. (See Harris and Tu \cite{HaT} for a different
approach to Theorem \ref{MultDet}.)

Next we show that $\Delta_t$ is shellable. This is a special case
of a more general theorem due to Bj\"orner \cite[Thm.7.1]{Bj}.

\begin{theorem}
\label{shellable} The simplicial complex $\Delta_t$ is shellable.
More precisely, the fa\-cets of $\Delta_t$ can be given a total
order such that $c(F)$ is the set of right-turns of $F$ for each
facet $F$ of $\Delta_t$.
\end{theorem}

\begin{proof}
First we give a partial order to the set of paths connecting
points $A,B$ with $A\leq B$. For two paths $P_1$ and $Q_1$ from
$A$ to $B$ we write $P_1< Q_1$ if $P_1$ is ``on the right'' of
$Q_1$ as one goes from $A$ to $B$ (in Cartesian as well as in
matrix notation). This is a partial order.

Let $A_i=(i,n)$ and $B_i=(m,i)$ for $i=1,\dots, t-1$. Given two
families of non-intersecting paths $P=P_1,\dots. P_{t-1}$ and
$Q=Q_1,\dots,\allowbreak Q_{t-1}$ from $A_1,\dots, A_{t-1}$ to
${B_1,\dots, B_{t-1}}$ we set $P<Q$ if $P_i<Q_i$ for the largest
$i$ such that $P_i\neq Q_i$. We extend this partial order on the
set of families arbitrarily to a total order.

To prove that he resulting total order is indeed a shelling, one
takes two families $Q$ and $P$ with $P<Q$ and lets $i$ denote the
largest index $j$ with $P_j\neq Q_j$. Then $Q_i$ is not on the
right of $P_i$. It is easy to see (just draw a picture) that there
exists a right-turn, say $H$, of $Q_i$ which is (strictly) on left
of $P_i$. By the choice of $i$, the point $H$ does not belong to
$P_j$ for all $j$. So it suffices to show that for every
right-turn $H=(x,y)$ of $Q_i$ there is a family $R$ which is $<Q$
in the total order such that $Q\setminus R=\{H\}$.

This is easy if either $(x-1,y-1)$ does not belong to $Q_{i-1}$ or
$i=1$: just replace $(x,y)$ with $(x-1,y-1)$ in $Q_i$ to get a
path $Q_i'$, and then set $R=R_1,\dots, R_{t-1}$ with $R_j=Q_j$ if
$j\neq i$ and $R_i=Q_i'$. By construction $R<Q$ in the total
order.

It is a little more complicated to define $R$ when (by bad luck)
the element $(x-1,y-1)$ belongs to $P_{i-1}$. But if this is the
case, then $(x-1,y-1)$ must by a right-turn of $P_{i-1}$ (draw a
picture). If $(x-2,y-2)$ does not belong to $P_{i-2}$, we may
repeat the construction above: define $R$ as the family obtained
form $Q$ by replacing $(x,y)$ with $(x-1,y-1)$ in $Q_i$ and
$(x-1,y-1)$ with $(x-2,y-2)$ in $Q_{i-1}$. The general case
follows by the same construction. \end{proof}

Theorem \ref{shellable} has two important consequences. The first
is

\begin{theorem} \label{CMforIt}
The algebras $K[\Delta_t]$ and $K[X]/I_t$ are Cohen-Macaulay.
\end{theorem}

\begin{proof} By Theorem \ref{shellable} $\Delta_t$ is shellable,
and hence $K[\Delta_t]$ is Cohen-Macaulay by Theorem \ref{simp1}.
Since $K[\Delta_t]$ is $K[X]/\ini(I_t)$, it follows from Theorem
\ref{reprbywei} that $K[X]/I_t$ is Cohen-Macaulay as well.
\end{proof}

The second consequence is a combinatorial interpretation of the
Hilbert series of determinantal rings. We need some more notation
for it. Given two sets of $s$ points ${\mathcal A}$ and ${\mathcal
B}$, let $\Paths({\mathcal A},{\mathcal B})_k$ denote the numbers
of families of non-intersecting paths from $\mathcal A$ to
$\mathcal B$ with exactly $k$ right turns, and set
$\Paths({\mathcal A},{\mathcal B},z)=\sum_k \Paths({\mathcal
A},{\mathcal B})_k z^k$. In the case of just one starting and one
ending point, say $A$ and $B$, we denote this polynomial simply by
$\Paths(A,B,z)$. We have

\begin{theorem}
\label{HFofIt} The Hilbert series $H_t(z)$ of $K[\Delta_t]$ and
$K[X]/I_t$ is of the form
$$
H_t(z)=\frac{\Paths({\mathcal A},{\mathcal B},z)}{(1-z)^d}
$$
where $d=(m+n+1-t)(t-1)$ is the Krull dimension, ${\mathcal
A}=(1,n),(2,n),\dots,\allowbreak (t-1,n)$ and ${\mathcal
B}=(m,1),(m,2), \dots, (m,t-1)$.
\end{theorem}

For $t=2$, i.e.\ one starting and one end point, the polynomial
$\Paths({\mathcal A},{\mathcal B},z)$ can be easily computed by
induction on $n$ and $m$ and this yields the following formula:
$$
H_2(z)=\frac{\sum_k {m-1 \choose k}{n-1 \choose k}
z^k}{(1-z)^{m+n-1}}.
$$
It can be obtained also from the interpretation of $K[X]/I_2$ as
the Segre product of two polynomial rings.

By analogy with the Gessel-Viennot formula one may wonder whether
the polynomial $\Paths({\mathcal A},{\mathcal B},z)$ has a
determinantal expression as
$$
\det \left( \Paths( A_i, B_j,z) \right)_{i,j=1,\dots,t-1}
\eqno{(3)}
$$
This is obviously true if there is just one starting and ending
point, but cannot be true in general since $\Paths({\mathcal
A},{\mathcal B},0)=1$ and $\Paths( A_i, B_j,0)=1$ for all $i,j$.
But, very surprisingly, equality holds after a shift of degree if
the starting points are consecutive integral points on a vertical
line and the end points are consecutive integral points on a
horizontal line. This is essentially the content of

\begin{theorem}
\label{HFofIt2} The Hilbert series $H_t(z)$ of $K[\Delta_t]$ and
$K[X]/I_t$ is
$$
H_t(z)=\frac{ \det \left( \sum_k {m-i \choose k}{n-j \choose k}
z^k \right)_{i,j=1,\dots,t-1}}{ z^{{t-1\choose
2}}(1-z)^{(m+n+1-t)(t-1)}}.
$$

\end{theorem}

\begin{proof} Krattenthaler \cite{Kr} proved a determinantal formula for
$\Paths({\mathcal A},{\mathcal B},z)$ for general ${\mathcal
A},{\mathcal B}$. If the starting and end points are those
specified in \ref{HFofIt}, one can show that the determinant in
Krattenthaler's formula is equal to
$$
z^{-{t-1\choose 2}} \det \left( \sum_k {m-i \choose k}{n-j \choose
k} z^k \right)
$$
For the proof one has to describe the transformations in the
corresponding matrices. Details are to be found in \cite{CH1}.
\end{proof}

Krattenthaler's determinantal formula for the enumeration of
families of non-intersection paths with a given number of right
turns holds no matter how the starting and end points are located.
But this is not equal to the polynomial (3) even if one allows a
shift in degree. So the formula of \ref{HFofIt2} should be
regarded as an ``accident'' while the combinatorial description of
\ref{HFofIt} holds more generally, for instance, for algebras
defined by $1$-cogenerated ideals. On the other hand Krattenthaler
and Prohaska \cite{KrP} were able to show that the same
``accident'' takes place if the paths are restricted to certain
subregions called one-sided ladders. This proves a conjecture of
Conca and Herzog on the Hilbert series of one-sided ladder
determinantal rings; see \cite{CH1}.

\begin{remark}\label{GorIt}\rm
Since the rings $K[X]/I_t$ are Cohen-Macaulay domains, the
Go\-renstein ones among them are exactly those with a symmetric
numerator polynomial in the Hilbert series \cite[4.4.6]{BH}. By a
tedious analysis of the formula for the Hilbert series (see
\cite{CH1}) one can prove that $K[X]/I_t$, $t\ge 2$, is Gorenstein
if and only if $m=n$, a result due to Svanes. It is however more
informative to determine the canonical module of $K[X]/I_t$ for
all shapes of matrices; see \cite[Section 8]{BV} or
\cite[7.3]{BH}.
\end{remark}

\begin{remark}\label{I2Kosz}\rm
It follows from Theorem \ref{reprbywei} that the ring $K[X]/I_2$
is Koszul. We can also conclude that the homogeneous coordinate
ring $K[\MM_m]$ (with $m\le n$) of the Grassmannian is Koszul. To
this end we represent it as the residue class ring of a polynomial
ring $S$ whose indeterminates are mapped to the $m$-minors of $X$.
Then we refine the partial order $\prec$ of $m$-minors to a linear
order, lift that order to the indeterminates of $S$, and choose
the Revlex term order on $S$. The elements of $S$ representing the
Pl\"ucker relations form a Gr\"obner basis of the ideal defining
$K[\MM_m]$.
\end{remark}

\begin{remark}\label{ladder}\rm
Ladder determinantal rings are an important  generalization of the
classical determinantal rings. They  are defined by the minors
coming from certain subregions, called \emph{ladders},  of a
generic matrix. These objects have been introduced by Abhyankar in
his study of the singularities of Schubert varieties of flag
manifolds and have been investigated by many authors, including
Conca, Ghorpade, Gonciulea, Herzog, Knutson, Krattenthaler,
Kulkarni, Lakshmibai, Miller, Mulay, Narasim\-han, Prohaska, Rubey
and Trung \cite{Co2, Co3, CH2, Gho, GL1, GL2, HT, KM, KrP, KrR,
Ku1, Mul}. Ladder deter\-mi\-nan\-tal rings share many property
with classical determinantal rings, for instance they are
Cohen-Macaulay normal domains, the Gorenstein ones are completely
characterized in terms of the shape of the ladder, and there are
determinantal for\-mu\-las for the Hilbert series and functions.
Many of these results are derived from the combinatorial structure
of the Gr\"obner bases of the ideals of definition.
\end{remark}

\begin{remark}\label{SymPf}\rm
The ideal of $t$-minors of a symmetric matrix of indeterminates
and the ideal $\operatorname{Pf}_t$ of $2t$-pfaffians of an
alternating matrix of indeterminates can also be treated by
Gr\"obner basis methods based on suitable variants of KRS.

For pfaffians the method was introduced by Herzog and Trung
\cite{HT}. They used it to compute the multiplicity of
$K[X]/\operatorname{Pf}_t$. A determinantal formula for the
Hilbert series can be found in De Negri \cite{DN}; see also
Ghorpade and Krattenthaler \cite{GKr}. Bae\c tica has extended the
results of Section \ref{SectAlgCM} to the pfaffian case.

Conca \cite{Co1} has transferred the method of Herzog and Trung to
the symmetric case, introducing a suitable version of KRS. He
derived formulas for the Hilbert series and the multiplicity (see
\cite{Co1} for the latter).

Conca \cite{Co5} has attacked another class of determinantal
ideals by Gr\"obner basis methods, the ideals of minors of a
Hankel matrix. This case, like that of maximal minors, is ``easy''
since different standard products of minors have different initial
terms so that the essential point is to define the standard
products.

In addition to the generic case, Harris and Tu \cite{HaT} give
formulas of type \ref{MultDet} also in the symmetric and the
alternating case.
\end{remark}

\section{Algebras of minors: Cohen-Macaulayness and normality}
\label{SectAlgCM}

In this section we consider three types of algebras: the Rees
algebra $\Rees(I)$ of a pro\-duct $I=I_{t_1}\cdots I_{t_u}$ of
determinantal ideals, the symbolic Rees algebra
$\Rees^{\textup{symb}}(I_t)$ of $I_t$, and the algebra of
$t$-minors $A_t$, namely the $K$-subalgebra $K[\MM_t]$ of $K[X]$
generated by the $t$-minors. By studying their initial algebras we
will show that these algebras are normal and Co\-hen-Macaulay
(under a suitable hypothesis on the characteristic of $K$). In all
the cases the initial algebra is a finitely generated normal
semigroup ring and its description as well as its normality are
essentially a translation of the results of Section
\ref{SectGB-KRS} into the algebra setting.

It is convenient to embed all these algebras into a common
polynomial ring $S$, obtained by adjoining a variable $T$ to
$K[X]$,
$$
S=K[X,T]=K[X][T].
$$
For an ideal $I$ of $K[X]$ the Rees algebra $\Rees(I)$ of $I$ can
be described as $\Rees(I)=\Dirsum_k I^kT^k\subseteq S$. The
symbolic Rees algebra of $I_t$ is
$\Rees^{\textup{symb}}(I_t)=\Dirsum_k I_t^{(k)}T^k\subseteq S$ and
the algebra of minors $A_t$ can be realized as the subalgebra of
$S$ generated by the elements of the form $\delta T$ where
$\delta$ is a minor of size $t$. (One only uses that all
$t$-minors have the same degree as elements of $K[X]$.)

Let us first discuss some simple and/or classical cases. (They are
included in the general discussion below.) The Rees algebra of the
polynomial $K[X]$ with respect to the ideal $I_1$, its irrelevant
maximal ideal, can be represented as a determinantal ring. In
fact, let $R=K[X_1,\dots,X_n]$ where the $X_i$ are pairwise
different indeterminates. Then the substitution $X_i\mapsto X_i$,
$Y_i\mapsto X_iT$, $i=1,\dots,n$, yields the isomorphism
$\Rees(X_1,\dots,X_n)\cong K[X,Y]/I_2(U)$ where
$$
U=\begin{pmatrix} X_1&\dots&X_n\\Y_1&\dots&Y_n\end{pmatrix}.
$$
For the isomorphism it is enough to note that the 2-minors of $U$
are mapped to $0$ by the substitution and that $I_2(U)$ is a prime
ideal of height $n-1$ so that $\dim K[X,Y]/I_2(U)=n+1=\dim
\Rees(X_1,\dots,X_n)$. It follows that the Rees algebra is a
normal Cohen-Macaulay domain. It is Gorenstein only in the cases
$n=1,2$.

The other extreme case $t=\min(m,n)$ is also much simpler than the
general one. Eisenbud and Huneke \cite{EH} have shown that
$\Rees(I_t)$ is an algebra with straightening law on a wonderful
poset. In particular it is Cohen-Macaulay. By Proposition
\ref{maxi} $\Rees(I_t)=\Rees^{\textup{symb}}(I_t)$. This implies
normality since symbolic powers of primes in $K[X]$ are integrally
closed. See \cite[Section 9]{BV} or Bruns, Simis and Trung
\cite{BST} for generalizations.

For $A_t$ the case $t=1$ is completely trivial, since $A_1=K[X]$.
In the opposite extreme case $t=\min(m,n)$, say $t=m\le n$, the
algebra $K[\MM_m]=A_m$ is the homogeneous coordinate ring of the
Grassmannian of $m$-dimensional subspaces of the vector space
$K^n$, as discussed in Section \ref{SectDet}. This algebra is a
factorial Gorenstein ring; see \cite{BV}.

We have seen in Corollary \ref{dimGrass} that $\dim A_m=m(n-m)+1$.
However, if $t<\min(m,n)$, then $\dim A_t=\dim K[X]=mn$. Indeed,
the indeterminates $X_{ij}$ are algebraic over the quotient field
of $A_t$. It is enough to show this for a $(t+1)\times(t+1)$
matrix $X$. The entries of the adjoint matrix $\tilde X$ of $X$
are in $A_t$. Therefore $(\det X)^t\in A_t$. It follows that the
entries of $X^{-1}=(\det X)^{-1}\tilde X$ are algebraic over
$\QF(A_t)$, and playing the same game again, we conclude
algebraicity for the entries of $X=(X^{-1})^{-1}$.

Incidentally, this discussion has revealed another simple case: If
$t=m-1=n-1$, then $A_t$ is generated by $mn=\dim A_t$ elements,
and so is isomorphic to a polynomial ring over $K$.

We turn to the general case. Powers and products of determinantal
ideals are intersections of symbolic powers; see Theorem
\ref{decomp}. It follows immediately that the Rees algebra of
$I_{t^1}\cdots I_{t_u}$ is the intersection of symbolic Rees
algebras of the various $I_t$ and their Veronese subalgebras. The
representation as an intersection is passed on to the initial
algebras: this is a consequence of the in-KRS property. To sum up:
the key part is the description of the initial algebra of the
symbolic Rees algebra of $I_t$. The rest, at least as far as
normality and Cohen-Macaulayness are concerned, will follow at
once.

So let us start with the symbolic Rees algebra of $I_t$. The
description of the symbolic powers in Proposition \ref{symb}
yields the following description of the symbolic Rees algebra:
$$
\Rees^{\textup{symb}}(I_t)=K[X]\left[I_tT,I_{t+1}T^2, \dots,
I_mT^{m-t+1}\right].
$$
Consider a diagonal term order on $K[X]$ and extend it arbitrarily
to a term order on $K[X,T]$. The initial algebra
$\ini(\Rees^{\textup{symb}}(I_t))$ of $\Rees^{\textup{symb}}(I_t)$
is then $\Dirsum_k \ini(I_t^{(k)} ) T^k$. The description of
$\ini\bigl(I_t^{(k)}\bigr)$ in Theorem \ref{GB-symbPow} yields the
following

\begin{lemma}
\label{KRSalg} The initial algebra
$\ini(\Rees^{\textup{symb}}(I_t))$ of the symbolic Rees algebra
$\Rees^{\textup{symb}}(I_t)$ is equal to
$$
K[X]\left[\ini(I_t)T,\ini(I_{t+1})T^2,\dots,\ini(I_m)
T^{m-t+1}\right].
$$
In particular, a monomial $MT^k$ is in
$\ini(\Rees^{\textup{symb}}(I_t))$ if and only if
$\gammah_t(M)\geq k$.
\end{lemma}

 The next step is to show that
$\ini(\Rees^{\textup{symb}}(I_t))$ is normal. This can be done
directly by using the convexity of the function $\gammah_t$ as in
\cite{BC1}. Instead we give a longer, but more informative
argument which involves the description of the initial algebra by
linear inequalities (for the exponent vectors of the monomials in
it). This description will be used in the next section to identify
the canonical modules of various algebras. The crucial fact is the
primary decomposition of $\ini(I_t^{(k)})$:

\begin{lemma}\label{symb-in1}
Let ${\bf F}_t$ denote the set of facets of $\Delta_t$, and, for
every $F\in {\bf F}_t$, let $P_F$ be the ideal generated by the
indeterminates $X_{ij}$ with $X_{ij}\not\in F$. Then
$$
\ini\bigl(I_t^{(k)}\bigr)=\Sect_{F \in {\bf F}_t} P_F^k .
$$
\end{lemma}

We have seen in Theorem \ref{GB-symbPow} that
$\ini\bigl(I_t^{(k)}\bigr)$ is generated by the monomials $M$ with
$\gammah_t(M)\geq k$. A monomial $M=\prod_{i=1}^s X_{a_ib_i}$ is
in $P_F^{k}$ if and only if the cardinality of $\{i : (a_i,b_i)
\not\in F\}$ is $\geq k$. Equivalently, $M$ is in $P_F^{k}$ if and
only if the cardinality of $\{i : (a_i,b_i) \in F\}$ is $\leq
\deg(M) - k$. As a measure we introduce
$$
\w_t(M)=\max\bigl\{|A| : A\subseteq [1,\dots,s] \mbox{ and }
\{(a_i,b_i) : i \in A\}\in \Delta_t\bigr\}.
$$
Then a monomial $M$ is in $\Sect_{F\in {\bf F}_t} P_F^{k}$ if and
only if $\w_t(M)\leq \deg(M) -k $, or, equivalently,
$\deg(M)-\w_t(M) \geq k$. Now Proposition \ref{symb-in1} follows
from

\begin{lemma}\label{symb-in2}
Let $M$ be a monomial. Then $\gammah_t(M)+\w_t(M)=\deg(M)$.
\end{lemma}

We reduce this lemma to a combinatorial statement on sequences of
integers. For such a sequence $b$ we set
$$
w_t(b)=\max\{ \length(c) : c \mbox{ is a subsequence of $b$ and }
\gammah_t(c)=0\}.
$$
Let $M=\prod_{i=1}^s X_{a_ib_i}$ be a monomial. We order the
indices as in the KRS correspondence, namely $a_i\leq a_{i+1}$ for
every $i$ and $b_{i+1}\geq b_{i}$ whenever $a_i=a_{i+1}$. By
Remark \ref{variedec} we have $\w_t(M)=\w_t(b)$. To sum up, it
suffices to prove

\begin{lemma}\label{symb-in3}
One has $\gammah_t(b)+\w_t(b)=\length(b)$ for every sequence $b$
of integers.
\end{lemma}

\begin{proof}
We use part (b) of Greene's theorem \ref{Greene}:  the sum
$\alpha_k^*(\Ins(b))$ of the lengths of the first $k$ columns of
the insertion tableau $\Ins(b)$ of $b$ is the length of the
longest subsequence of $b$ that can be decomposed into $k$
non-increasing subsequences.

It follows that a sequence $a$ has no increasing subsequence of
length $t$ if and only if it can be decomposed into $t-1$
non-increasing subsequences. In fact, the sufficiency of the
condition is obvious, whereas its necessity follows from
Schensted's theorem \ref{schen} and the just quoted result of
Greene: if $a$ has no increasing subsequence of length $t$, then
all the rows in the insertion tableau $\Ins(a)$ have length at
most $t-1$. So $\alpha_{t-1}^*(\Ins(a))$ is the length of $a$, and
$a$ can be decomposed into $t-1$ non-increasing subsequences by
Greene's theorem.

Consequently $w_t(b)$ is the maximal length of a subsequence of
$b$ that can be decomposed into $t-1$ decreasing subsequences.
Applying Greene's theorem once more, we see that
$\w_t(b)=\alpha_{t-1}^*(\Ins(b))$. On the other hand,
$\gammah_t(b)=\gamma_t(\Ins(b))$ by Theorem \ref{gamma}. Since
$\gamma_t(\Ins(b))$ is the sum of the lengths of the columns of
$P$ of index $\geq t$, one has $w_t(b)+\gammah_t(b)=\length(b)$.
\end{proof}

Now we are ready to describe the linear inequalities supporting
$\ini(\Rees^{\textup{symb}}(I_t))$. To simplify notation we
identify monomials of $S$ with their exponent vectors in
$\RR^{mn+1}$. For every subset $F$ of $\{1,\dots, m\}\times
\{1,\dots, n\}$ we define a linear form $\ell_F$ on $\RR^{mn}$ by
setting $\ell_F(X_{ij})=1$ if $(i,j)\not\in F$ and $0$ otherwise.

\begin{theorem}\label{NaCMsym}
We extend $\ell_F$ for $F\in\mathbf{F}_t$ to a linear form $L_F$
on $\RR^{mn}\dirsum \RR$ by setting $L_F(T)=-1$. Then:
\begin{itemize}
\item[(a)] A monomial $MT^k$ is in the initial algebra
$\ini(\Rees^{\textup{symb}}(I_t))$ iff it has non-negative
exponents and $L_F(MT^k)\geq 0$ for all $F\in{\bf F}_t$.
\item[(b)] The initial algebra $\ini(\Rees^{\textup{symb}}(I_t))$
is normal and Cohen-Macaulay.
\item[(c)] The symbolic Rees algebra $\Rees^{\textup{symb}}(I_t)$
is normal and Cohen-Macaulay.
\end{itemize}

\end{theorem}

\begin{proof} (a) is a restatement of \ref{symb-in1}. Part (b) follows from (a)
and \cite[6.1.2, 6.1.4, 6.3.5]{BH}. Finally (c) follows from (b)
and Theorem \ref{reprbywei}. \end{proof}

As already mentioned, Theorem \ref{NaCMsym} has several
consequences. The first is

\begin{theorem}
\label{NaCMprod} Suppose that $\chara K=0$ or $\chara
K>\min(t_i,m-t_i,n-t_i)$ for all~$i$. Then
\begin{itemize}
\item[(a)] $\ini(\Rees(I_{t_1}\cdots I_{t_r}))$ is finitely
generated and normal,

\item[(b)] $\Rees(I_{t_1}\cdots I_{t_r})$ is Cohen--Macaulay and
normal.
\end{itemize}
\end{theorem}

\begin{proof}
Set $J=I_{t_1}\cdots I_{t_r}$. One has
$\ini(\Rees(J))=\Dirsum_{k\geq 0} \ini(J^k)T^k$, and, by Theorem
\ref{in-prod}, $\ini(J^k)=\Sect_{1\leq j\leq m}
\ini\bigl(I_j^{(kg_j)}\bigr)$. Hence
$$
\ini(\Rees(J))=\Sect_{1\leq j\leq m}\Dirsum_{k\geq 0}
\ini\bigl(I_j^{(kg_j)}\bigr)T^k.
$$
The monomial algebra $\Dirsum_{k\geq 0}
\ini\bigl(I_j^{(kg_j)}\bigr)T^k$ is isomorphic to the $g_j$th
Veronese subalgebra of the monomial algebra
$\ini(\Rees^{\textup{symb}}(I_j))$ (in the relevant case $g_j>0$
and equal to $K[X,T]$ otherwise). By \ref{NaCMsym} the latter is
normal and finitely generated, and therefore $\Dirsum_{k\geq 0}
\ini\bigl(I_j^{(kg_j)}\bigr)T^k$ is a normal, finitely generated
monomial algebra. Thus $\ini(\Rees(J))$ is finitely generated and
normal. In fact, the intersection of a finite number of finitely
generated normal monomial algebras is finitely generated and
normal. (This follows easily from standard results about normal
affine semigroup rings; see Bruns and Herzog \cite[6.1.2 and
6.1.4]{BH}.) For (b) one applies Corollary \ref{mainc2} again.
\end{proof}

We single out the most important case.

\begin{theorem}
\label{NaCMpow} Suppose that $\chara K=0$ or $\chara
K>\min(t,m-t,n-t)$. Then $\Rees(I_t)$ is Cohen--Macaulay and
normal.
\end{theorem}

\begin{remark}\label{OlCM}\rm
(a) The Cohen-Macaulayness of the Rees algebra of $I_t$ in the
case of maximal minors has been proved by Eisenbud and Huneke in
\cite{EH}, as pointed out above. For arbitrary $t$ and $\chara
K=0$, Bruns \cite{Br} has shown that $\Rees(I_t)$ and $A_t$ are
Cohen-Macaulay.

(b) In \ref{NaCMpow} the hypothesis on the characteristic is
essential. If $m=n=4$ and $\chara K=2$ then $\Rees(I_2)$ has
dimension $17$ and depth $1$; see \cite{Br}.

(c) In order to obtain a version of \ref{NaCMprod} that is valid
in arbitrary characteristic one must replace the Rees algebra by
its integral closure. The integral closure is always equal to the
intersection of symbolic Rees algebras that in non-exceptional
characteristic gives the Rees algebra itself (see \cite{Br}).
\end{remark}

\begin{remark}\label{hyperpro}\rm
One can describe the hyperplanes defining the initial algebra of
the Rees algebra of a product of $I_{t_1}\cdots I_{t_r}$ in terms
 of proper extensions of the linear forms $\ell_F$.
For every $j$ set $g_j=\gamma_j(t_1,\dots, t_r)$ and for every
$F\in {\bf F}_j$ extend $\ell_F$ to $L_F$ by setting
$L_F(T)=-g_j$. Then the initial algebra of the Rees algebra of
$I_{t_1}\cdots I_{t_r}$ is given by the inequalities
$L_F(MT^k)\geq 0$ for all $F\in {\bf F}_j$ and for all~$j$ (and
the non-negativity of the exponents of $MT^k$).
\end{remark}

For the algebra of minors $A_t$ we have

\begin{theorem}\label{AlgMin1}
Suppose that $\chara K=0$ or $>\min(t,m-t,n-t)$. Then the initial
algebra $\ini(A_t)$ is finitely generated and normal. Hence $A_t$
is a normal Cohen--Macaulay ring.
\end{theorem}

\begin{proof}
Let $V_t$ be the subalgebra of $S$ generated by the monomials of
the form $MT$ with $\deg M=t$, i.e. $V_t$ is (isomorphic to) the
$t$-Veronese subalgebra of $K[X]$. By construction
$A_t=\Rees(I_t)\sect V_t$. This clearly implies
$\ini(A_t)=\ini(\Rees(I_t))\sect V_t$. Since $\ini(\Rees(I_t))$ is
normal by Theorem \ref{NaCMprod} and $V_t$ is normal by obvious
reasons, $\ini(A_t)$ is normal. \end{proof}

\begin{remark}
\label{SagbiforAt}\rm (a) As for the other cases one can give a
description of the initial algebra of $A_t$ by linear inequalities
and equations.

(b) Although we have proved that the initial algebras of
$\Rees(I_t)$ and $A_t$ are finitely generated, we cannot specify a
finite Sagbi basis: we do not know what largest degree occurs in a
system of generators for their initial algebras.
\end{remark}

\section{Algebras of minors: the canonical module }
\label{SectAlgCan}

The goal of this section is to describe the canonical modules of
the algebras $\Rees(I_t)$ and $A_t$. The first step is to find the
canonical modules of their initial algebras. The characteristic of
the field $K$ will be either $0$ or $>min(t,m-t,n-t)$ throughout.

Recall that $V_t$ is the subalgebra of $S$ generated by all the
monomials of the form $MT$ where $M$ is a monomial in $K[X]$ of
degree $t$ (thus $V_t$ is isomorphic to the Veronese subalgebra of
$K[X]$). Part (a) of the following lemma is just a restatement of
\ref{decomp} and part (b) is a restatement of \ref{in-prod}.
However, (c) and (d) contain a somewhat surprising simplification
for $A_t$ and its initial algebra.

\begin{lemma} \label{algbygamma}\nopagebreak
\begin{itemize}
\item[(a)] A $K$-basis of $\Rees(I_t)$ is given by the set of the
elements $\Sigma T^k$ where $\Sigma$ is a standard bitableau with
$\gamma_i(\Sigma)\geq k(t+1-i)$ for all $i=1,\dots,t$.

\item[(b)] A $K$-basis of $\ini(\Rees(I_t))$ is given by the set of
the elements $M T^k$ where $M$ is a monomial of $K[X]$ with
$\gammah_i(M)\geq k(t+1-i)$ for all $i=1,\dots,t$.

\item[(c)] A $K$-basis of $A_t$ is given by the set of the elements
$\Delta T^k$ where $\Delta$ is a standard bitableau with
$\gamma_2(\Delta)\geq k(t-1)$ and $\deg(\Delta)=tk$.

\item[(d)] A $K$-basis of $\ini(A_t)$ is given by the set of the
elements $M T^k$ where $M$ is a monomial of $K[X]$ with
$\gammah_2(M)\geq k(t-1)$ and $\deg(M)=tk$.
\end{itemize}
\end{lemma}

\begin{proof}
As pointed out above, only (c) and (d) still need a proof. Since
$A_t=V_t\sect \Rees(I_t)$, a $K$ basis of $A_t$ is given by the
elements $\Delta T^k$ in the basis of $\Rees(I_t)$ with
$\deg(\Delta)=kt$. Similarly, since $\ini(A_t)=V_t\sect
\ini(\Rees(I_t))$, a $K$ basis of $\ini(A_t)$ is given by the
elements $MT^k$ in the basis of $\ini(\Rees(I_t))$ with
$\deg(M)=kt$. Now (c) and (d) result from the following statement:
if $\lambda$ is a shape such that $\sum_i \lambda_i=kt$ and
$\gamma_2(\lambda)\geq k(t-1)$, then $ \gamma_i(\lambda)\geq
k(t+1-i)$ for all $i=1,\dots,t$. We leave the proof to the reader;
it is to be found in \cite{BC4}.\end{proof}

\begin{remark}\label{valu}\rm
One can extend the valuations $\gamma_i$, $i=1,\dots,t$, to
$K[X,T]$ by choosing $\gamma_i(T)=-t+i-1$. Therefore Lemma
\ref{algbygamma}(a) contains a description of $\Rees(I_t)$ as an
intersection of $K[X,T]$ with discrete valuation domains. This
aspect is discussed in \cite{BC4} and \cite{BR}. Part (c) has a
similar interpretation.

However, since the equation
$\gammah_t(MN)=\gammah_t(M)+\gammah_t(N)$ does not always hold,
the functions $\gammah_t$ cannot be interpreted as valuations.
\end{remark}

We know that $\ini(\Rees(I_t))$ and $\ini(A_t)$ are normal. Hence
their canonical modules are the vector spaces spanned by all
monomials represented by integral points in the relative interiors
of the corresponding cones (see Bruns and Herzog
\cite[Ch.~6]{BH}). We have seen in \ref{hyperpro} how to describe
the semigroup of $\ini(\Rees(I_t))$ in terms of linear homogeneous
inequalities using the linear forms $L_F$ defined as follows: For
every $i=1,\dots, t$ and for every facet $F$ of $\Delta_i$ we
extend $\ell_F$ to a linear form $L_F$ on $\RR^{mn}\dirsum \RR$ by
setting $L_F(T)=-(t+1-i)$. We obtain

\begin{lemma}\label{caninRbylineq}
The canonical module of $\ini(\Rees(I_t))$ is the ideal of
$\ini(\Rees(I_t))$ whose $K$-basis is the set of the monomials $N$
of $S$ with all exponents $\geq 1$ and $L_F(N)\geq 1$ for every $F
\in {\bf F}_i$ and for $i=1,\dots, t$.
\end{lemma}

Let $\XX$ denote the product of all the variables $X_{ij}$ with
$(i,j)\in \{1,\dots,m\}\times \{1,\dots,n\}$. We can give a
description of the canonical module $\omega(\ini(\Rees(I_t)))$ in
terms of $\XX$ and the functions $\gammah_i$:

\begin{lemma}
\label{caninRe} The canonical module $\omega(\ini(\Rees(I_t)))$ of
$\ini(\Rees(I_t))$ is the ideal of $\ini(\Rees(I_t))$ whose
$K$-basis is the set of the monomials $MT^k$ of $S$ where $M$ is a
monomial of $K[X]$ with $\XX T \mid MT^k$ in $S$ and
$\gammah_i(M)\geq (t+1-i)k+1$ for all $i=1,\dots, t$.
\end{lemma}

\begin{proof}
It suffices to show that the conditions given define the monomials
described in \ref{caninRbylineq}. Let $N=MT^k$ be a monomial,
$M\in K[X]$. Then, for a given $i$, one has $L_F(N)\geq 1$ for
every $F \in {\bf F}_i$ if and only if $\ell_F(M)\geq k(t+1-i)+1$
for every $F \in {\bf F}_i$. By \ref{symb-in1} this is equivalent
to $M\in \ini\bigl(I_i^{(k(t+1-i)+1)}\bigr)$, which in turn is
equivalent to $\gammah_i(M)\geq k(t+1-i)+1$. To sum up,
$L_F(N)\geq 1$ for every $F \in {\bf F}_i$ and $i=1,\dots,t$ if
and only if $\gammah_i(M)\geq (t+1-i)+1$ for every $i=1,\dots,t$.
\end{proof}

Similarly the canonical module $\omega(\ini(A_t))$ has a
description in terms of the function $\gammah_2$:

\begin{lemma}
\label{caninAt} The canonical module $\omega(\ini(A_t))$ of
$\ini(A_t)$ is the ideal of $\ini(A_t)$ whose $K$-basis is the set
of the monomials $MT^k$ of $V_t$ where $M$ is a monomial of $K[X]$
with $\XX T \mid MT^k$ in $S$ and $\gammah_2(M)\geq (t-1)k+1$.
\end{lemma}

For ``de-initialization'' the following lemma is necessary. Its
part (b) asserts that $\XX$ is a ``linear'' element for the
functions $\gammah_i$.

\begin{lemma}\label{rhoXX}\nopagebreak
\begin{itemize}
\item[(a)] $\gammah_i(\XX)=(m-i+1)(n-i+1)$.
\item[(b)] Let $M$ be a monomial in $K[X]$. Then
$\gammah_i(\XX M)=\gammah_i(\XX)+\gammah_i(M)$ for every
$i=1,\dots,\min(m,n)$.
 \end{itemize}
\end{lemma}

\begin{proof}
Let $M$ be a monomial in the $X_{ij}$'s. We know that
$\gammah_i(M)\geq k$ if and only if $M\in
\ini\bigl(I_i^{(k)}\bigr)$. From \ref{symb-in1} we deduce that
$$
\gammah_i(M)=\inf\{ \ell_F(M) : F\in {\bf F}_i\}.
$$
Note that $\Delta_i$ is a pure simplicial complex of dimension one
less the dimension of $K[X]/I_i$. Thus
$\ell_F(\XX)=(m-i+1)(n-i+1)$ for every facet $F$ of $\Delta_i$. In
particular, $\gammah_i(\XX)=(m-i+1)(n-i+1)$.

Since $\ell_F(NM)=\ell_F(N)+\ell_F(M)$ for all monomials $N,M$ and
for every $F$, we have $\gammah_i(MN)\geq
\gammah_i(N)+\gammah_i(M)$. Conversely, let $G$ be a facet of
$\Delta_i$ such that $\gammah_i(M)=\ell_G(M)$. Then $\ell_G(\XX
M)= \ell_G(\XX)+ \ell_G(M)=\gammah_i(\XX)+\gammah_i(M)$. Therefore
$\gammah_i(MN)\leq \gammah_i(N)+\gammah_i(M)$, too.
\end{proof}

Now we apply the above results to the canonical modules of
$\Rees(I_t)$ and $A_t$. Assume for simplicity that $m\leq n$. Let
us try to find a product of minors $D$ such that $\ini(D)=\XX$ and
$\gamma_i(D)=\gammah_i(\XX)$ for all $i$. Since we have already
computed $\gammah_i(\XX)$ (see \ref{rhoXX}), we can determine the
shape of $D$, which turns out to be $1^2,2^2,\dots,\allowbreak
{(m-1)^2},\allowbreak m^{(n-m+1)}$. In other words, $D$ must be
the product of $2$ minors of size $1$, $2$ minors of size $2$,
$\dots$, $2$ minors of size $m-1$ and $n-m+1$ minors of size $m$.
Now it is not difficult to show that $D$ is uniquely determined,
the $1$-minors are $[m|1]$ and $[1|n]$, the $2$-minors are
$[m-1,m|1,2]$ and $[1,2|n-1,n]$ and so on.

\begin{theorem}
\label{canmodRA} Let $H$ be the $K$-subspace of $S$ whose
$K$-basis is the set of the elements of the form $\Delta T^k$
where $\Delta$ is a standard tableau with $\gamma_i(D\Delta)\geq
(k+1)(t+1-i)$ for all $=1\dots, t$.

Let $H_1$ be the $K$-subspace of $S$ whose $K$-basis is the set of
the elements of the form $\Delta T^k$ where $\Delta$ is a standard
tableau with $\gamma_2(D\Delta)\geq (k+1)(t-1)$ and
$\deg(D\Delta)=t(k+1)$. Set $J= DTH$ and $J_1=DTH_1$. Then we
have:
\begin{itemize}
\item[(a)] $J$ is an ideal of $\Rees(I_t)$. Furthermore $J$ is the
canonical module of $\Rees(I_t)$.
\item[(b)] $J_1$ is an ideal of $A_t$. Furthermore $J_1$ is the
canonical module of $A_t$.
\end{itemize}

\end{theorem}

\begin{proof}
That $J$ and $J_1$ are indeed ideals in the corresponding algebras
follows by the evaluation of shapes and by the description
\ref{algbygamma} of the algebras. Next we show that $\ini(J)$ and
$\ini(J_1)$ are the canonical modules of $\ini(\Rees(I_t))$ and
$\ini(A_t)$ respectively. It is enough to check that $\ini(J)$ is
exactly the ideal described in \ref{caninRe} and $\ini(J_1)$ is
the ideal described in \ref{caninAt}. Note that
$\ini(J)=\ini(DT)\ini(H)=\XX T\ini(H)$ and
$\ini(J_1)=\ini(DT)\ini(H_1)=\XX T\ini(H_1)$. Furthermore, by
virtue of \ref{rhoXX}, the canonical module of $\ini(\Rees(I_t))$
can be written as $\XX T G$ where $G$ is the space with basis the
set of the monomials $MT^k$ such that
$\gammah_i(M)+\gammah_i(\XX)\geq (k+1)(t+1-i)$ for all
$i=1,\dots,t$. Similarly, the canonical module of $\ini(A_t)$ can
be written as $\XX T G_1$ where $G_1$ is the space with basis the
set of the monomials $MT^k$ such that
$\gammah_2(M)+\gammah_2(\XX)\geq (k+1)(t-1)$ and
$\deg(M\XX)=t(k+1)$.

The spaces $H$ and $H_1$ are defined by the same inequalities
involving the $\gamma$ functions for bitableaux. As pointed out in
Remark \ref{Messina}(a), such vector spaces are in-KRS. This
implies $\ini(G)=H$, and similarly $\ini(G_1)=H_1$. As just
proved, $\ini(J)$ and $\ini(J_1)$ are the canonical modules of
$\ini(\Rees(I_t))$ and $\ini(A_t)$. Now the claim follows from
Theorem \ref{CanIni}. \end{proof}

\begin{remark}\label{class}\rm
(a) In \cite{BC4} we have translated the combinatorial description
of the canonical module into a divisorial one. If $t<\min(m,n)$,
the divisor class group of $\Rees(I_t)$ is free of rank $t$,
generated by the classes of prime ideals $P_1,\dots,P_t$ where a
$K$-basis of $P_i$ is given by all products $\Sigma T^k$, $\Sigma$
a standard bitableau with $\gamma_i(\Sigma)\ge k(t-i+1)+1$,
$i=1,\dots,t$. Then the canonical module of $\Rees(I_t)$ has
divisor class
$$
\sum_{i=1}^t \bigl(2-(m-i+1)(n-i+1)+t-i\bigr)
\cl(P_i)=\cl(I_t\Rees)+\sum_{i=1}^t (1-\height I_i)\cl(P_i)
$$

(b) If $t=\min(m,n)$, we may suppose that $t=m$. If even $t=m=n$,
then $I_t$ is a principal ideal, and $\Rees(I_t)$ is isomorphic to
a polynomial ring over $K$.

Let $t=m<n$. In this case a theorem of Herzog and Vasconcelos
\cite{HV} yields that the divisor class group of $\Rees(I_m)$ is
free of rank 1, generated by the extension $P$ of $I_m$ to
$\Rees(I_m)$. Moreover it implies that the canonical module has
class $(2-(n-m+1))P$.

(c) in Section \ref{SectAlgCM} we have seen that $A_t$ is factorial and,
hence, Gorenstein in the following cases: $t=1$, $t=\min(m,n)$
(the Grassmannian), and $t=m-1=n-1$.

(d) In all the cases different from those in (c), the ring $A_t$
is not factorial. Its divisor class group is free of rank $1$,
generated by the class of a single prime ideal $\qq$ that can be
chosen as $\qq=(f)S[T]\sect A_t$ where $f$ is a $(t+1)$-minor of
$X$. Then
$$
(mn-mt-nt)\cl(\qq).
$$
is the class of the canonical module; see \cite{BC4}.

(e) The expression for the class of the canonical module given in
(a) can be generalized to a larger class of Rees algebras; see
Bruns and Restuccia \cite{BR}. In particular one obtains results
for the algebras of minors of symmetric matrices of indeterminates
algebras generated by Pfaffians of alternating such matrices, and
algebras of minors of Hankel matrices. The latter case has been
treated by ``initial methods'' in \cite{BC4}.
\end{remark}

As a corollary we have

\begin{theorem}
\label{GorAt} The ring $A_t$ is Gorenstein if and only if one of
the following conditions is satisfied:
\begin{itemize}
\item[(a)] $t=1$; in this case $A_t=K[X]$.
\item[(b)] $t=\min(m,n)$; in this case $A_t$ is the coordinate ring
of a Grassmannian.
\item[(c)] $t=m-1$ and $m=n$; in this case $A_t$ is isomorphic to
a polynomial ring.
\item[(d)] $mn=t(m+n)$.
\end{itemize}
\end{theorem}

\begin{proof}
The cases (a), (b) and (c) are those discussed in Remark
\ref{class}(c). So we may assume that $1<t<\min(m,n)$ and $t\neq
m-1$ if $m=n$. Now Remark \ref{class}(d) completes the proof.

Since we have not discussed divisorial methods in detail, let us
indicate how to prove the theorem by combinatorial methods. In
view of \ref{GorIni} it makes no difference whether one works in
$A_t$ or $\ini(A_t)$. We choose $\ini(A_t)$.

Suppose that $mn=t(m+n)$. Then $\XX T^{m+n}$ not only belongs to
$\ini(A_t)$, but even to the ideal $\omega(\ini(A_t))$. Moreover,
$\gammah_2(\XX)-(t-1)(m+n)=1$. Using the ``linearity'' of $\XX$,
it is now easy to see that the ideal $\omega(\ini(A_t))$ is
generated by $\XX T^{m+n}$. So the canonical module is isomorphic
to $\ini(A_t)$, and $\ini(A_t)$ is Gorenstein.

In all the cases not covered by (a)--(d) one has to show that
$\ini(J_1)$ is not a principal ideal. We leave this to the reader
as an exercise. \end{proof}


\vspace{-0.1cm}
\begin{thebibliography}{99}

\bibitem{Abh} S.~S.~Abhyankar. {\em Enumerative combinatorics of Young
 tableaux.} M. Dekker 1988.

\bibitem{AK} S.~S.~Abhyankar and D.~M.~Kulkarni. {\em On Hilbertian
 ideals.} Linear Algebra Appl. {\bf 116} (1989), 53--79.

\bibitem{ABW} K. Akin, D. A. Buchsbaum, and J. Weyman. {\em Schur functors and
Schur complexes}. Adv. Math. {\bf 44} (1982), 207--278.

\bibitem{BF} J. Backelin and R. Fr\"oberg.
{\em Koszul algebras, Veronese subrings, and rings with linear
resolutions.} Rev. Roumaine Math. Pures Appl. {\bf 30} (1985),
85--97.

\bibitem{Ba} C. Bae\c tica. {\em Rees algebra of ideals generated by
pfaffians.} Commun. Algebra {\bf 26} (1998), 1769--1778.

\bibitem{BeZ} D. Bernstein and A. Zelevinsky.
{\em Combinatorics of maximal minors.} J. Algebr. Comb. {\bf 2}
(1993), 111--121.

\bibitem{Bj} A.~Bj\"orner. {\em Shellable and Cohen-Macaulay partially
ordered sets.} Trans. Amer. Math. Soc. {\bf 260} (1980), 159--183.

\bibitem{Br} W. Bruns. {\em Algebras defined by powers of
determinantal ideals}. J. Algebra {\bf 142} (1991), 150--163.

\bibitem{BC1} W. Bruns and A. Conca. {\em KRS
and powers of determinantal ideals.} Compositio Math. {\bf 111}
(1998), 111--122.

\bibitem{BC2} W. Bruns and A. Conca. {\em
The F-rationality of determinantal rings and their Rees rings.}
Mich. Math. J. {\bf 45} (1998), 291--299.

\bibitem{BC3} W. Bruns and A. Conca. {\em KRS and determinantal
rings}. In: J. Herzog, G. Restuccia (Eds.), Geometric and
combinatorial aspects of commutative algebra. Lecture Notes in
Pure and Applied Mathematics {\bf 217}. M. Dekker 2001, pp.\
67--87.

\bibitem{BC4} W. Bruns and A. Conca. {\em Algebras of minors.}
J. Algebra {\bf 246} (2001), 311--330.

\bibitem{BH} W. Bruns and J. Herzog. {\em Cohen-Macaulay rings},
Rev. Ed., Cambridge University Press 1998.

\bibitem{BK} W. Bruns and M. Kwieci\'nski. {\em Generic graph construction
ideals and Greene's theorem}. Math. Z. {\bf 233} (2000), 115--126.

\bibitem{BR} W. Bruns and G. Restuccia. {\em The canonical module
of a Rees algebra.} Preprint.

\bibitem{BST}  W. Bruns, A. Simis, and Ng\^o Vi\d{\^e}t Trung. {\em
Blow-up of straightening-closed ideals in ordinal Hodge algebras.}
Trans. Amer. Math. Soc. {\bf 326} (1991), 507--528.

\bibitem{BV} W.~Bruns and U.~Vetter. {\em Determinantal
rings.} Lect.~Notes Math. {\bf 1327}, Springer 1988.

\bibitem{CGG} L.~Caniglia, J.~A.~Guccione, and J.~J.~Guccione. {\em
Ideals of generic minors.} Commun.~Algebra {\bf 18} (1990),
2633--2640.

\bibitem{Co1} A.~Conca. {\em Gr\"obner bases of ideals of minors of
a symmetric matrix.}  J.~Algebra {\bf 166} (1994), 406--421.

\bibitem{Co2} A. Conca. {\em Ladder determinantal rings.} J. Pure Appl.
Algebra {\bf 98} (1995), 119--134.

\bibitem{Co3} A. Conca. {\em Gorenstein ladder determinantal rings.} J.
London Math. Soc. (2) {\bf 54} (1996), 453--474.

\bibitem{Co4} A.~Conca. {\em Gr\"obner bases of powers of ideals of maximal
minors.} J. Pure Appl. Alg. {\bf 121} (1997), 223--231.

\bibitem{Co5} A. Conca. {\em Straightening law and powers of determinantal
ideals of Hankel matrices.} Adv. Math. {\bf 138} (1998), 263--292.

\bibitem{CH1} A.~Conca and J.~Herzog. {\em On the Hilbert function
of determinantal rings and their canonical module.}
Proc.~Amer.~Math.~Soc. {\bf 122}, 677--681, (1994).

\bibitem{CH2} A. Conca and J. Herzog. {\em Ladder determinantal rings have
rational singularities.} Adv. Math. {\bf 132} (1997), 120--147.

\bibitem{CHV} A. Conca, J. Herzog, and G. Valla. {\em Sagbi bases
and application to blow-up algebras.} J. Reine Angew. Math. {\bf
474} (1996), 113--138.

\bibitem{DEP1} C.~De Concini, D.~Eisenbud, and C.~Procesi. {\em
Young diagrams and determinantal varieties.} Invent. math. {\bf
56} (1980), 129--165.

\bibitem{DEP2} C.~De Concini, D. Eisenbud, and C. Procesi. {\em Hodge algebras.}
Ast\'erisque {\bf 91}, Soc. Math. de France 1982.

\bibitem{DeN} E. De Negri. {\em ASL and Gr\"obner bases theory for Pfaffians and
monomial algebras.} Dissertation, Universit\"at Essen 1996.

\bibitem{DP} C.~De Concini, C.~Procesi. {\em A characteristic
free approach to invariant theory.} Adv. Math. {\bf 21}, (1976),
330--354.

\bibitem{DN} E. De Negri. {\em Some results on Hilbert series and
$a$-invariant of Pfaffian ideals.} Math. J. Toyama Univ. {\bf 24}
(2001), 93--106.

\bibitem{DRS} P. Doubilet, G.-C. Rota, and J. Stein. {\em On the
foundations of combinatorial theory: {IX}, {C}ombinatorial methods
in invariant theory.} Stud. Appl. Math. {\bf 53} (1974), 185--216.

\bibitem{Eis} { D.~Eisenbud.} {\em Commutative algebra with a view
towards algebraic geometry.} Springer 1995.

\bibitem{EH}  D. Eisenbud and C. Huneke. {\em Cohen-Macaulay Rees
algebras and their specialization.} J. Algebra {\bf 81} (1983),
202--224.

\bibitem{Fr} R. Fr\"oberg. {\em Determination of a class of Poincar\'e
series.} Math. Scand. {\bf 37} (1975), 29--39.

\bibitem{Fu} W. Fulton. {\em Young tableaux.} Cambridge University Press 1997.

\bibitem{GV} I.~M.~Gessel and G.~Viennot. {\em Binomial
determinants, paths, and hook length formulae.} Adv. Math. {\bf
58} (1985), 300--321.

\bibitem{GKr} S. R. Ghorpade and C. Krattenthaler. {\em The
Hilbert series of Pfaffian rings.} In: C. Bajaj, C. Christensen,
A. Sathaye, and G. Sundaram (Eds.), Algebra and Algebraic Geometry
with Applications. Springer, to appear.

\bibitem{Gho} S. Ghorpade. {\em Hilbert functions of ladder
determinantal varieties.} Discrete Math. {\bf 246} (2002),
131--175.

\bibitem{GL1} N. Gonciulea and V. Lakshmibai. {\em Schubert varieties,
toric varieties, and ladder determinantal varieties.} Ann. Inst.
Fourier {\bf 47} (1997), 1013--1064.

\bibitem{GL2} N. Gonciulea and V. Lakshmibai. {\em Singular loci of
ladder determinantal varieties and Schubert varieties.} J. Algebra
{\bf 229} (2000), 463--497.

\bibitem{Gre} C. Greene. {\em An extension of Schensted's theorem.}
Adv. Math. {\bf 14} (1974), 254--265.

\bibitem{HaT} J. Harris and L. Tu. {\em On symmetric and skew-symmetric
determinantal varieties.} Topology {\bf 23} (1984), 71--84.

\bibitem{HT} J. Herzog and Ng\^o
Vi\d{\^e}t Trung. {\em Gr\"obner bases and multiplicity of
determinantal and pfaffian ideals.} Adv. Math. {\bf 96} (1992),
1--37.

\bibitem{HV} J. Herzog and W. Vasconcelos. {\em On the divisor class
group of Rees-algebras.} J. Algebra {\bf 93} (1985), 182--188.

\bibitem{Ho} M.~Hochster. {\em Rings of invariants of tori,
Cohen-Macaulay rings generated by monomials, and polytopes.} Ann.
Math. {\bf 96} (1972), 318--337.

\bibitem{HE} M.~Hochster and
J.A.~Eagon. {\em Cohen-Macaulay rings, invariant theory, and the
generic perfection of determinantal loci.} Amer.~J.~Math. {\bf 93}
(1971), 1020--1058.

\bibitem{Hodge} W.V.D. Hodge. {\em Some enumerative results in the
theory of forms.} Proc. Camb. Philos. Soc. {\bf 39} (1943),
22--30.

\bibitem{HP} W.V.D. Hodge and D. Pedoe. {\em Methods of algebraic
geometry, Vol. II.} Cambridge University Press 1952.

\bibitem{Kn} D.~E.~Knuth. {\em Permutations, matrices, and
generalized Young tableaux.} Pacific J. Math. {\bf 34} (1970),
709--727.

\bibitem{Kn2} D.~E.~Knuth. {\em The art of computer programming,
Vol. 3.} Addison-Wesley 1975.

\bibitem{KM} A. Knutson and E. Miller. {\em Gr\"obner geometry of
Schubert polynomials.} Ann. Math., to appear.

\bibitem{Kr} C. Krattenthaler. {\em Counting nonintersecting
lattice paths with respect to weighted turns.} S\'emin. Lothar.
Comb. {\bf 34} (1995), B34i, 17 pp.

\bibitem{KrP} C. Krattenthaler and  M. Prohaska. {\em A remarkable
formula for counting nonintersecting lattice paths in a ladder
with respect to turns.} Trans. Am. Math. Soc. {\bf 351} (1999),
1015--1042.

\bibitem{KrR} C. Krattenthaler and M. Rubey. {\em A determinantal
formula for the Hilbert series of one-sided ladder determinantal
rings.} In: C. Bajaj, C. Christensen, A. Sathaye, and G. Sundaram
(Eds.), Algebra and Algebraic Geometry with Applications.
Springer, to appear.

\bibitem{KRo} M. Kreuzer and L. Robbiano. {\em Computational
commutative algebra 1.} Springer 2000.

\bibitem{Ku1} D. M. Kulkarni. {\em Hilbert polynomial of a certain
ladder-determinantal ideal.} J. Algebraic Combin. {\bf 2} (1993),
57--71.

\bibitem{Ku2} D. M. Kulkarni. {\em Counting of paths and
coefficients of Hilbert polynomial of a determinantal ring.}
 Discrete Math. {\bf 154} (1996), 141--151.

\bibitem{MaY} Y. Ma. {\em On the minors defined by a generic
matrix.} J. Symb. Comput. {\bf 18} (1994), 503--518.

\bibitem{Mul} S. B. Mulay. {\em Determinantal loci and the flag
variety.} Adv. Math. {\bf 74} (1989), 1--30.

\bibitem{Na}
H.~Narasimhan. {\em The irreducibility of ladder determinantal
varieties.} J.~Algebra {\bf 102} (1986), 162--185.

\bibitem{RS} L.~Robbiano and M.~Sweedler. {\em Subalgebra bases.}
In: W.~Bruns, A.~Simis (Eds.), Commutative Algebra, Lect.~Notes
Math. {\bf 1430}, Springer 1990, pp. 61--87.

\bibitem{Sa} B.~Sagan. {\em The symmetric group.} Second Ed.,
Springer 2000.

\bibitem{Sch} C.~Schensted. {\em Longest increasing and decreasing
subsequences.} Can.~J.~Math. {\bf 13} (1961), 179--191.

\bibitem{Sij} A. Schrijver. {\em Theory of linear and integer
programming.} Wiley-Interscience 1986.

\bibitem{Stu1} B.~Sturmfels. {\em Gr\"obner bases and Stanley
decompositions of determinantal rings.} Math.~Z. {\bf 205} (1990),
137--144.

\bibitem{Stu2} B.~Sturmfels.
{\em Gr\"obner bases and convex polytopes.}
 Amer. Math. Soc. 1996.

\bibitem{Sta2} R. P. Stanley. {\em Combinatorics and Commutative
Algebra.} Second Ed., Birkh\"auser 1996.

\bibitem{Sta3}  R. P. Stanley. {\em Enumerative combinatorics, Vol. 2.}
Cambridge University Press 1999.

\bibitem{Tr} Ng\^o Vi\d{\^e}t Trung. {\em On the symbolic powers of
determinantal ideals.} J. Algebra {\bf 58} (1979), 361--369.

\bibitem{Va} W. V. Vasconcelos. {\em Computational methods in
commutative algebra and algebraic geometry.} Springer 1998.


\end{thebibliography}
\end{document}